\documentclass[11pt,british,reqno]{amsart}

\usepackage[T1]{fontenc}

\usepackage{babel,eucal,url,amssymb,stmaryrd,booktabs,enumerate,amscd,paralist,cite}

\theoremstyle{plain}
\newtheorem{lemma}{Lemma}[section]

\newtheorem{proposition}[lemma]{Proposition}
\newtheorem{corollary}[lemma]{Corollary}
\newtheorem{theorem}[lemma]{Theorem}
\newtheorem{remark}[lemma]{Remark}

\newtheorem*{ack}{Acknowledgements}

\DeclareMathOperator{\Ric}{Ric}

\newcommand{\Lie}[1]{\operatorname{\textsl{#1}}}
\newcommand{\lie}[1]{\operatorname{\mathfrak{#1}}}

\newcommand{\SO}{\Lie{SO}}

\newcommand{\Un}{\Lie{U}}
\newcommand{\un}{\lie{u}}
\newcommand{\Gtwo}{\ifmmode{{\rm G}_2}\else{${\rm G}_2$}\fi}

\newcommand{\sym}[1]{{S^{#1}\mskip-2mu}}

 \newcommand{\cyclic}{\mathop{\kern0.9ex{{+}\kern-2.2ex\raise-.28ex\hbox{\Large\hbox
 {$\circlearrowright$}}}}}

\newcommand{\real}[1]{\left\llbracket #1 \right\rrbracket}

\newcommand{\talt}{\tilde{\mathbf a}}

\newcommand{\Ricac}{\Ric^{\textup{ac}}}

\def\sideremark#1{\ifvmode\leavevmode\fi\vadjust{\vbox to0pt{\vss
 \hbox to 0pt{\hskip\hsize\hskip1em
 \vbox{\hsize2.5cm\tiny\raggedright\pretolerance10000
 \noindent #1\hfill}\hss}\vbox to8pt{\vfil}\vss}}}%

\newfont{\eusm}{eusm10 scaled \magstep1}
\newfont{\eusmiii}{eusm10 scaled \magstep3}

\newcommand{\comp}{\makebox[7pt]{\raisebox{1.5pt}{\tiny $\circ$}}}

\title{Harmonic almost contact structures via the intrinsic torsion}

\author{J.~C.~Gonz{\'a}lez-D{\'a}vila}
\address[J.~C.~Gonz{\'a}lez-D{\'a}vila]{Department of Fundamental Mathematics\\
  University of La Laguna\\ 38200 La Laguna, Tenerife, Spain}
\email{jcgonza@ull.es}

\author{F.~Mart\'\i n~Cabrera}
\address[F.~Mart\'\i n~Cabrera]{Department of Fundamental Mathematics\\ CSIC Associated Unity\\
  University of La Laguna\\ 38200 La Laguna, Tenerife, Spain}
\email{fmartin@ull.es}

\date{\today}
\begin{document}

\maketitle

\begin{abstract}{\indent}
We go further on the study of harmonicity for almost contact metric
structures already initiated by Vergara-Díaz and Wood. By using the
intrinsic torsion, we characterise harmonic almost contact metric
structures in several equivalent ways and show conditions relating
harmonicity and classes of almost contact metric structures.
Additionally, we study the harmonicity of such structures  as a map
into the quotient bundle of the oriented orthonormal frames by the
action of the structural group $\Lie{U}(n) \times 1$. Finally, by
using a Bochner type formula proved by Bör and Hernández Lamoneda,
we display some examples which give the absolute minimum for the
energy.

\vspace{2mm}

  \vspace{2mm}

\end{abstract}

\tableofcontents

\section{Introduction}{\indent} \setcounter{equation}{0}
For an oriented Riemannian manifold $M$ of dimension $n$, given
 a Lie subgroup $G$ of $\SO(n)$, $M$ is said to  be equipped with a
\emph{$G$-structure}, if there exists a subbundle $\mathcal{G}(M)$
of the oriented orthonormal frame bundle $\mathcal{SO}(M)$ with
structural group $G$. For a fixed $G$, a natural question arises,
'which are the best $G$-structures on M?'. An approach to answer
this question is based on the notion of the energy  of a
$G$-structure which is a particular case of the energy of a map
between Riemannian manifolds. Such a functional has been widely
studied by diverse authors \cite{EeLe1,EeLe2,Ur}. The corresponding
critical points are called \emph{harmonic maps} and have been
characterised by Eells and Sampson \cite{EeSa}.

For principal $G$-bundles $Q \to M$ over a Riemannian manifold,
Wood in \cite{Wood2} considers global sections $\sigma \,: \; M
\to Q/H$ of the quotient bundle $\pi : Q/H \to M$, where $H$ is a
Lie subgroup of $G$ such that $G/H$ is reductive. Such global
sections are in one-to-one correspondence with the $H$-reductions
of the $G$-bundle $Q \to M$. Likewise, a connection on $Q \to M$
and a $G$-invariant metric on $G/H$ are fixed. Thus, $Q/H$ can be
equipped in a natural way with a metric defined by using the
metrics on $M$ and $G/H$.  In such conditions, Wood regards
harmonic sections as generalisations of harmonic maps from $M$
into $Q/H$, deriving the corresponding characterising conditions,
called {\it harmonic sections equations}.

The situation described in the previous paragraph arises when the
Riemannian manifold $M$ is equipped with some additional geometric
structure, viewed as reduction of the structure group of the tangent
bundle. Thus, in \cite{GDMC} it is  considered the particular
situation for $G$-structures defined on an oriented Riemannian
manifold $M$ of dimension $n$, where $\Lie{G}$ is  a closed and
connected subgroup of $\SO(n)$. Since the existence of a
$G$-structure on $M$ is equivalent to the existence of a global
section $\sigma  : M \to \mathcal{SO}(M)/G$ of the quotient bundle,
the energy  of a $G$-structure is defined as the energy of the
corresponding map $\sigma$. If $\xi^G$ denotes the intrinsic torsion
of the $G$-structure, such an energy functional is essentially
determined by the \emph{total bending} given by $
 B(\sigma)=\tfrac12 \int_{M}\| \xi^G\|^{2} \, dv.
$ As a consequence, the notion of \emph{harmonic $G$-structure},
introduced by Wood in \cite{Wood2}, is given in terms of the
intrinsic torsion. This analysis makes possible to go further in
the study of relations between harmonicity and classes of
$G$-structures. Thus, the  study of harmonic almost Hermitian
structures initiated in \cite{Wood1,Wood2} is completed in
\cite{GDMC} with additional results.

Our purpose in the present work is going on the study of
harmonicity for almost contact metric structures initiated by
Vergara-Díaz and Wood in  \cite{VergaraWood}.  Almost contact
metric structures can be seen as $\Lie{U}(n)$-structures defined
on Riemannian manifolds of dimension $2n+1$.

After characterising harmonic almost contact metric structures in
several equivalent ways (Theorem \ref{characharmalmcontact}), we
show conditions relating harmonicity and   Chinea and
González-Dávila's  classes \cite{ChineaGonzalezDavila} of almost
contact metric structures. This is exposed in Theorem
\ref{harmclasscontact}. Finally, the harmonicity  of an almost
contact metric structure as a map into $\mathcal{SO}(M)/(\Lie{U}(n)
\times 1)$ is studied in Theorem \ref{mapharmcontact}.

As a relevant remark, we point out the  rôle played by the
identities given in Lemmas \ref{xxx} and \ref{ec1bisbisbisxx}.
They are consequences of the trivial identity $d^2 F=0$, where $F$
is the fundamental two-form of the almost contact metric structure
(see Section \ref{sect:almcont}). Such identities are deduced by
firstly expressing $d^2 F$ in terms of the intrinsic torsion and
the minimal connection, and then extracting certain
$\Lie{U}(n)$-components. An analog identity for almost Hermitian
structures was deduced  in \cite{FMCAS}. In the proofs of Theorems
\ref{harmclasscontact} and \ref{mapharmcontact}, the use of these
identities beside the harmonicity criteria is fundamental.

In the final part of this work, we  apply the simple and elegant
argument used by Bör et al. in \cite{BHLS2} to almost Hermitian
structures,  here  to almost contact metric structures. Thus, by
making use of a Bochner type formula proved in \cite{BHLS1} (see
the equation \eqref{mainBoch} below), we show that in some
situations there is an absolute minimum.  This and other results
are applied in the last Section, where we display several
examples. In particular, some of them give the absolute minimum
for the energy.

\begin{ack} {\rm  The authors are supported by  grants from MEC (Spain),
projects  MTM2007-65852 and MTM2007-66375.}
\end{ack}

\section{Preliminaries}{\indent} \setcounter{equation}{0}
On an $n$-dimensional oriented Riemannian manifold $(M ,\langle
\cdot , \cdot \rangle),$ we consider the principal $SO(n)$-bundle
$\pi_{\Lie{SO}(n)} : \mathcal{SO}(M) \to M$ of the oriented
orthonormal frames with respect to the metric $\langle \cdot,
\cdot \rangle.$ Given a closed and connected subgroup $G$ of
$\SO(n),$ a $G$-structure on $(M ,\langle \cdot , \cdot \rangle)$
is a reduction $\mathcal{G}(M)\subset \mathcal{SO}(M)$ to $G.$ In
the present Section we briefly recall some notions relative to
$G$-structures (see \cite{GDMC,Wood2} for  more details).

Let ${\mathcal S\mathcal O}(M)/G$ be the orbit space under the
action of $G$ on ${\mathcal S\mathcal O}(M)$ on the right as
subgroup of $SO(n).$ Then the $G$-orbit map $\pi_{G}:{\mathcal
S\mathcal O}(M)\to {\mathcal S\mathcal O}(M)/G$ is a principal
$G$-bundle and we have $\pi_{SO(n)} = \pi\comp \pi_{G},$ where
$\pi: {\mathcal S\mathcal O}(M)/G\to M$ is a fibre bundle with
fibre $SO(n)/G,$ which is naturally isomorphic to the associated
bundle ${\mathcal S\mathcal O}(M)\times_{SO(n)}SO(n)/G.$ The map
$\sigma:M\to {\mathcal S\mathcal O}(M)/G$ given by $\sigma(m) =
\pi_{G}(p),$ for all $p\in {\mathcal G}(M)$ with $\pi_{SO(n)}(p) =
m,$ is a smooth section and  determines a one-to-one
correspondence between the totally of $G$-structures and the
manifold $\Gamma^{\infty}({\mathcal S\mathcal O}(M)/G)$ of all
global sections of ${\mathcal S\mathcal O}(M)/G.$ In what sequel,
we shall also denote by $\sigma$ the $G$-structure determined by
the section $\sigma$.

The presence of a reduced subbundle $\mathcal{G}(M)$ yields to
express the bundle of endomorphisms $\mbox{End}(\mbox{T} M)$ on the
fibers in the tangent bundle $\mbox{T}M$ as the associated vector
bundle $\mathcal{G}(M) \times_{G} \mbox{End}(\mathbb{R}^n)$. We
restrict our attention on the subbundle $\lie{so}(M)$ of $\mbox{End}
(\mbox{T} M )$ of skew-symmetric endomorphisms $\varphi_m$, for all
$m \in M$, i.e. $\langle \varphi_m X , Y \rangle= -\langle \varphi_m
Y , X \rangle$. Note that this subbundle $\lie{so}(M)$ is expressed
as $ \lie{so}(M) = \mathcal{SO}(M) \times_{\Lie{SO}(n)} \lie{so}(n)
= \mathcal{G}(M) \times_{G} \lie{so}(n) $. Furthermore, because
$\lie{so}(n)$ is decomposed into the $G$-modules $\lie{g},$ the Lie
algebra of $G,$ and the orthogonal complement $\lie{m}$ on
$\lie{so}(n)$ with respect to the inner product $\langle \cdot ,
\cdot \rangle$ given by $\langle X,Y\rangle = -\mbox{trace}\; XY,$
the bundle $\lie{so}(M)$ is also decomposed into $\lie{so}(M) =
\lie{g}_{\sigma} \oplus \lie{m}_{\sigma}$, where $\lie{g}_{\sigma} =
\mathcal{G}(M) \times_G \lie{g}$ and $\lie{m}_{\sigma} =
\mathcal{G}(M) \times_G \lie{m}$.

Under the conditions above fixed, if $M$ is equipped with a
$G$-structure, then there exists a $G$-connection
$\widetilde{\nabla}$ defined on $M$. Doing the difference
$\widetilde{\xi}_X = \widetilde{\nabla}_X - \nabla_X$, where
$\nabla_X$ is the Levi-Civita connection of $\langle \cdot , \cdot
\rangle$, a tensor $\widetilde{\xi}_X \in \lie{so}(M)$ is
obtained.  Decomposing $\widetilde{\xi}_X = ( \widetilde{\xi}_X
)_{\lie{g}_{\sigma}} + ( \widetilde{\xi}_X )_{\lie{m}_{\sigma}}$,
$( \widetilde{\xi}_X )_{\lie{g}_{\sigma}} \in \lie{g}_{\sigma}$
and $( \widetilde{\xi}_X )_{\lie{m}_{\sigma}} \in
\lie{m}_{\sigma}$, a new $G$-connection $\nabla^G$, defined by
$\nabla^G_X = \widetilde{\nabla}_X - (\tilde{\xi}_X
)_{\lie{g}_{\sigma}}$, can be considered. Because the difference
between two $G$-connections must be in $\lie{g}_{\sigma}$,
$\nabla^G$ is the unique $G$-connection on $M$ such that its
torsion satisfies the condition $\xi^G_X = ( \widetilde{\xi}_X
)_{\lie{m}_{\sigma}} = \nabla^{G}_X - \nabla_X$ is in
$\lie{m}_{\sigma}$. $\nabla^G$ is called the {\it minimal
connection} and $\xi^G$ is referred to as the {\it intrinsic
torsion} of the $G$-structure $\sigma$
\cite{CleytonSwann:torsion}. If $\xi^G=0$, the $G$-structure is
said  to be {\it integrable}. In such a case, the Riemannian
holonomy group of $M$ is contained in $G$.

\begin{remark} \label{clasif:struc}
{\rm  A natural way of classifying $G$-structures arises by
decomposing the space $\mathcal W = \mbox{T}^* M \otimes
\lie{m}_{\sigma}$ of possible intrinsic torsion into irreducible
$G$-modules. This was initiated by Gray and Hervella
\cite{Gray-H:16} for almost Hermitian structures. In this
particular case, $G = \Lie{U}(n)$ and the space $\mathcal W$ is
decomposed into four irreducible $\Lie{U}(n)$-modules. Therefore,
$2^4=16$ classes of almost Hermitian structures were obtained. }
\end{remark}

Along  the present paper, we will consider the natural extension
of the metric $\langle \cdot , \cdot \rangle$ to $(r,s)$-tensors
on $M$. Such an extension is defined by
\begin{equation} \label{extendedmetric}
\langle \Psi,\Phi \rangle = \Psi^{i_{1}\dots i_{r}}_{j_{1}\dots
j_{s}} \Phi^{i_{1}\dots
 i_{r}}_{j_{1}\dots j_{s}},
  \end{equation}
where the summation convention is used and  $\Psi^{i_{1}\dots
i_{r}}_{j_{1}\dots j_{s}}$ and $\Phi^{i_{1}\dots i_{r}}_{j_{1}\dots
j_{s}}$ are the components of the $(r,s)$-tensors  $\Psi,\Phi\in
\mbox{T}_{s\, m}^{r} M$, with respect to an orthonormal frame over
$m \in M$. Likewise, we will make reiterated use of the {\it musical
isomorphisms} $\flat : \mbox{\rm T}M \to \mbox{\rm T}^* M$ and
$\sharp : \mbox{\rm T}^* M \to \mbox{\rm T} M$, induced by the
metric $\langle \cdot , \cdot \rangle$ on $M$, defined respectively
by $X^{\flat} = \langle X , \cdot \rangle$ and $\langle
\theta^{\sharp} , \cdot \rangle = \theta $. \vspace{2mm}

If $\omega$ is  the connection one-form associated to $\nabla$,
then $\mbox{T} \mathcal{SO}(M) = \ker \pi_{\Lie{SO}(n)\ast} \oplus
\ker \omega$. Now considering  the projection $\pi_{G} \colon
\mathcal{SO}(M) \to \mathcal{SO}(M)/G$,   the tangent bundle of
$\mathcal{SO}(M)/G$ is decomposed into $\mbox{T} \mathcal{SO}(M)/G
= \mathcal{V} \oplus \mathcal{H}$, where $\mathcal{V} =
\pi_{G\ast} (\ker \pi_{\Lie{SO}(n)\ast} )$ and $\mathcal{H} =
\pi_{G\ast} (\ker \omega )$.
 For the projection
$\pi \, : \, \mathcal{SO}(M)/G \to M$, $\pi(pG) =
\pi_{\Lie{SO}(n)} (p)$,  the {\it vertical} and {\it horizontal}
distributions $\mathcal{V}$ and $\mathcal{H}$ are such that
$\pi_{\ast} \mathcal{V} =0$ and $\pi_{\ast} \mathcal{H} =
\mbox{\rm T}M$.

 Moreover, it is  considered
the bundle $\pi^* \lie{so}(M)$ on $\mathcal{SO}(M)/G$  consisting
of those pairs $(pG, \breve{\varphi}_m)$, where $\pi(pG)=m$ and
 $\breve{\varphi}_m \in \lie{so}(M)_m$.
  Alternatively, $\pi^* \lie{so}(M)$ is described as the bundle
$\pi^* \lie{so}(M)  = \mathcal{SO}(M) \times_G \lie{so}(n)=
\lie{g}_{\mathcal{SO}(M)} \oplus \lie{m}_{\mathcal{SO}(M)}$, where
$\lie{g}_{\mathcal{SO}(M)}= \mathcal{SO}(M) \times_G \lie{g} $ and
 $\lie{m}_{\mathcal{SO}(M)}= \mathcal{SO}(M) \times_G \lie{m}$.
A metric on each fiber of $\pi^* \lie{so}(M)$ is defined by
$$
\langle (pG , \breve{\varphi}_m) , (pG , \breve{\psi}_m) \rangle =
\langle \breve{\varphi}_m , \breve{\psi}_m \rangle,
$$
where $\langle \cdot , \cdot \rangle$ in the right side is the
extension to $(1,1)$-tensors of the metric on $M$ given by
\eqref{extendedmetric}. With respect to this metric,  the
decomposition $\pi^* \lie{so}(M) = \lie{g}_{\mathcal{SO}(M)}
\oplus \lie{m}_{\mathcal{SO}(M)}$ is orthogonal.

There is a canonical isomorphism between $\mathcal{V}$ and the
bundle $\lie{m}_{\mathcal{SO}(M) }$. In fact,  elements in
$\lie{m}_{\mathcal{SO}(M)}$ can be seen as pairs $(pG ,
\breve{\varphi}_m)$ such that if $\breve{\varphi}_m$ is expressed
with respect to  $p$, then  it is obtained a matrix $( a_{ji}) \in
\lie{m}$.  For all $a \in \lie{m}$, we have the fundamental vector
field $a^*$ on $\mathcal{SO}(M)$ given by
$$
a^*_p = \frac{d}{dt}_{|t=0} p . \exp t a \in \ker
\pi_{\Lie{SO}(n)*p} \subseteq \mbox{T}_p \mathcal{SO}(M).
$$
Any vector in $\mathcal{V}_{pG}$ is given by $\pi_{G*p} (a^*_p)$,
for some $a =(a_{ji}) \in \lie{m}$. The isomorphism $\phi_{|
\mathcal{V}_{pG}} \colon \mathcal{V}_{pG} \to \left(
\lie{m}_{\mathcal{SO}(M)} \right)_{pG}$  is defined by
$$
\phi_{| \mathcal{V}_{pG}} ( \pi_{G*p} (a^*_p)) = (pG, a_{ji} \,
p(u_i)^{\flat} \otimes p(u_j)).
$$
Next it is extended the map $\phi_{|\mathcal{V}} : \mathcal{V} \to
\lie{m}_{\mathcal{SO}(M)}$ to $\phi : \mbox{T} \,
\mathcal{SO}(M)/G \to \lie{m}_{\mathcal{SO}(M)}$ by saying that
$\phi (A) =0$, for all $A \in \mathcal{H}$,  and $\phi (V) =
\phi_{|\mathcal{V}}(V)$, for all $V \in \mathcal{V}$. This is used
to define a metric $\langle \cdot , \cdot
\rangle_{\mathcal{SO}(M)/G}$ on $\mathcal{SO}(M)/G$ by
\begin{equation} \label{metricquo}
\langle A , B \rangle_{\mathcal{SO}(M)/G} = \langle \pi_{\ast} A ,
\pi_{\ast} B \rangle + \langle \phi (A) , \phi (B) \rangle.
\end{equation}
For this metric, the projection $\pi \, : \, \mathcal{SO}(M)/G \to
M$ is a Riemannian submersion with totally geodesic fibres (see
\cite{Vilms} and \cite[page 249]{Besse:Einstein}).

Now, it is considered the set of all possible $G$-structures on a
closed and oriented Riemannian manifold $M$ which are compatible
with the metric $\langle \cdot , \cdot \rangle$. Such a set is
identified with the manifold $\Gamma^{\infty}(\mathcal{SO}(M)/G)$ of
all possible global sections $\sigma \, : \, M \to
\mathcal{SO}(M)/G$. With   respect to the metrics $\langle \cdot ,
\cdot \rangle$ and  $ \langle \cdot , \cdot
\rangle_{\mathcal{SO}(M)/G}$, the {\it energy} of $\sigma$ is the
integral
\begin{equation}\label{1}
 {\mathcal E}(\sigma)=\frac{\textstyle 1}{\textstyle2}\int_{M}\|\sigma_{\ast}\|^{2}dv,
\end{equation}
where $\|\sigma_{\ast}\|^{2}$ is the norm of the differential
$\sigma_{\ast}$ of $\sigma$ and $dv$ denotes the volume form on
$(M,\langle \cdot , \cdot \rangle)$. On the domain of a local
orthonormal frame field $\{e_1, \dots , e_n \}$ on $M$,
$\|\sigma_{*}\|^{2}$ can be locally expressed as
$\|\sigma_{*}\|^{2} = \langle
\sigma_{*}e_{i},\sigma_{*}e_{i}\rangle_{\mathcal{SO}(M)/G}$.
Furthermore,  using \eqref{metricquo}, from \eqref{1}  it is
obtained that the energy $\mathcal{E}(\sigma)$ of $\sigma$ is
given by
\[
{\mathcal E}(\sigma) = \frac{\textstyle n}{\textstyle 2} {\rm
Vol}(M) + \frac{\textstyle 1}{\textstyle 2}\int_{M}\|  \phi \,
\sigma_*\|^{2}dv.
\]
The relevant part of this formula $B(\sigma)=\frac{\textstyle
1}{\textstyle 2}\int_{M}\| \phi \, \sigma_*\|^{2}dv$  is called
the {\it total bending} of the $G$-structure $\sigma$. In
\cite{GDMC}, it was shown that $\phi \, \sigma_{\ast} = - \xi^G$.
Therefore,
$$
B(\sigma)=\frac{\textstyle 1}{\textstyle 2}\int_{M}\| \xi^G\|^{2}
\, dv.
$$

 To study critical points
  of the functional $\mathcal{E}$ on
$\Gamma^{\infty}(\mathcal{SO}(M)/G)$,   smooth variations
$\sigma_t \in \Gamma^{\infty} (\mathcal{SO}(M)/G)$ of
$\sigma=\sigma_0$ are considered. The corresponding {\it variation
fields} $m
 \to \varphi(m) = \frac{d}{dt}_{|t=0} \sigma_t(m)$  are sections of the
induced bundle $\sigma^* \mathcal{V}$ on $M$.
  Furthermore, by using $\phi$, we will have $\phi
  \mbox{pr}_2^{\sigma} \sigma^* \mathcal{V}   \cong \sigma^* \lie{m}_{\mathcal{SO}(M)}
\cong \lie{m}_{\sigma}$. Thus, the tangent space
$\mbox{T}_{\sigma} \Gamma^{\infty} (\mathcal{SO}(M)/G)$ is firstly
identified with the space $\Gamma^{\infty} (\sigma^* \mathcal V)$
of global sections of $\sigma^* \mathcal V$ \cite{Ur}. A second
identification is $\Gamma^{\infty} (\sigma^* \mathcal V ) \cong
\Gamma^{\infty} (\lie{m}_{\sigma})$ as global sections of
$\lie{m}_{\sigma}$.

In following results, we will consider the coderivative $d^*
\xi^G$ of the intrinsic torsion $\xi^G$, which is given by
$$
d^* \xi^G_m = - (\nabla_{e_i} \xi^G)_{e_i}  = - (\nabla^{G}_{e_i}
\xi^G)_{e_i} - \xi^G_{\xi^G_{e_i} e_i} \in \lie{m}_{\sigma \, m} ,
$$
where $\{ e_1, \dots , e_n \}$ is any orthonormal frame on $m\in
M$. Therefore,  $d^* \xi^G$ is a global section of
$\lie{m}_{\sigma}$.

\begin{theorem}[\cite{GDMC}]
If  $\,G$ is a closed and connected subgroup of $\Lie{SO}(n)$,
$(M,\langle \cdot , \cdot\rangle)$ a closed and oriented
Riemannian manifold and $\sigma$ a global section of
$\mathcal{SO}(M)/G$, then:
\begin{enumerate}
\item[{\rm (i)}] $($The first variation formula$)$. For the energy
functional $\mathcal{E}: \Gamma^{\infty}(\mathcal{SO}(M)/G) \to
\mathbb{R}$ and for all $\varphi \in \Gamma^{\infty}
(\lie{m}_{\sigma}) \cong \mbox{\rm T}_{\sigma} \Gamma^{\infty}
(\mathcal{SO}(M)/G)$, we have
$$
d \mathcal{E}_{\sigma} (\varphi) = - \int_M \langle \xi^G , \nabla
\varphi \rangle dv =  - \int_M \langle d^* \xi^G , \varphi \rangle
dv.
$$
\item[{\rm (ii)}] $($The second variation formula$)$.  The Hessian
form $({\rm Hess}\;{\mathcal E})_{\sigma}$ on $
\Gamma^{\infty}(\lie{m}_{\sigma} )$ is given by
\begin{eqnarray*}
({\rm Hess}\; {\mathcal E})_{\sigma}\varphi & = & \int_{M} \left(
\|\nabla  \varphi \|^{2}   - \tfrac12 \| [ \xi^G ,
\varphi]_{\lie{m}_{\sigma}}
  \|^2   +   \langle \nabla
\varphi , 2[ \xi^G , \varphi] -  [ \xi^G ,
\varphi]_{\lie{m}_{\sigma}} \rangle  \right) dv.
\end{eqnarray*}
\end{enumerate}
\end{theorem}

As a consequence of this Theorem the following notion is
introduced: for general Riemannian manifolds $(M,\langle \cdot
,\cdot \rangle)$ not necessarily closed and oriented,
 a $G$-structure $\sigma$ is said to be  {\em harmonic}, if it satisfies $d^* \xi^G=0$ or, equivalently,
 $(\nabla^G_{e_i} \xi)_{e_i} = -  \xi^G_{\xi^G_{e_i} e_i}$.

Given a $G$-structure $\sigma$ on a closed Riemannian manifold
$(M,\langle\cdot ,\cdot \rangle),$ the map $\sigma:(M, \langle
\cdot, \cdot \rangle) \mapsto (\mathcal{SO}(M)/G,
\langle\cdot,\cdot\rangle_{\mathcal{SO}(M)/G})$ is harmonic, i.e.
$\sigma$ is a critical point for the energy functional on
$\mathcal{C}^{\infty}(M,\mathcal{SO}(M)/G)$ if and only if its {\it
tension field} $\tau(\sigma)= \left( \nabla_{e_i}
\sigma_*\right)(e_i)$ vanishes \cite{Ur}. Here, $\nabla \sigma_{*}$
is defined by $\left( \nabla_X \sigma_* \right)(Y) =
\nabla^q_{\sigma_* X} \sigma_* Y - \sigma_* (\nabla_X Y),$ where
$\nabla^q$ denotes the induced connection by the Levi-Civita
connection $\nabla^q$ of the metric in $\mathcal{SO}(M)/G$.
According with \cite{GDMC,Wood2}, harmonic sections $\sigma$ are
characterised by the vanishing of the vertical component of
$\tau(\sigma)$ and the horizontal component of $\tau(\sigma)$ is
determined by the horizontal lift of the vector field metrically
equivalent to the one-form $\nu_{\sigma},$ defined by
\begin{equation}
 \label{harmmap1} \nu_{\sigma}(X) = \langle \xi^G_{e_i} ,
R_{e_i,X}\rangle.
\end{equation}
Hence, the map $\sigma:(M, \langle \cdot, \cdot \rangle) \mapsto
(\mathcal{SO}(M)/G, \langle\cdot,\cdot\rangle_{\mathcal{SO}(M)/G})$
turns out to be a harmonic map if and only if $\sigma$ is a harmonic
$G$-structure such that $\nu_{\sigma}$ vanishes for all vectors.
\vspace{1mm}

 If $R_{X,Y\lie{m}_{\sigma}}=0$, the $G$-structure $\sigma$ is referred to as a {\it flat
$G$-structure}. The intrinsic torsion of a flat $G$-structure has
not contributions in the $G$-components of $R$ orthogonal to
$\sym{2} \lie{g}_{\sigma}$. Thus, $R$ is in the space of algebraic
curvature tensors for manifolds with an integrable $G$-structure.

Relevant types of diverse $G$-structures are characterised by saying
that their intrinsic torsion $\xi^G$ is metrically equivalent to a
skew-symmetric three-form, that is, $\xi^G_X Y = - \xi^G_Y X$. Now
we will recall some facts satisfied by  such $G$-structures.

\begin{proposition}[\cite{GDMC}] \label{pro:skew}  For   a $G$-structure $\sigma$  such
that $\xi^G_X Y = - \xi^G_{Y} X$, we have:
\begin{enumerate}
  \item[{\rm (i)}]
 If $[\xi^G_X  , \xi^G_Y] \subseteq \lie{g}_\sigma$, for all $X,Y
 \in \mathfrak X (M)$,
  then
    $\langle R_{X,Y  \lie{m}_{\sigma}} X , Y \rangle  =  2\langle
\xi^G_{X} Y , \xi^G_X Y \rangle$. Therefore, $\sigma$ is
integrable if and only if  $\sigma$ is flat.
 \item[{\rm (ii)}] If $\sigma$ is a harmonic $G$-structure, then
$\sigma$ is also a harmonic map.
\end{enumerate}
\end{proposition}

In Section \ref{harmonicalmostcontact}, we will  study harmonicity
of almost contact metric structures. Such structures are examples of
$G$-structures defined by means of one or several $(r,s)$-tensor
fields $\Psi$ which are stabilised under the action of $G$, i.e. $g
\cdot \Psi= \Psi$, for all $g \in G$. Moreover, it will be possible
characterise the harmonicity of such $G$-structures by conditions
given in terms of those tensors $\Psi$. The \emph{connection
Laplacian} (or {\it rough Laplacian}) $\nabla^* \nabla \Psi$ will
play a relevant rôle  in such conditions. We recall that
$$
\nabla^* \nabla \Psi = -  \left( \nabla^2 \Psi \right)_{e_i,e_i},
$$
where $\{ e_1, \dots , e_n \}$ is an orthonormal  frame field and
$(\nabla^2\Psi)_{X,Y} = \nabla_X (\nabla_Y \Psi) -
\nabla_{\nabla_XY}\Psi$.  If a Riemannian manifold $(M , \langle
\cdot , \cdot \rangle )$ of dimension $n$ is equipped by a
$G$-structure, where $G \subseteq \Lie{SO}(n)$ is closed and
connected, and
 $\Psi$ is a $(r,s)$-tensor field on $M$ which is stabilised under the action of $G$, then
\begin{equation}\label{lapstaten}
\nabla^* \nabla \Psi =  \left( \nabla^{G}_{e_i}
\xi^{G}\right)_{e_i} \Psi + \xi^{G}_{\xi^{G}_{e_i}e_i} \Psi -
\xi^{G}_{e_i} (\xi^{G}_{e_i} \Psi).
\end{equation}
Moreover, if the $G$-structure is harmonic,  $ \nabla^* \nabla \Psi
= -\xi^{G}_{e_i} (\xi^{G}_{e_i} \Psi). $

\section{Almost contact metric structures}{\indent}
\label{sect:almcont} \setcounter{equation}{0}
 An almost contact
metric manifold is a $2n+1$-dimensional Riemannian manifold
$(M,\langle \cdot, \cdot \rangle)$  equipped with a $(1,1)$-tensor
field $\varphi$ and a unit vector field $\zeta$,  called the {\it
characteristic vector field} of the structure,  such that
\[
\varphi^{2} = -I + \eta\otimes \zeta,\;\;\;\;\; \langle \varphi
X,\varphi Y\rangle = \langle X,Y\rangle - \eta(X)\eta(Y),
\]
where $\eta = \zeta^{\flat}$. Associated with the almost contact
metric structure the two-form $F= \langle \cdot, \varphi
\cdot\rangle$, called the {\em fundamental two-form}, is usually
considered. Using $F$ and $\eta$, $M$ can be oriented fixing a
constant multiple of $F^n \wedge \eta = F \wedge
\stackrel{(n)}{\dots}\wedge F \wedge \eta$ as volume form. Likewise,
the presence of an almost contact metric structure is equivalent to
say that $M$ is equipped with a $\Lie{U}(n) \times 1$-structure. It
is well known that $\Lie{U}(n)\times 1$ is a closed and connected
subgroup of $\SO(2n+1)$ and $\SO(2n+1)/ (\Lie{U}(n) \times 1)$ is
reductive. In this case, the cotangent space on each point
$\mbox{T}^{*}_m M$ is not irreducible under the action of the group
$\Lie{U}(n) \times 1$. In fact, $\mbox{T}^* M = \eta^{\perp} \oplus
\mathbb{R} \eta$ and
$$
\lie{so}(2n+1) \cong \Lambda^{2} \mbox{T}^* M = \Lambda^2
\eta^{\perp} \oplus \eta^{\perp} \wedge \mathbb R \eta.
$$
From now on, we will denote $X_{\zeta^{\perp}} = X - \eta(X) \zeta$,
for all $X  \in \mathfrak{X}(M)$. Since $\Lambda^2 \eta^{\perp} =
\lie{u}(n) \oplus \lie{u}(n)^{\perp}_{|\zeta^{\perp}} $, where
$\lie{u}(n)$ (resp., $\lie{u}(n)^{\perp}_{|\zeta^{\perp}} $)
consists of those two-forms $b$ such that $b(\varphi X , \varphi Y)
= b( X_{\zeta^{\perp}} ,  Y_{\zeta^{\perp}})$ (resp., $b(\varphi X ,
\varphi Y) = - b( X_{\zeta^{\perp}} , Y_{\zeta^{\perp}})$), we have
$$
\lie{so}(2n+1) = \lie{u}(n) \oplus \lie{u}(n)^{\perp}, \quad
\mbox{with } \lie{u}(n)^{\perp} =\lie{u}(n)^{\perp}_{|\zeta^{\perp}}
 \oplus \eta^{\perp} \wedge
\mathbb R \eta.
$$
Therefore, for the space $\mbox{T}^* M \otimes \lie{u}(n)^{\perp}$
of possible intrinsic $\Lie{U}(n) \times 1$-torsions,  we obtain
$$
\mbox{T}^* M \otimes  \lie{u}(n)^{\perp} = (\eta^{\perp} \otimes
\lie{u}(n)^{\perp}_{|\zeta^{\perp}}) \oplus (\eta \otimes
\lie{u}(n)^{\perp}_{|\zeta^{\perp}}) \oplus ( \eta^{\perp} \otimes
\eta^{\perp} \wedge   \eta  ) \oplus (\eta \otimes
  \eta^{\perp} \wedge   \eta).
$$
Chinea  and Gonz{\'a}lez-D{\'a}vila \cite{ChineaGonzalezDavila}
showed that $\mbox{T}^* M \otimes \lie{u}(n)^{\perp}$ is
decomposed into twelve irreducible $\Lie{U}(n)$-modules $\mathcal
C_1, \dots , \mathcal C_{12}$, where
\begin{eqnarray*}
\eta^{\perp} \otimes \lie{u}(n)^{\perp}_{|\zeta^{\perp}}
 & = &
\mathcal
C_1 \oplus \mathcal C_2 \oplus \mathcal C_3 \oplus \mathcal C_4, \\
\eta^{\perp} \otimes \eta^{\perp} \wedge  \eta  & = & \mathcal C_5
\oplus \mathcal C_8 \oplus \mathcal C_9 \oplus \mathcal C_6 \oplus
\mathcal C_7  \oplus \mathcal C_{10}, \\
 \eta \otimes \lie{u}(n)^{\perp}_{|\zeta^{\perp}}
& = & \mathcal C_{11} , \\
 \eta \otimes \eta^{\perp} \wedge   \eta & = & \mathcal
C_{12}.
\end{eqnarray*}
The modules  $\mathcal C_1, \dots , \mathcal C_4$ are isomorphic
to  the Gray and Hervella's $\Lie{U}(n)$-modules mentioned in
Remark \ref{clasif:struc}. Furthermore, note that $\varphi$
restricted to $\zeta^{\perp}$ works as an almost complex structure
and, if one considers the $\Lie{U}(n)$-action on the bilinear
forms $\otimes^2 \eta^{\perp}$, then we have the decomposition
$$
\textstyle \otimes^2 \eta^{\perp} = \mathbb R \langle \cdot ,
\cdot \rangle_{|\zeta^{\perp}} \oplus \lie{su}(n)_s \oplus
\real{\sigma^{2,0}} \oplus \mathbb R F  \oplus \lie{su}(n)_a
\oplus \lie{u}(n)^{\perp}_{|\zeta^{\perp}}.
$$
The modules  $\lie{su}(n)_s$ (resp., $\lie{su}(n)_a$) consists of
Hermitian symmetric (resp., skew-symmetric) bilinear forms
orthogonal to $\langle \cdot , \cdot \rangle_{|\zeta^{\perp}}$
(resp., $F$),
 and  $\real{\sigma^{2,0}}$ (resp., $\lie{u}(n)^{\perp}_{|\zeta^{\perp}}
 $) is the space of  anti-Hermitian
symmetric (resp., skew-symmetric) bilinear forms. With respect to
the modules $\mathcal C_i$, one has $\eta^{\perp} \otimes
\eta^{\perp} \wedge  \mathbb R \eta \cong   \otimes^2
\eta^{\perp}$ and, using the $\Lie{U}(n)$-map $\xi^{\Lie{U}(n)}
\to - \xi^{\Lie{U}(n)} \eta = \nabla \eta$, it is obtained
$$
\mathcal C_5 \cong  \mathbb R \langle \cdot , \cdot
\rangle_{|\zeta^{\perp}} , \quad \mathcal C_8 \cong \lie{su}(n)_s,
\quad \mathcal C_9 \cong \real{\sigma^{2,0}}, \quad \mathcal C_6
\cong \mathbb R F, \quad \mathcal C_7 \cong  \lie{su}(n)_a, \quad
\mathcal C_{10} \cong \lie{u}(n)^{\perp}_{|\zeta^{\perp}}.
$$
In summary,  the space of possible intrinsic torsions
  $\mbox{T}^* M\otimes  \lie{u}(n)^{\perp}$ consists of those tensors
$\xi^{\Lie{U}(n)}$ such that
\begin{equation} \label{inttorcar}
  \varphi \xi^{\Lie{U}(n)}_X Y + \xi^{\Lie{U}(n)}_X \varphi Y =
  \eta (Y) \varphi \xi^{\Lie{U}(n)}_X \zeta +
\eta ( \xi^{\Lie{U}(n)}_X \varphi Y) \zeta
\end{equation}
and, under the action of $U(n)\times 1$, is decomposed into:
\begin{enumerate}
\item if $n=1$, $ \xi^{\Lie{U}(1)} \in \mbox{T}^* M \otimes
\un(1)^\perp = \mathcal C_{5} \oplus \mathcal C_{6} \oplus
\mathcal C_{9} \oplus \mathcal C_{12}$; \item if $n=2$, $
\xi^{\Lie{U}(2)} \in \mbox{T}^* M \otimes \un(2)^\perp = \mathcal
C_{2} \oplus \mathcal C_{4} \oplus \dots \oplus \mathcal C_{12}$;
\item if $n \geqslant 3$, $ \xi^{\Lie{U}(n)} \in \mbox{T}^* M
\otimes \un(n)^\perp =
  \mathcal C_{1}  \oplus \dots  \oplus \mathcal C_{12}$.
\end{enumerate}

Now, we recall how  some of these classes are referred to  by
diverse authors \cite{Bl,ChineaGonzalezDavila}:

 $\{ \xi^{\Lie{U}(n)}=0 \}=$ cosymplectic manifolds,
 $\;\mathcal C_1=$ nearly-K-cosymplectic manifolds,
 $\;\mathcal C_5=$ $\alpha$-Kenmotsu manifolds,
 $\;\mathcal C_6=$ $\alpha$-Sasakian manifolds,
  $\;\mathcal C_5\oplus \mathcal C_6=$ trans-Sasakian manifolds,
 $\;\mathcal C_2 \oplus \mathcal C_9=$ almost cosymplectic manifolds,
  $\;\mathcal C_6 \oplus \mathcal C_7=$ quasi-Sasakian manifolds,
 $\;\mathcal C_1 \oplus \mathcal C_5 \oplus \mathcal C_6=$ nearly-trans-Sasakian manifolds,
 $\;\mathcal C_1 \oplus \mathcal C_2 \oplus \mathcal C_9 \oplus \mathcal C_{10}=$ quasi-K-cosymplectic manifolds,
 $\;\mathcal C_3 \oplus \mathcal C_4 \oplus \mathcal C_5 \oplus \mathcal C_{6}
  \oplus \mathcal C_7 \oplus \mathcal C_{8}=$ normal manifolds.
  \vspace{2mm}

The minimal $\Un(n)$-connection is given by $\nabla^{\Lie{U}(n)} =
\nabla + \xi^{\Lie{U}(n)}$, with
\begin{eqnarray} \label{torsion:xialmcont}
  \xi^{\Lie{U}(n)}_X  & = &  -  \tfrac12  \varphi \circ \nabla_X \varphi +
\nabla_X \eta \otimes \zeta - \tfrac12 \eta \otimes  \nabla_X \zeta \\
& = & \tfrac12 (\nabla_X \varphi) \circ \varphi   + \tfrac12
\nabla_X \eta \otimes \zeta - \eta \otimes \nabla_X \zeta.
\nonumber
\end{eqnarray}
In fact, one can firstly check that $\nabla^{\Lie{U}(n)}$ is metric
and $\nabla^{\Lie{U}(n)} \varphi = \nabla^{\Lie{U}(n)} \eta=0$.
Therefore,  $\nabla^{\Lie{U}(n)}$ is a $\Lie{U}(n)$-connection.
Finally, it is direct to see that $\xi^{\Lie{U}(n)}_X \in
\un(n)^\perp$. Note that if the almost contact metric structure is
of type $\mathcal{C}_5 \oplus \dots \oplus \mathcal{C}_{10} \oplus
\mathcal{C}_{12}$, then the expression for the intrinsic torsion is
reduced to
\begin{eqnarray} \label{torsion:xialmcont1}
  \xi^{\Lie{U}(n)}_X  & = & \nabla_X \eta \otimes \zeta - \eta \otimes \nabla_X
  \zeta.
\end{eqnarray}
For sake of simplicity, we will write $\xi = \xi^{\Lie{U}(n)}$ in
the sequel. Likewise, $\xi_{(i)}$ will denote the component of
$\xi$  obtained by the $\Lie{U}(n)$-isomorphism  $(\nabla F)_{(i)}
=(-\xi F)_{(i)}
 \in \mathcal C_i \to \xi_{(i)}$. In this way we
are using the same terminology  used in
\cite{ChineaGonzalezDavila} by Chinea and González-Dávila when we
are referring to classes.

For studying certain $\Lie{U}(n)$-components of the Riemannian
curvature tensor $R$ of an almost contact metric manifold, it is
necessary to consider a Ricci type tensor $\Ricac$,  called the {\it
almost contact Ricci tensor}, associated to the almost contact
metric structure. Such a tensor is defined by $\Ricac (X,Y) =
\langle R_{ e_i,X} \varphi e_i , \varphi Y \rangle$.

 In general, $\Ricac$ is not symmetric. However, since $\Ricac$  satisfies
 the identities
$$
\Ricac(\varphi X, \varphi Y) =
 \Ricac( Y_{\zeta^{\perp}}, X_{\zeta^{\perp}}), \qquad
 \Ricac(X, \zeta ) = 0,
$$
it can be claimed that
$$
\Ricac \in \mathbb{R} \langle \cdot , \cdot \rangle \oplus
\lie{su}(n)_s \oplus \lie{u}(n)^{\perp}_{|\zeta^{\perp}}
 \oplus
\eta^{\perp}_d \subseteq \mathbb{R} \langle \cdot , \cdot \rangle
\oplus  \lie{su}(n)_s \oplus
 \eta \odot
\eta^{\perp} \oplus \lie{u}(n)^{\perp}_{|\zeta^{\perp}}
 \oplus \eta \wedge \eta^{\perp},
$$
where $ \eta^{\perp}_d= \{ 2 \eta \odot \theta +  \eta \wedge
\theta \, | \, \theta \in \eta^{\perp} \} \cong \eta^{\perp}$ and
we follow the convention $a \odot b = \tfrac12 (a\otimes b + b
\otimes a)$.

The skew-symmetric part $\Ricac_{\rm alt}$ of
 $\Ricac$ will play a special r{\^o}le in the present work. Relative to  $\Ricac_{\rm alt}$  we
have the following result.
\begin{lemma}\label{astricciacm1}
  Let $(M,\langle \cdot ,\cdot \rangle, \varphi , \zeta)$ be a $2n+1$-dimensional
  almost contact metric manifold.  Then  the almost contact Ricci curvature
satisfies
  \begin{equation*}
    \begin{split}
  \Ricac_{\mbox{\rm \footnotesize alt}}( X_{\zeta^{\perp}} , Y_{\zeta^{\perp}})  =
 &
  \langle  (\nabla^{\Lie{U}(n)}_{e_i} \xi)_{\varphi e_i} \varphi X_{\zeta^{\perp}} , Y_{\zeta^{\perp}} \rangle
  + \langle \xi_{\xi_{e_i} \varphi e_i} \varphi X_{\zeta^{\perp}} , Y_{\zeta^{\perp}}
  \rangle,
   \\
 \Ricac( \zeta , X ) =
 &
 ( (\nabla^{\Lie{U}(n)}_{e_i} \xi)_{\varphi e_i} \eta )(\varphi X)
  +  ( \xi_{\xi_{e_i} \varphi e_i} \eta )(\varphi X),
    \end{split}
\end{equation*}
for all $X,Y \in \mathfrak{X}(M)$. Furthermore,  if $n>1$,  we
have:
\begin{enumerate}
\item[{\rm (i)}] The restriction $\Ricac_{{\rm alt} |
\zeta^\perp}$ of $\Ricac_{\mbox{\rm \footnotesize alt}}$ to the
space $\zeta^\perp$ is in $\lie{u}(n)^{\perp}_{|\zeta^{\perp}}$
and determines a $\Lie{U}(n)$-component of the Weyl curvature
tensor $W$.
  \item[{\rm (ii)}] The one-form $\Ricac( \zeta , \cdot  )$
   is in $\eta^\perp$ and determines another
   $\Lie{U}(n)$-component of $W$.
\end{enumerate}
As a consequence, if $(M,\langle \cdot ,\cdot \rangle)$ is
conformally flat, i.e. $W=0$, and  $n>1$, then $\Ricac_{{\rm alt}
| \zeta^\perp}=0$
 and $\Ricac( \zeta ,
\cdot  )=0$, or equivalently, $\Ricac_{\mbox{\rm \footnotesize
alt}}=0$.
\end{lemma}
\begin{proof}
The so-called Ricci formula \cite[p.~26]{Besse:Einstein} implies
  \begin{equation*}
  - (R_{e_i,\varphi e_i} F)(X,Y) =
    \talt(\nabla^2F)_{e_i,\varphi e_i}(X,Y),
  \end{equation*}
  where $\talt \colon T^* M \otimes T^* M \otimes \Lambda^2 T^* M \to
  \Lambda^2 T^* M \otimes \Lambda^2 T^* M$ is the skewing
  mapping.
  On one hand, by making use of first Bianchi's identity, it is
  straightforward to check that
  $$
  - (R_{e_i,\varphi e_i} F)(X,Y) = 4   \Ricac_{\mbox{\footnotesize alt}}( X , Y).
  $$
On the other hand, it is relatively direct to check that
\begin{eqnarray*}
 \talt(\nabla^2F)_{e_i,\varphi
 e_i}(X,Y) & = &-  2 \langle \varphi (\nabla^{\Lie{U}(n)}_{e_i} \xi)_{\varphi e_i}
 X ,   Y \rangle + 2 \langle (\nabla^{\Lie{U}(n)}_{e_i} \xi)_{\varphi e_i}
\varphi  X ,  Y \rangle \\
&& - 2 \langle \varphi \xi_{\xi_{e_i} \varphi e_i}
 X ,
  Y \rangle + 2 \langle \xi_{\xi_{e_i} \varphi e_i}
\varphi  X ,  Y \rangle.
\end{eqnarray*}
Now, using equation \eqref{inttorcar}, we will obtain the following
expression for $\Ricac_{\mbox{\rm \footnotesize alt}}( X , Y)$:
\begin{eqnarray} \label{otraricsac1}
     \Ricac_{\mbox{\footnotesize alt}}( X , Y) & = &   \langle (\nabla^{\Lie{U}(n)}_{e_i} \xi)_{\varphi e_i}
\varphi  X ,  Y \rangle
  +  \langle \xi_{\xi_{e_i} \varphi e_i} \varphi  X ,  Y \rangle
  \\
&&
    +  \eta \odot ((\nabla^{\Lie{U}(n)}_{e_i} \xi)_{\varphi e_i} \eta) \circ\varphi (X,Y)
   +  \eta \odot (\xi_{\xi_{e_i} \varphi e_i} \eta) \circ\varphi
   (X,Y). \nonumber
\end{eqnarray}
 Note that
$\nabla^{\Lie{U}(n)}_{e_i} \xi$ and $\xi$ are tensors of the same
type because $\nabla^{\Lie{U}(n)}$ is a $\Lie{U}(n)$-connection.
 Now, by replacing
$X=X_{\zeta^{\perp}}$ and $Y=Y_{\zeta^{\perp}}$ in equation
\eqref{otraricsac1}, we will obtain the first required identity.
Likewise, by replacing $X=\zeta $ and $Y=X$ in equation
\eqref{otraricsac1}, the second required identity follows.

Next we prove the final assertions in the Lemma. Denoting the
space of algebraic curvature tensors  by $\mathcal{R}$, we
consider the $\Lie{U}(n)$-map $\Phi_1  :
\lie{u}(n)^{\perp}_{|\zeta^{\perp}} \to \mathcal{R}$  defined by
\begin{equation*}
 \Phi_1 (b) = -\tfrac{1}{4(n+1)} \left( 6 (\varphi_{(1)} + \varphi_{(2)}) b \odot  F
 - (\varphi_{(1)}+\varphi_{(2)}) b \wedge F \right),
\end{equation*}
where $(\varphi_{(1)} b)(X,Y) = - b(\varphi X,Y)$ and
$(\varphi_{(2)} b)(X,Y) = - b(X,\varphi Y)$. It is easy to check
that
$$
\Ric \Phi_1 (b) = 0, \qquad \Ric^{\rm ac}  \Phi_1 (b) =b.
$$
Therefore,  $\Phi_1 ( \lie{u}(n)^{\perp}_{|\zeta^{\perp}})$ is
contained in the space $\mathcal{W} \subseteq \mathcal{R}$ of
algebraic Weyl curvature tensors and is isomorphic to the
$\Lie{U}(n)$-module denoted by $\real{\lambda^{2,0}}$ in
references. This notation is described below in Section
\ref{harmonicalmostcontact}. The tensor $\Phi_1( \Ricac_{{\rm alt}
| \zeta^\perp})$ is the component of $W$ included in
$\real{\lambda^{2,0}}\subseteq \mathcal{W}$.

Next we consider the the $\Lie{U}(n)$-map $\Phi_2  : \eta^\perp
\to \mathcal{R}$  given  by
$$
\Phi_2 (\theta) = 6 \left( \eta \wedge (\theta \circ
\varphi)\right) \odot
                F -  \eta \wedge (\theta \circ \varphi) \wedge F
                - \tfrac{6}{2n-1}  (\eta \odot \theta)
 \ovee \langle \cdot , \cdot \rangle,
$$
where $\ovee$ denotes the usual Kulkarni-Nomizu product
\cite{Besse:Einstein} defined by
\begin{equation} \label{KulNom}
(a \ovee b)(x,y,z,w)=a(x,z)b(y,w) -a(y,z)b(x,w) +a(y,w)b(x,z)
-a(x,w)b(y,z).
\end{equation}
One can check that
$$
\Ric \Phi_2(\theta) =0, \qquad  \Ric^{{\rm ac}} \Phi_2(\theta) =
\tfrac{4(n^2-1)}{2n-1} \eta \otimes \theta.
$$
Therefore,  $\Phi_1 \left( \eta^\perp\right)$ is contained in
$\mathcal{W} \subseteq \mathcal{R}$ and is isomorphic to the
$\Lie{U}(n)$-module  $\real{\lambda^{1,0}}$. The tensor $\Phi_2(
\tfrac{2n-1}{4(n^2-1)} \Ricac (\zeta, \cdot))$ is the component of
$W$ included in $\real{\lambda^{1,0}}\subseteq \mathcal{W}$.
\end{proof}

The vector field $\xi_{e_i} \varphi
 e_i$ involved in  $\Ricac$ is given by
$$
- 2 \xi_{e_i} \varphi e_i =  (d^* F)^\sharp  + d^* F(\zeta) \zeta
+ \varphi \nabla_\zeta \zeta.
$$
Thus, this vector field is contributed by the components of $\xi$ in
$\mathcal C_4$ and   $\mathcal C_6$. In fact, $2 \xi_{{(4)}e_i}
\varphi e_i = - (d^* F)^{\sharp} - \varphi \nabla_\zeta \zeta + d^*
F(\zeta) \zeta$ and $\xi_{{(6)}e_i} \varphi e_i = - d^* F(\zeta)
\zeta$. On the other hand, the vector field $\xi_{e_i} e_i$ which
takes part in the harmonicity criterion is given by
$$
- 2  \xi_{e_i} e_i  =  \varphi (d^* F)^{\sharp} + 2 d^* \eta \;
\zeta +  \nabla_{\zeta}\zeta.
$$
Because  $2 \xi_{{(4)}e_i} e_i = - \varphi (d^* F)^{\sharp} +
\nabla_{\zeta}\zeta $, $ \xi_{{(5)}e_i} e_i = - d^* \eta \; \zeta$
and $ \xi_{{(12)}e_i} e_i = -
 \nabla_{\zeta}\zeta$, in this case we have a vector field contributed by
$\mathcal C_4$, $\mathcal C_5$ and $\mathcal C_{12}$.


\section{Harmonic almost contact structures}{\indent}
\setcounter{equation}{0} \label{harmonicalmostcontact}
  In this section we will show
conditions relating harmonicity and Chinea and González-Dávila's
classes of almost contact metric structures. Firstly, we
characterise harmonic almost contact structures by means of the
rough Laplacian of stabilised tensors under the action of
$\Lie{U}(n)$.
\begin{theorem} \label{characharmalmcontact}
If  $(M,  \langle \cdot, \cdot \rangle,\varphi, \zeta )$ is  a
$2n+1$-dimensional  almost contact metric manifold  with fundamental
two-form $F$, then the following conditions are equivalent:
 \begin{enumerate}
 \item[{\rm (i)}] $(M,  \langle \cdot, \cdot \rangle,\varphi, \zeta )$  is harmonic.
\vspace{1mm}

\item[{\rm (ii)}]  $\nabla^* \nabla \varphi \in \lie{u}(n) +
\zeta^{\perp}_c$, where $\zeta^{\perp}_c = \{a \otimes \zeta -
\eta \otimes a^{\sharp} \,|\, a \in \eta^{\perp} \} \cong
\eta^{\perp} \wedge \eta$, and  $ \nabla^{*} \nabla \zeta
=-\xi_{e_i} \xi_{e_i} \zeta$.

  \item[{\rm (iii)}]
   $[\nabla^* \nabla \varphi, \varphi] =
  3 \eta \otimes \nabla^* \nabla \zeta - 3 \nabla^* \nabla \eta \otimes \zeta $, where $[,]$
  denotes the commutator bracket for endomorphisms, and $ \nabla^{*} \nabla \zeta
=-\xi_{e_i} \xi_{e_i} \zeta$. \vspace{1mm}

 \item[{\rm (iv)}]$\nabla^* \nabla F \in \lie{u}(n) +
\eta^{\perp} \wedge \eta$, i.e. $\nabla^* \nabla F(\varphi X,
\varphi Y ) = \nabla^* \nabla F(X_{\zeta^{\perp}},
Y_{\zeta^{\perp}}  )$, and $\nabla^{*} \nabla \eta = -
\xi_{e_i}(\xi_{e_i} \eta)$. \vspace{1mm}

\item[{\rm (v)}] $\nabla^* \nabla F(\varphi X, \varphi Y ) =
\nabla^* \nabla F( X,  Y ) - 3 \eta \wedge \left( \nabla^* \nabla
\eta \right) \circ \varphi ( X , Y)$ and $ \nabla^{*} \nabla \eta
=-\xi_{e_i} (\xi_{e_i} \eta)$.

 \item[{\rm (vi)}]$\nabla^* \nabla F (X,Y)  =
 - 4 F( \xi_{e_i}X ,\xi_{e_i}Y) + (
 \xi_{e_i} \eta) \wedge (
 \xi_{e_i} \eta) \circ \varphi (X,Y)+ \eta \wedge (\xi_{e_i}  (\xi_{e_i} \eta)) \circ
\varphi  (X,Y).$

\item[{\rm (vii)}] For all $X, Y \in \mathfrak X ( M)$, we  have
$\langle  (\nabla^{\Lie{U}(n)}_{e_i} \xi)_{ e_i}  X_{\zeta^{\perp}} , Y_{\zeta^{\perp}} \rangle
  + \langle \xi_{\xi_{e_i}  e_i}  X_{\zeta^{\perp}} , Y_{\zeta^{\perp}}
  \rangle = 0\,$ and $\,(\nabla^{\Lie{U}(n)}_{e_i} \xi)_{e_i} \eta +   \xi_{\xi_{e_i} e_i} \eta
   =0$.
\end{enumerate}
In  particular, an almost contact metric structure of type
$\mathcal{C}_{5} \oplus \ldots \oplus \mathcal{C}_{10} \oplus
\mathcal{C}_{12}$ is harmonic if and only if $\nabla^* \nabla \zeta
= \| \nabla \zeta \|^2 \zeta$, that is, the characteristic vector
field $\zeta$ is a harmonic unit vector field {\rm (see
\cite{Wie1,GMS} for this notion)}.
 \end{theorem}
 \begin{remark} {\rm
Vergara-Díaz and Wood \cite{VergaraWood} characterise harmonic
almost contact metric structures by the conditions  $\nabla^*
\nabla \varphi \in \lie{u}(n) + \zeta^{\perp}_c$,  and $\nabla^{*}
\nabla \zeta = \|\nabla \zeta\|^2 \zeta  - \frac12 \varphi
T(\varphi)$. In terms given here, one can check that $ \varphi
T(\varphi) = 2 (\|\nabla \zeta\|^2 \zeta + \xi_{e_i} \xi_{e_i}
\zeta)$}.
 \end{remark}

\begin{proof}
For (i) implies (vi). Because $F$ is stabilised by $\Lie{U}(n)$,
by $(\nabla^{\Lie{U}(n)}_{e_i} \xi)_{ e_i} = -  \xi_{\xi_{e_i}
e_i}$ and \eqref{lapstaten}, (i) implies
$$
(\nabla^* \nabla F)(X,Y) = - 2 F( \xi_{e_i}X ,\xi_{e_i}Y) +
\langle \xi_{e_i} X , \xi_{e_i} \varphi Y \rangle - \langle
\xi_{e_i} \varphi X , \xi_{e_i}  Y \rangle.
$$
Now, using \eqref{inttorcar} we have (vi).

 For (vi) implies (v). In general, it is
satisfied  $ (\xi_X F)(\zeta , \varphi Y) = (\xi_X \eta)(Y)$.
Taking this into account, it is direct to check
\begin{equation} \label{nanaeta}
(\nabla^* \nabla F)(\zeta , \varphi X) = (\nabla^* \nabla \eta)(X)
+ 4 (\xi_{e_i} (\xi_{e_i} \eta))(X) +3 \|\nabla \eta\|^2 \eta(X).
\end{equation}
Then, replacing $X=\zeta$ and $Y= \varphi X$ in (vi) and  using
again \eqref{inttorcar}, it follows $\nabla^{*} \nabla \eta = -
\xi_{e_i}(\xi_{e_i} \eta)$. The remaining identity for $(\nabla^*
\nabla F)(\varphi X , \varphi Y)$ is directly deduced from (vi).

For (v) implies (i). Since  $ \nabla^{*} \nabla \eta =-\xi_{e_i}
(\xi_{e_i} \eta)$, we have $( \nabla^{\Lie{U}(n)}_{e_i} \xi)_{e_i}
+ \xi_{\xi_{e_i}e_i}) \eta =0$. Now, using \eqref{inttorcar},  we
can point out that, in general,
\begin{eqnarray*}
 \left( \xi_{e_i} (\xi_{e_i}
 F)\right) (\varphi X , \varphi Y) & = &  \left( \xi_{e_i} (\xi_{e_i}
 F)\right) ( X ,  Y) - 3  \eta \wedge (\xi_{e_i}\xi_{e_i}
 \eta)
 \circ \varphi (X,Y).
\end{eqnarray*}
Making use of this identity, from (v) one deduces
$$
- 2  (( \nabla^{\Lie{U}(n)}_{e_i} \xi)_{e_i} + \xi_{\xi_{e_i}e_i})
F  = \eta \wedge \left( ( \nabla^{\Lie{U}(n)}_{e_i} \xi)_{e_i} +
\xi_{\xi_{e_i}e_i}) \eta \right) \circ \varphi =0.
$$
But the map $A \to - F(A \cdot, \cdot) - F (\cdot, A \cdot) $,
from $\lie{u}(n)^{\perp} \subseteq \lie{so}(2n+1)$ to
 $\lie{u}(n)^{\perp} \subseteq \Lambda^2 T^* M$,  is an
 $U(n)$-isomorphism. Therefore, $( \nabla^{\Lie{U}(n)}_{e_i} \xi)_{e_i} +
\xi_{\xi_{e_i}e_i}=0$. Thus we have (i).

For the equivalence between (iv) and  (v). It is immediate that
(v) implies (iv). Conversely, using \eqref{nanaeta}, it is not
hard to see that (iv) implies (v).

For the equivalence between (iii) and (v). Because $(\nabla_X
F)(Y,Z) = \langle Y , (\nabla_X \varphi) Z \rangle$, it is direct
to see that $(\nabla^* \nabla F)(X,Y) = \langle X , (\nabla^*
\nabla \varphi)(Y) \rangle$. Then using this and \eqref{nanaeta},
the equivalence follows. Finally, the equivalences between (ii)
and (iii)  and between (i) and (vii) are immediate.

The final remark included in the Theorem is an immediate consequence
of the computation of $d^* \xi$ by using  equation
\eqref{torsion:xialmcont1}.
\end{proof}

In next two results, for certain types of almost contact metric
structures, we will deduce  conditions characterising harmonic
almost contact structures.
\begin{theorem} \label{harmclasscontact}
For a $2n+1$-dimensional almost contact metric manifold
$(M,\langle \cdot ,\cdot \rangle, \varphi , \zeta)$,  we have:
\begin{enumerate}
 \item[{\rm (i)}] If $M$ is
  of type $\mathcal D$, where $\mathcal D = \mathcal C_1 \oplus \mathcal C_2 \oplus  \mathcal C_5
\oplus  \mathcal C_{6}\oplus  \mathcal C_7 \oplus  \mathcal
C_{8}$, $
 \mathcal
C_1 \oplus \mathcal C_2 \oplus  \mathcal C_9 \oplus  \mathcal
C_{10}$, then the almost contact metric structure is harmonic if and
only if $ \Ricac_{\rm alt}( X_{\zeta^{\perp}} , Y_{\zeta^{\perp}} )
= 0$ and $ \Ricac( \zeta , X ) =  0$, for all $X,Y \in \mathfrak
X(M)$.

 \item[{\rm (ii)}]
 For $n \neq 2$, if $M$ is
 of type  $\mathcal C_1 \oplus  \mathcal C_4 \oplus \mathcal C_5
\oplus \mathcal C_6 \oplus \mathcal C_7 \oplus \mathcal C_8$, then
the almost contact metric structure is harmonic if and only if
\begin{equation*}
     \begin{split}
 (n-1)(n-5) \Ricac_{\rm alt}( X_{\zeta^{\perp}} , Y_{\zeta^{\perp}})  =
   & \;
   2(n+1)(n-3) \langle  \xi_{ \xi_{e_i} e_i} X_{\zeta^{\perp}} , Y_{\zeta^{\perp}}
   \rangle,
  \\
 \Ricac( \zeta , X ) = &\; - 2 ( \xi_{ \xi_{e_i}  e_i} \eta )( X),
    \end{split}
 \end{equation*}
 for all $X,Y \in \mathfrak X(M)$.

 \item[{\rm (iii)}]
For $n \neq 2$, if $M$ is
 of type  $\mathcal C_1 \oplus  \mathcal C_4 \oplus \mathcal C_9
\oplus \mathcal C_{10}$, then almost contact metric structure is
harmonic if and only if
\begin{gather*}
 \quad \quad \quad  (n-1)(n-5) \Ricac_{\rm alt}( X_{\zeta^{\perp}} , Y_{\zeta^{\perp}})  =
   2(n+1)(n-3) \langle  \xi_{ \xi_{e_i} e_i} X_{\zeta^{\perp}} , Y_{\zeta^{\perp}}
   \rangle,
  \\
 \Ricac( \zeta , X ) = 2 ( \xi_{ \xi_{e_i}  e_i} \eta )( X_{\zeta^{\perp}}),
  \end{gather*}
 for all $X,Y \in \mathfrak X(M)$.
 \item[{\rm (iv)}] For $n \neq 2$, if $M$ is
 of type  $ \mathcal C_2 \oplus  \mathcal C_4 \oplus \mathcal C_5
\oplus \mathcal C_6  \oplus \mathcal C_7 \oplus \mathcal C_8$,
then the almost contact metric structure is harmonic if and only
if
\begin{equation*}
     \begin{split}
(n-1)  \Ricac_{\rm alt}( X_{\zeta^{\perp}} , Y_{\zeta^{\perp}}) =
   & \;  2n \langle  \xi_{ \xi_{e_i} e_i} X_{\zeta^{\perp}} , Y_{\zeta^{\perp}}
   \rangle,
  \\
 \Ricac( \zeta , X ) = & - 2  ( \xi_{ \xi_{e_i}  e_i} \eta )( X),
    \end{split}
\end{equation*}
 for all $X,Y \in \mathfrak X(M)$.

 \item[{\rm (v)}] For $n \neq 2$, if $M$ is
 of type  $\mathcal C_2 \oplus  \mathcal C_4 \oplus \mathcal C_9
\oplus \mathcal C_{10}$, then the almost contact metric structure
is harmonic if and only if
\begin{equation*}
    \begin{split}
  (n-1) \Ricac_{\rm alt}( X_{\zeta^{\perp}} , Y_{\zeta^{\perp}})  =
   &  \; 2n \langle  \xi_{ \xi_{e_i} e_i} X_{\zeta^{\perp}} , Y_{\zeta^{\perp}}
   \rangle,
   \qquad
 \Ricac( \zeta , X ) =  2( \xi_{ \xi_{e_i}  e_i} \eta )( X),
    \end{split}
\end{equation*}
 for all $X,Y \in \mathfrak X(M)$.

 \item[{\rm (vi)}]
If $M$ is
 normal $( \mathcal C_3 \oplus  \mathcal C_4 \oplus \mathcal C_5
\oplus \mathcal C_6  \oplus \mathcal C_7 \oplus \mathcal C_8)$,
then the almost contact metric structure is harmonic if and only
if
\begin{equation*}
    \begin{split}
  \Ricac_{\rm alt}( X_{\zeta^{\perp}} , Y_{\zeta^{\perp}}
)  = &  - 2 \langle  \xi_{ \xi_{e_i} e_i } X_{\zeta^{\perp}} ,
Y_{\zeta^{\perp}} \rangle,   \qquad
 \Ricac( \zeta , X ) =   -  2 ( \xi_{ \xi_{e_i}  e_i} \eta )( X)
 ,
    \end{split}
\end{equation*}
for all $X,Y \in \mathfrak X(M)$.

\item[{\rm (vii)}]
 If $M$ is
 of type  $\mathcal C_3 \oplus  \mathcal C_4 \oplus \mathcal C_9
\oplus \mathcal C_{10}$,
 then the almost contact
metric structure is harmonic if and only if
\begin{equation*}
     \begin{split}
  \Ricac_{\rm alt}( X_{\zeta^{\perp}} , Y_{\zeta^{\perp}}
)  = &  - 2 \langle  \xi_{ \xi_{e_i} e_i } X_{\zeta^{\perp}} ,
Y_{\zeta^{\perp}} \rangle,  \qquad \Ricac( \zeta , X ) = 2 ( \xi_{
\xi_{e_i}  e_i} \eta )( X),
    \end{split}
\end{equation*}
for all $X,Y \in \mathfrak X(M)$.
    \item[{\rm (viii)}]
    If $M$ is of type $\mathcal D$, where
 $\mathcal D =\mathcal C_1 \oplus  \mathcal C_5 \oplus \mathcal C_9$,
  $\mathcal C_1 \oplus  \mathcal C_6 \oplus \mathcal C_8$,
then the almost contact metric structure is harmonic if and only if
$ \Ricac( \zeta , X ) =  0$, for all $X \in \mathfrak X(M)$.

  \item[{\rm (ix)}] For $n \neq 2$, if $M$ is of type $\mathcal D$, where
 $\mathcal D =\mathcal C_4 \oplus  \mathcal C_5 \oplus \mathcal C_6$,
  $\mathcal C_4 \oplus  \mathcal C_5 \oplus \mathcal C_7$,
   $\mathcal C_4 \oplus  \mathcal C_5 \oplus \mathcal C_9$,
    $\mathcal C_4 \oplus  \mathcal C_8$,
then the almost contact metric structure is harmonic if and only if
$ \Ricac( \zeta , X ) =  0$, for all $X \in \mathfrak X(M)$.
\end{enumerate}
In particular:
 \begin{enumerate}
 \item[{\rm (i)$^*$}]
  If $M$ is of type $\mathcal D$, where
 $\mathcal D =\mathcal C_1 \oplus  \mathcal C_5$,
  $\mathcal C_1 \oplus  \mathcal C_8 $,
  $\mathcal C_1 \oplus  \mathcal C_9 $,
  $\mathcal C_3 \oplus \mathcal C_6$,
$\mathcal C_3 \oplus \mathcal C_7$,
 $ \mathcal C_3 \oplus \mathcal C_{10}$,
 $ \mathcal C_5 \oplus \mathcal C_{6} \oplus \mathcal C_{7}$,
 $ \mathcal C_5 \oplus \mathcal C_{8}$,
 $ \mathcal C_5 \oplus \mathcal C_{9}$,
 $ \mathcal C_5 \oplus \mathcal C_{10}$,
 $ \mathcal C_6 \oplus \mathcal C_{7} \oplus \mathcal C_{8}$,
 $ \mathcal C_6 \oplus \mathcal C_7 \oplus \mathcal C_{10}$,
 $ \mathcal C_8 \oplus \mathcal C_{9}$,
 $ \mathcal C_9 \oplus \mathcal C_{10}$,
then the almost contact metric structure is harmonic.
  \item[{\rm (ii)$^*$}]
  For $n \neq 2$, if $M$ is of type $\mathcal D$, where
   $\mathcal D = \mathcal C_4 \oplus \mathcal C_5$,
  $\mathcal C_4 \oplus \mathcal C_6 $,
   $\mathcal C_4 \oplus \mathcal C_7$,
   $\mathcal C_4 \oplus  \mathcal C_9$,
    then the almost contact metric structure is harmonic.
  \end{enumerate}
\end{theorem}

\begin{corollary} $\,$ \label{corharmonic}
\begin{enumerate}
\item[{\rm (a)}]If an almost contact metric structure is of type
$\mathcal{D}$, for $\mathcal{D} =  \mathcal C_5 \oplus \mathcal
C_{6} \oplus \mathcal C_{7}$,
 $ \mathcal C_5 \oplus \mathcal C_{8}$,
 $ \mathcal C_5 \oplus \mathcal C_{9}$,
 $ \mathcal C_5 \oplus \mathcal C_{10}$,
 $ \mathcal C_6 \oplus \mathcal C_{7} \oplus \mathcal C_{8}$,
 $ \mathcal C_6 \oplus \mathcal C_7 \oplus \mathcal C_{10}$,
 $ \mathcal C_8 \oplus \mathcal C_{9}$,
 $ \mathcal C_9 \oplus \mathcal C_{10}$, then the characteristic vector
field $\zeta$ is a harmonic unit vector field.
 \item[{\rm (b)}] For a conformally flat
manifold $(M,\langle \cdot ,\cdot \rangle)$ of dimension $2n+1$ with
$n>1$:
\begin{enumerate}
\item[{\rm (i)}] If an almost contact  structure compatible with
$\langle \cdot ,\cdot \rangle$ is of type $\mathcal D$, where
$\mathcal D = \mathcal C_1 \oplus \mathcal C_2 \oplus \mathcal C_5
\oplus \mathcal C_{6}\oplus \mathcal C_7 \oplus \mathcal C_{8}$, $
 \mathcal
C_1 \oplus \mathcal C_2 \oplus  \mathcal C_9 \oplus  \mathcal
C_{10}$,
 $\mathcal C_1 \oplus  \mathcal C_5 \oplus \mathcal C_9$,
 then it is harmonic.
 \item[{\rm (ii)}] For $n>2$, if an  almost contact structure compatible with  $\langle \cdot
,\cdot \rangle$  is of type $\mathcal D$,  where
 $\mathcal D =\mathcal C_4 \oplus  \mathcal C_5 \oplus \mathcal C_6$,
  $\mathcal C_4 \oplus  \mathcal C_5 \oplus \mathcal C_7$,
   $\mathcal C_4 \oplus  \mathcal C_5 \oplus \mathcal C_9$,
    $\mathcal C_4 \oplus  \mathcal C_8$, then it is harmonic.
 \end{enumerate}
 \end{enumerate}
\end{corollary}

 In order to prove these last results,   we will show
some identities for almost contact metric structures which are
analog to the one satisfied for almost Hermitian structures (see
\cite[p. 182]{FMCAS}
 or \cite[Equation (4.2)]{GDMC}). They are consequences from the fact $d^2 F=0$ and arise
 when we write $d^2 F$ by means of $\nabla^{\Lie{U}(n)}$ and
 $\xi$.
  For sake of simplicity, we will
deduce such identities for almost contact metric  structures of type
$\mathcal C_1 \oplus \ldots \oplus \mathcal C_{10}$, i.e.
$\xi_{\zeta}=0$. In fact, let us describe an expression for $d^2 F$
in terms of the intrinsic torsion. Since
$$
\tfrac12 \left( \nabla_Y F \right) (Z,W) =  \langle \xi_Y Z ,
\varphi W \rangle -  \eta \odot (\xi_Y \eta) \circ \varphi (Z,W),
$$
by alternating this identity,  we get
\begin{equation*}
\begin{array}{rl}
\tfrac12 \, d F (Y,Z,W) =&
  \langle \xi_Y Z , \varphi W \rangle
   + \langle \xi_W Y , \varphi Z \rangle
   + \langle \xi_Z W , \varphi Y \rangle
    -  \eta \odot (\xi_Y \eta) \circ \varphi (Z,W)
\\
  &
 - \eta \odot (\xi_W \eta) \circ \varphi (Y,Z)
  - \eta \odot (\xi_Z \eta) \circ \varphi (W,Y).
\end{array}
\end{equation*}
 Now, taking $\nabla = \nabla^{\Lie{U}(n)} - \xi$ into account,
we have $d^2 F = \mbox{alt} (\nabla^{\Lie{U}(n)}  d F) - \mbox{alt}
(\xi d F)$, where $\mbox{alt}$ is the alternation map. On the one
hand, it is direct to check that $\mbox{alt} (\nabla^{\Lie{U}(n)} d
F)$ is twice the alternation of the tensor $H$  given by
\begin{equation} \label{ec1}
H(X,Y,Z,W) = \langle \left( \nabla^{\Lie{U}(n)}_X \xi \right)_Y Z,
\varphi W \rangle + \eta \odot \left( (\nabla^{\Lie{U}(n)}_X
\xi)_Y \eta \right) \circ \varphi (Z,W).
\end{equation}
That is,
\begin{equation} \label{ec1bis}
{\footnotesize  \begin{array}{l}
\hspace{-.5cm} \tfrac12\mbox{alt} (\nabla^{\Lie{U}(n)}  d F)(X,Y,Z,W) =  \\
 \langle \left( \nabla^{\Lie{U}(n)}_X \xi \right)_Y Z, \varphi W
\rangle
 + \langle \left( \nabla^{\Lie{U}(n)}_X \xi \right)_Z W, \varphi Y \rangle
 + \langle \left( \nabla^{\Lie{U}(n)}_X \xi \right)_W Y, \varphi Z \rangle
  - \langle \left( \nabla^{\Lie{U}(n)}_Y \xi \right)_X Z, \varphi W \rangle
 \\[1mm]
  - \langle \left( \nabla^{\Lie{U}(n)}_Y \xi \right)_Z W, \varphi X \rangle
  - \langle \left( \nabla^{\Lie{U}(n)}_Y \xi \right)_W X, \varphi Z \rangle
  + \langle \left( \nabla^{\Lie{U}(n)}_Z \xi \right)_X Y, \varphi W \rangle
  + \langle \left( \nabla^{\Lie{U}(n)}_Z \xi \right)_Y W, \varphi X \rangle
 \\[1mm]
  + \langle \left( \nabla^{\Lie{U}(n)}_Z \xi \right)_W X, \varphi Y \rangle
   - \langle \left( \nabla^{\Lie{U}(n)}_W \xi \right)_X Y, \varphi Z \rangle
  - \langle \left( \nabla^{\Lie{U}(n)}_W \xi \right)_Y Z, \varphi X \rangle
  - \langle \left( \nabla^{\Lie{U}(n)}_W \xi \right)_Z X, \varphi Y \rangle
 \\[2mm]
 + \eta \odot \left( (\nabla^{\Lie{U}(n)}_X \xi )_Y \eta \right) \circ \varphi (Z,W)
 + \eta \odot \left( (\nabla^{\Lie{U}(n)}_X \xi )_Z \eta \right) \circ \varphi (W,Y)
 + \eta \odot \left( (\nabla^{\Lie{U}(n)}_X \xi )_W \eta \right) \circ \varphi (Y,Z)
 \\[1mm]
 - \eta \odot \left( (\nabla^{\Lie{U}(n)}_Y \xi )_X \eta \right) \circ \varphi (Z,W)
 - \eta \odot \left( (\nabla^{\Lie{U}(n)}_Y \xi )_Z \eta \right) \circ \varphi (W,X)
 - \eta \odot \left( (\nabla^{\Lie{U}(n)}_Y \xi )_W \eta \right) \circ \varphi (X,Z)
 \\[1mm]
 + \eta \odot \left( (\nabla^{\Lie{U}(n)}_Z \xi )_X \eta \right) \circ \varphi (Y,W)
 + \eta \odot \left( (\nabla^{\Lie{U}(n)}_Z \xi )_Y \eta \right) \circ \varphi (W,X)
 + \eta \odot \left( (\nabla^{\Lie{U}(n)}_Z \xi )_W \eta \right) \circ \varphi (X,Y)
 \\[1mm]
  - \eta \odot \left( \nabla^{\Lie{U}(n)}_W \xi )_X \eta \right) \circ \varphi (Y,Z)
  - \eta \odot \left( \nabla^{\Lie{U}(n)}_W \xi )_Y \eta \right) \circ \varphi (Z,X)
  - \eta \odot \left( \nabla^{\Lie{U}(n)}_W \xi )_Z \eta \right) \circ \varphi
  (X,Y).
\end{array} }
\end{equation}

 On the other hand, it is direct to see that
\begin{equation} \label{ec2}
{\footnotesize \begin{split}
&
 \hspace{-1cm}  - \tfrac12  \mbox{alt} (\xi d
F)(X,Y,Z,W) = \\
 & \langle \xi_X Y, e_i \rangle \langle
\xi_{e_i} Z , \varphi W \rangle +\langle \xi_X Z, e_i \rangle
\langle \xi_{e_i} W , \varphi Y \rangle + \langle \xi_X W, e_i
\rangle \langle \xi_{e_i} Y, \varphi Z \rangle
\\[1mm]
 &-\langle \xi_Y X, e_i \rangle \langle \xi_{e_i} Z , \varphi W \rangle
 -\langle \xi_Y Z, e_i \rangle \langle \xi_{e_i} W ,\varphi X \rangle
   -\langle \xi_Y W, e_i \rangle \langle \xi_{e_i} X , \varphi Z \rangle
\\[1mm]
 &+ \langle \xi_Z W, e_i \rangle \langle \xi_{e_i} X, \varphi Y \rangle
  + \langle \xi_Z X, e_i \rangle \langle \xi_{e_i} Y , \varphi W \rangle
 + \langle \xi_Z Y, e_i \rangle \langle \xi_{e_i} W , \varphi X \rangle
\\[1mm]
 &-\langle \xi_W X, e_i \rangle \langle \xi_{e_i} Y , \varphi Z \rangle
  -\langle \xi_W Y, e_i \rangle \langle \xi_{e_i} Z , \varphi X \rangle
 -\langle \xi_W Z, e_i \rangle \langle \xi_{e_i} X , \varphi Y \rangle
  \\[1mm]
 &+  (\xi_Y \eta)(Z) (\xi_X \eta)\circ \varphi (W)
 +  (\xi_Z \eta)(W) (\xi_X \eta)\circ \varphi (Y)
 +  (\xi_W \eta)(Y) (\xi_X \eta)\circ \varphi (Z)
 \\[1mm]
 &-  (\xi_X \eta)(Z) (\xi_Y \eta)\circ \varphi (W)
 -  (\xi_Z \eta)(W) (\xi_Y \eta)\circ \varphi (X)
 -  (\xi_W \eta)(X) (\xi_Y \eta)\circ \varphi (Z)
 \\[1mm]
 &+  (\xi_X \eta)(Y) (\xi_Z \eta)\circ \varphi (W)
 +  (\xi_Y \eta)(W) (\xi_Z \eta)\circ \varphi (X)
 +  (\xi_W \eta)(X) (\xi_Z \eta)\circ \varphi (Y)
 \\[1mm]
 &-  (\xi_X \eta)(Y) (\xi_W \eta)\circ \varphi (Z)
 -  (\xi_Y \eta)(Z) (\xi_W \eta)\circ \varphi (X)
 -  (\xi_Z \eta)(X) (\xi_W \eta)\circ \varphi (Y)
  \\[2mm]
  &-  (\xi_X (\xi_Y \eta)) \circ \varphi (Z) \eta(W)
  -  (\xi_X (\xi_Z \eta)) \circ \varphi (W) \eta(Y)
  -  (\xi_X (\xi_W \eta)) \circ \varphi (Y) \eta(Z)
  \\[1mm]
  &+  (\xi_Y (\xi_X \eta)) \circ \varphi (Z) \eta(W)
  +  (\xi_Y (\xi_Z \eta)) \circ \varphi (W) \eta(X)
  +  (\xi_Y (\xi_W \eta)) \circ \varphi (X) \eta(Z)
  \\[1mm]
  &-  (\xi_Z (\xi_X \eta)) \circ \varphi (Y) \eta(W)
  -  (\xi_Z (\xi_Y \eta)) \circ \varphi (W) \eta(X)
  -  (\xi_Z (\xi_W \eta)) \circ \varphi (X) \eta(Y)
  \\[1mm]
  &+  (\xi_W (\xi_X \eta)) \circ \varphi (Y) \eta(Z)
  +  (\xi_W (\xi_Y \eta)) \circ \varphi (Z) \eta(X)
  +  (\xi_W (\xi_Z \eta)) \circ \varphi (X) \eta(Y).
\end{split}}
\end{equation}
Now, for our purposes, we need to  note that
 $ d^2 F \in  \Lambda^4T^*M $ and  $\Lambda^4T^*M$ is equal to
 $
 \real{\lambda^{4,0}} \oplus \real{\lambda^{3,1}}\oplus \real{\lambda^{2,0}} \wedge F \oplus
  [\lambda^{2,2}_0] \oplus [\lambda^{1,1}_0]\wedge  F \oplus \mathbb R F\wedge F
  \oplus \real{\lambda^{3,0}} \wedge \eta \oplus \real{\lambda^{2,1}_0}
\wedge \eta
   \oplus \real{\lambda^{1,0}} \wedge  F \wedge \eta.
   $
We recall that $\lambda_0^{p,q}$ is a complex irreducible
$\Un(n)$-module coming from the $(p,q)$-part of the complex
exterior algebra and its corresponding dominant weight in standard
coordinates is given by $(1, \dots , 1,0, \dots , 0 , -1 , \dots ,
-1)$, where $1$ and $-1$ are repeated $p$ and $q$ times,
respectively.   The notation~$\real V$ means the real vector space
underlying a complex vector space~$V$ and $[W]$~denotes a real
vector space which admits~$W$ as its complexification.

From now on, we are assuming  that the almost contact metric
structure is of type $\mathcal C_1 \oplus \ldots \oplus \mathcal
C_{10}$, i.e. $\xi_{\zeta}=0$. Our purpose is to compute the
$\real{\lambda^{2,0}}$-part  of $d^2 F$.
 This can be obtained from the two-form
$\pi_{\real{\lambda^{2,0}}} \circ F_{1\,2} ( d^2 F)$, where
$F_{1\,2}$ is the  $1\, 2$-contraction by $F$ and
$\pi_{\real{\lambda^{2,0}}}$ is the projection $\Lambda^2 T^* M
\to \real{\lambda^{2,0}}$. Thus, from \eqref{ec1bis} and
\eqref{ec2} we obtain
\begin{equation} \label{ec1bisbisbis}
{\small  \begin{array}{rl} & \hspace{-.8cm}  \tfrac12 F_{1\,2} (d^2
F) (X,Y) =
  \\[1mm]
& - 2 \langle ( \nabla^{\Lie{U}(n)}_{e_i} \xi )_{\varphi e_i} X,
\varphi Y \rangle
 + 2 \langle ( \nabla^{\Lie{U}(n)}_{ \varphi e_i} \xi )_X Y,   { \varphi e_i} \rangle
  - 2 \langle ( \nabla^{\Lie{U}(n)}_{\varphi  e_i} \xi )_Y X, \varphi  e_i \rangle
  + 2 \langle ( \nabla^{\Lie{U}(n)}_X \xi )_{ \varphi e_i} {\varphi e_i},  Y \rangle
\\[1mm]
  &
  - 2 \langle ( \nabla^{\Lie{U}(n)}_Y \xi )_{\varphi e_i} { \varphi e_i},  X \rangle
  - 2  \langle \xi_{ \xi_{e_i} {\varphi e_i}} X , \varphi Y \rangle
  + 2 \langle \xi_X e_i  , \xi_{e_i} Y , \rangle
  - 2  \langle  \xi_Y { e_i} , \xi_{e_i} X  \rangle
\\[1mm]
    &
    - \eta(X) \left( (\nabla^{\Lie{U}(n)}_Y \xi )_{e_i} \eta \right) (e_i)
    + \eta(Y) \left( (\nabla^{\Lie{U}(n)}_X \xi )_{e_i} \eta \right) (e_i)
 - 2 \eta \odot \left( (\nabla^{\Lie{U}(n)}_{e_i} \xi )_{\varphi e_i} \eta \right) \circ \varphi (X,Y)
\\[1mm]
&
 +  \eta(X) \left( (\nabla^{\Lie{U}(n)}_{e_i} \xi )_{Y} \eta \right) (e_i)
 +  \eta(Y) \left( (\nabla^{\Lie{U}(n)}_{e_i} \xi )_{X} \eta \right) (e_i)
 +  2 \eta(X) (\xi_{\xi_{e_i} Y} \eta) (e_i)
 \\[1mm]
 &
 + 2 (\xi_{e_i} \eta)(X) (\xi_{\varphi e_i} \eta)\circ \varphi (Y)
 - 2 (\xi_{e_i} (\xi_Y \eta))  ({e_i}) \eta(X)
 + 2 (\xi_X \eta)(Y) (\xi_{ e_i} \eta)({e_i})
 \\[1mm]
 &
 + 2 (\xi_Y \eta)({e_i}) (\xi_{\varphi e_i} \eta)\circ \varphi (X)
 + 2 (\xi_{e_i} (\xi_{\varphi e_i} \eta)) \circ \varphi (X) \eta(Y)
 -  (\xi_{ e_i} \eta)({\varphi e_i}) (\xi_X \eta)\circ \varphi (Y)
 \\[1mm]
 &
 +  (\xi_{ e_i} \eta)({\varphi e_i}) (\xi_Y \eta)\circ \varphi (X)
   -  3 (\xi_{\xi_X \zeta} \eta)(Y)
 +  3 (\xi_{\xi_Y \zeta} \eta)(X).
\end{array} }
\end{equation}
Therefore, $\pi_{\real{\lambda^{2,0}}} \circ F_{1\,2} ( d^2 F)$ is
such that
  {\small
      \begin{gather*}
      \begin{split}
      & \hspace{-.5cm} \tfrac12 \,   \pi_{\real{\lambda^{2,0}}}  \circ  \,
F_{1\,2}  (d^2 F)(X_{\zeta^{\perp}},Y_{\zeta^{\perp}}) = \\
 & - 2 \langle (\nabla^{\Lie{U}(n)}_{e_i} \xi )_{\varphi e_i} X_{\zeta^{\perp}}, \varphi Y_{\zeta^{\perp}} \rangle
  +\langle (\nabla^{\Lie{U}(n)}_{\varphi e_i} \xi )_{X_{\zeta^{\perp}}} Y_{\zeta^{\perp}}, \varphi e_i \rangle
  -\langle (\nabla^{\Lie{U}(n)}_{\varphi e_i} \xi )_{Y_{\zeta^{\perp}}} X_{\zeta^{\perp}},\varphi e_i \rangle
   \\
  & -\langle (\nabla^{\Lie{U}(n)}_{ \varphi e_i} \xi )_{\varphi X_{\zeta^{\perp}}} \varphi Y_{\zeta^{\perp}}, \varphi e_i \rangle
   +\langle (\nabla^{\Lie{U}(n)}_{e_i} \xi )_{\varphi Y_{\zeta^{\perp}}} \varphi X_{\zeta^{\perp}}, e_i \rangle
   + \langle (\nabla^{\Lie{U}(n)}_{X_{\zeta^{\perp}}} \xi )_{e_i} e_i, Y_{\zeta^{\perp}} \rangle
   \\
   &   - \langle (\nabla^{\Lie{U}(n)}_{Y_{\zeta^{\perp}}} \xi )_{e_i} e_i, X_{\zeta^{\perp}} \rangle
      -  \langle (\nabla^{\Lie{U}(n)}_{ \varphi X_{\zeta^{\perp}}} \xi )_{e_i} e_i, \varphi Y_{\zeta^{\perp}} \rangle
       + \langle (\nabla^{\Lie{U}(n)}_{\varphi Y_{\zeta^{\perp}}} \xi )_{e_i} e_i, \varphi X_{\zeta^{\perp}} \rangle
    \\
    &    - 2\langle \xi_{\xi_{e_i} \varphi e_i} X_{\zeta^{\perp}}, \varphi Y_{\zeta^{\perp}} \rangle
       + \langle  \xi_{X_{\zeta^{\perp}}} e_i,  \xi_{e_i} Y_{\zeta^{\perp}} \rangle
       - \langle  \xi_{Y_{\zeta^{\perp}}} e_i,  \xi_{e_i}  X_{\zeta^{\perp}} \rangle
        \\
    & - \langle  \xi_{\varphi X_{\zeta^{\perp}}} e_i,  \xi_{e_i} \varphi Y_{\zeta^{\perp}} \rangle
     + \langle \xi_{\varphi Y_{\zeta^{\perp}}} e_i,  \xi_{e_i}  \varphi X_{\zeta^{\perp}} \rangle
     +  (\xi_{e_i} \eta)(X_{\zeta^{\perp}}) (\xi_{\varphi e_i} \eta)\circ \varphi (Y_{\zeta^{\perp}})
     \\
    & - (\xi_{e_i} \eta)(Y_{\zeta^{\perp}}) (\xi_{\varphi e_i} \eta)\circ \varphi (X_{\zeta^{\perp}})
   + \tfrac12  (\xi_{X_{\zeta^{\perp}}} \eta)(Y_{\zeta^{\perp}}) (\xi_{ e_i} \eta)({e_i})
  - \tfrac12  (\xi_{Y_{\zeta^{\perp}}} \eta)(X_{\zeta^{\perp}}) (\xi_{ e_i} \eta)({e_i})
  \\
  &- \tfrac12 (\xi_{\varphi X_{\zeta^{\perp}}} \eta)(\varphi Y_{\zeta^{\perp}}) (\xi_{ e_i} \eta)({e_i})
   + \tfrac12  (\xi_{\varphi Y_{\zeta^{\perp}}} \eta)(\varphi X_{\zeta^{\perp}}) (\xi_{ e_i} \eta)({e_i})
  - \tfrac12 (\xi_{Y_{\zeta^{\perp}}} \eta)({\varphi  e_i}) (\xi_{e_i} \eta)\circ \varphi (X_{\zeta^{\perp}})
   \\
  &+ \tfrac12 (\xi_{X_{\zeta^{\perp}}} \eta)({\varphi e_i}) (\xi_{e_i} \eta)\circ \varphi (Y_{\zeta^{\perp}})
  + \tfrac12 (\xi_{\varphi Y_{\zeta^{\perp}}} \eta)({e_i})(\xi_{\varphi e_i} \eta) (X_{\zeta^{\perp}})
 - \tfrac12 (\xi_{\varphi X_{\zeta^{\perp}}} \eta)({e_i}) (\xi_{\varphi e_i} \eta)(Y_{\zeta^{\perp}})
  \\
   & - \tfrac12  (\xi_{ e_i} \eta)({\varphi e_i}) (\xi_{X_{\zeta^{\perp}}} \eta)\circ \varphi (Y_{\zeta^{\perp}})
     + \tfrac12 (\xi_{ e_i} \eta)({\varphi e_i}) (\xi_{Y_{\zeta^{\perp}}} \eta)\circ \varphi (X_{\zeta^{\perp}})
 - \tfrac12 (\xi_{ e_i} \eta)({\varphi e_i}) (\xi_{\varphi X_{\zeta^{\perp}}} \eta) (Y_{\zeta^{\perp}})
  \\
  &+ \tfrac12 (\xi_{ e_i} \eta)({\varphi e_i}) (\xi_{\varphi Y_{\zeta^{\perp}}} \eta) (X_{\zeta^{\perp}})
  - \tfrac32 (\xi_{ X_{\zeta^{\perp}}} \eta) (e_i) (\xi_{e_i} \eta)( Y_{\zeta^{\perp}})
 + \tfrac32 (\xi_{ Y_{\zeta^{\perp}}} \eta)(e_i) (\xi_{e_i} \eta)( X_{\zeta^{\perp}})
    \\
  &+ \tfrac32 (\xi_{\varphi X_{\zeta^{\perp}}} \eta) (e_i) (\xi_{e_i} \eta)(\varphi Y_{\zeta^{\perp}})
 - \tfrac32  (\xi_{\varphi Y_{\zeta^{\perp}}} \eta)(e_i) (\xi_{e_i} \eta)(\varphi
 X_{\zeta^{\perp}}).
  \end{split}
  \end{gather*}  }

\noindent From this expression,  taking into account the properties
of the components $\xi_{(1)}, \dots , \xi_{{(10)}}$ of $\xi$  (see
\cite{ChineaGonzalezDavila}), one can obtain
$\pi_{\real{\lambda^{2,0}}} \circ F_{1\,2} ( d^2 F)$ given in terms
of $\xi_{(1)}, \dots , \xi_{{(10)}}$. We recall that if
$\nabla^{\Lie{U}(n)}$ is a $\Un(n)$-connection, then
$\nabla^{\Lie{U}(n)} \xi_{(i)} \in \mathcal{C}_i$. Thus, we get
$\frac14 \pi_{\real{\lambda^{2,0}}} \circ F_{1\,2} (d^2 F)$ which is
the right side of next equation. Then we have
\begin{lemma} \label{xxx}
 For almost contact metric structures of type $\mathcal C_1 \oplus \ldots \oplus \mathcal C_{10}$, the
following identity is satisfied 

\begin{equation*}
{\rm  \small
     \begin{array}{rl}
     0= &
      3 \langle (\nabla^{\Lie{U}(n)}_{e_i} \xi_{(1)} )_{e_i} X_{\zeta^{\perp}}, Y_{\zeta^{\perp}} \rangle
  - \langle (\nabla^{\Lie{U}(n)}_{e_i} \xi_{(3)} )_{e_i} X_{\zeta^{\perp}}, Y_{\zeta^{\perp}} \rangle
  + (n-2) \langle (\nabla^{\Lie{U}(n)}_{e_i} \xi_{(4)} )_{e_i} X_{\zeta^{\perp}},Y_{\zeta^{\perp}}
   \rangle
  \\[2mm]
&
  +  \langle \xi_{{(3)}X_{\zeta^{\perp}}} e_i,  \xi_{{(1)}e_i} Y_{\zeta^{\perp}} \rangle
  -  \langle  \xi_{{(3)}Y_{\zeta^{\perp}}} e_i,  \xi_{{(1)}e_i}  X_{\zeta^{\perp}} \rangle
  +  \langle  \xi_{{(3)}X_{\zeta^{\perp}}} e_i,  \xi_{{(2)}e_i} Y_{\zeta^{\perp}} \rangle
   -  \langle  \xi_{{(3)}Y_{\zeta^{\perp}}} e_i,  \xi_{{(2)}e_i}  X_{\zeta^{\perp}} \rangle
    \\[2mm]
   &
 \displaystyle   - \tfrac{n-5}{n-1}  \langle \xi_{{(1)} \xi_{{(4)}e_i} e_i }X_{\zeta^{\perp}}, Y_{\zeta^{\perp}} \rangle
    - \tfrac{n-2}{n-1}  \langle  \xi_{{(2)} \xi_{{(4)}e_i} e_i } X_{\zeta^{\perp}}, Y_{\zeta^{\perp}} \rangle
      +  \langle \xi_{{(3)} \xi_{(4)e_i} e_i } X_{\zeta^{\perp}}, Y_{\zeta^{\perp}} \rangle
     \\[2mm]
     &
     + (n- 2)  (   \xi_{{(5)} \; e_i} \eta) \wedge (  \xi_{{(10)} \; e_i}
     \eta)(  X_{\zeta^{\perp}},  Y_{\zeta^{\perp}})
      - 2  (   \xi_{{(8)} \; e_i} \eta) \wedge (  \xi_{{(10)} \; e_i}
     \eta)(  X_{\zeta^{\perp}},  Y_{\zeta^{\perp}})
        \\[2mm]
     &
     +(n - 2)  (   \xi_{{(6)} \; e_i} \eta) \wedge (  \xi_{{(10)} \; e_i}
     \eta)(  X_{\zeta^{\perp}},  Y_{\zeta^{\perp}})
      - 2  (   \xi_{{(7)} \; e_i} \eta) \wedge (  \xi_{{(10)} \; e_i}
     \eta)(  X_{\zeta^{\perp}},  Y_{\zeta^{\perp}}).
   \end{array}}
  \end{equation*}
  \end{lemma}

In order to compute the $\real{\lambda^{1,0}}$-part of $d^2 F$, we
replace $X=\zeta$ and $Y = Y_{\zeta_{\perp}}$ in
\eqref{ec1bisbisbis}. Thus, making extensive use of  the
properties of $\xi_{(1)}, \dots , \xi_{{(10)}}$, we get
\begin{equation} \label{ec1bisbisbis1}
{\footnotesize  \begin{array}{rl}
 \tfrac14   F_{1\,2}  (d^2F)& \hspace{-.3cm}  (
 \zeta,Y_{\zeta_{\perp}})=
  \\[1mm]
&
      \langle ( \nabla^{\Lie{U}(n)}_{\zeta} \xi_{(4)} )_{  e_i} { e_i},  Y_{\zeta_{\perp}} \rangle
    + (n - 2) \left( (\nabla^{\Lie{U}(n)}_{e_i} \xi_{(5)} )_{e_i} \eta \right)  (Y_{\zeta_{\perp}})
    -  \left( (\nabla^{\Lie{U}(n)}_{e_i} \xi_{(6)} )_{e_i} \eta \right)  (Y_{\zeta_{\perp}})
\\[1mm]
  &
    -  \left( (\nabla^{\Lie{U}(n)}_{e_i} \xi_{(7)} )_{e_i} \eta \right) (Y_{\zeta_{\perp}})
    - 2 \left( (\nabla^{\Lie{U}(n)}_{e_i} \xi_{(8)} )_{ e_i} \eta \right)  (Y_{\zeta_{\perp}})
    +  \left( (\nabla^{\Lie{U}(n)}_{e_i} \xi_{(9)} )_{ e_i} \eta \right) (Y_{\zeta_{\perp}})
\\[1mm]
  &
    + 2 \left( (\nabla^{\Lie{U}(n)}_{e_i} \xi_{{(10)}} )_{e_i} \eta \right)  (Y_{\zeta_{\perp}})
    + 2  (\xi_{{(9)}e_i} \eta) (\xi_{{(1)}e_i} Y_{\zeta_{\perp}})
    +    (\xi_{{(9)}e_i} \eta) (\xi_{{(2)}e_i} Y_{\zeta_{\perp}})
    -  3  (\xi_{{(10)}e_i} \eta) (\xi_{{(2)}e_i} Y_{\zeta_{\perp}})
    \\[1mm]
  &
   +    (\xi_{{(5)}e_i} \eta) (\xi_{{(3)}e_i} Y_{\zeta_{\perp}})
   -    (\xi_{{(6)}e_i} \eta) (\xi_{{(3)}e_i} Y_{\zeta_{\perp}})
    -    (\xi_{{(7)}e_i} \eta) (\xi_{{(3)}e_i} Y_{\zeta_{\perp}})
      +    (\xi_{{(8)}e_i} \eta) (\xi_{{(3)}e_i} Y_{\zeta_{\perp}})
         \\[1mm]
         &
  -   (\xi_{{(9)}e_i} \eta) (\xi_{{(3)}Y} { e_i})
  -   (\xi_{{(10)}e_i} \eta) (\xi_{{(3)}Y} { e_i})
    +   \tfrac{1}{2n}   (\xi_{{(6)}e_i} \eta)(\varphi e_i)     \langle \xi_{{(4)}e_i} e_i ,
    \varphi Y_{\zeta_{\perp}} \rangle
\\[1mm]
    &
    - \tfrac{1}{n-1}   (\xi_{{(7)}Y_{\zeta_{\perp}}} \eta) (\xi_{{(4)}e_i} e_i)
    + \tfrac{n}{n-1}   (\xi_{{(8)}Y} \eta) (\xi_{{(4)}e_i} e_i)
   + \tfrac{n-3}{n-1}  (\xi_{{(10)}Y_{\zeta_{\perp}}} \eta) (\xi_{{(4)}e_i} {  e_i}).
\end{array} }
\end{equation}

 On the other hand, by replacing  $X = Y_{\zeta_{\perp}}$
and $Y=\zeta$, we will obtain
\begin{equation} \label{ec1bisbisbis2}
{\footnotesize  \begin{array}{rl}
 \tfrac14   F_{1\,2}  (d^2
F)& \hspace{-.3cm}  (Y_{\zeta_{\perp}},\zeta) =
  \\[1mm]
&
  -  \langle ( \nabla^{\Lie{U}(n)}_{\zeta} \xi_{(4)} )_{ e_i} { e_i}, Y_{\zeta_{\perp}} \rangle
 -(n-1)   \left( (\nabla^{\Lie{U}(n)}_{e_i} \xi_{(5)} )_{e_i} \eta \right)  (Y_{\zeta_{\perp}})
  - 2  \left( (\nabla^{\Lie{U}(n)}_{e_i} \xi_{(6)} )_{e_i} \eta \right)  (Y_{\zeta_{\perp}})
\\[1mm]
&
 - 2 \left( (\nabla^{\Lie{U}(n)}_{e_i} \xi_{(7)} )_{e_i} \eta \right)  (Y_{\zeta_{\perp}})
 +   \left( (\nabla^{\Lie{U}(n)}_{e_i} \xi_{(8)} )_{e_i} \eta \right)  (Y_{\zeta_{\perp}})
  + 2 (( \nabla^{\Lie{U}(n)}_{ e_i} \xi_{(9)} )_{X} \eta)   {  e_i}
\\[1mm]
&
 -  \left( (\nabla^{\Lie{U}(n)}_{e_i} \xi_{{(10)}})_{e_i} \eta \right) (Y_{\zeta_{\perp}})
 -   (\xi_{{(6)}e_i} \eta)(\xi_{{(1)} e_i} (Y_{\zeta_{\perp}}))
  -   (\xi_{{(7)} e_i} \eta)(\xi_{{(1)} e_i}  (Y_{\zeta_{\perp}}))
\\[1mm]
&
  - 2  (\xi_{{(9)} e_i} \eta)  ( \xi_{{(1)} e_i} Y_{\zeta_{\perp}})
  -   (\xi_{{(10)}e_i} \eta)  ( \xi_{{(1)}e_i} Y_{\zeta_{\perp}})
   -   (\xi_{{(6)} e_i} \eta)(\xi_{{(2)} e_i}  (Y_{\zeta_{\perp}}))
  -   (\xi_{{(7)} e_i} \eta)(\xi_{{(2)} e_i}  (Y_{\zeta_{\perp}}))
\\[1mm]
&
  -   (\xi_{{(9)} e_i} \eta)(\xi_{{(2)}e_i}  (Y_{\zeta_{\perp}}))
   + 2  (\xi_{{(10)}e_i} \eta)  ( \xi_{{(2)}e_i}Y_{\zeta_{\perp}} )
   +   (\xi_{{(6)} e_i} \eta)(\xi_{{(3)} e_i} Y_{\zeta_{\perp}})
  +   (\xi_{{(7)} e_i} \eta)(\xi_{{(3)} e_i} Y_{\zeta_{\perp}})
  \\[1mm]
  &
    +   (\xi_{{(9)}e_i} \eta)  ( \xi_{{(3)}Y_{\zeta_{\perp}}} e_i)
   +   (\xi_{{(10)}e_i} \eta)  ( \xi_{{(3)}Y_{\zeta_{\perp}}} e_i)
  +   (\xi_{{(9)}e_i} \eta)(\xi_{{(3)}e_i} Y_{\zeta_{\perp}})
 \\[1mm]
&
    - \tfrac{1}{2n}   (\xi_{{(6)} e_i} \eta) (\varphi e_i)  \langle \xi_{{(4)}e_i}
  e_i, \varphi Y_{\zeta_{\perp}} \rangle
  + \tfrac{1}{n-1}   (\xi_{{(7)}Y_{\zeta_{\perp}}} \eta)  ( \xi_{{(4)}e_i} e_i)
    + \tfrac{2}{n-1}   (\xi_{{(10)}Y_{\zeta_{\perp}}} \eta)  ( \xi_{{(4)}e_i} e_i).
\end{array} }
\end{equation}

Since $d^2 F=0$, from these two identities  the last  two
identities included in next Lemma are immediate consequences. The
first identity in the Lemma is obtained by summing the equations
\eqref{ec1bisbisbis1} and \eqref{ec1bisbisbis2}.

\begin{lemma} \label{ec1bisbisbisxx}
 For almost contact metric manifolds of type
$\mathcal C_1 \oplus \ldots \oplus \mathcal C_{10}$, the following
identity is satisfied
  \begin{equation*}
{\rm
 \begin{array}{rl}
 0=&
    -  \left( (\nabla^{\Lie{U}(n)}_{e_i} \xi_{(5)} )_{e_i} \eta \right)  (Y_{\zeta_{\perp}})
    - 3 \left( (\nabla^{\Lie{U}(n)}_{e_i} \xi_{(6)} )_{e_i} \eta \right)  (Y_{\zeta_{\perp}})
    - 3 \left( (\nabla^{\Lie{U}(n)}_{e_i} \xi_{(7)} )_{e_i} \eta \right) (Y_{\zeta_{\perp}})
    \\[1mm]
  &
    -  \left( (\nabla^{\Lie{U}(n)}_{e_i} \xi_{(8)} )_{ e_i} \eta \right)  (Y_{\zeta_{\perp}})
    + 3 \left( (\nabla^{\Lie{U}(n)}_{e_i} \xi_{(9)} )_{ e_i} \eta \right) (Y_{\zeta_{\perp}})
    +  \left( (\nabla^{\Lie{U}(n)}_{e_i} \xi_{{(10)}} )_{e_i} \eta \right)  (Y_{\zeta_{\perp}})
\\[2mm]
    &
    -   (\xi_{{(6)}e_i} \eta) (\xi_{{(1)}e_i} Y_{\zeta_{\perp}})
     -   (\xi_{{(7)}e_i} \eta) (\xi_{{(1)}e_i} Y_{\zeta_{\perp}})
      -  (\xi_{{(10)}e_i} \eta) (\xi_{{(1)}e_i} Y_{\zeta_{\perp}})
\\[2mm]
    &
    -    (\xi_{{(6)}e_i} \eta) (\xi_{{(2)}e_i} Y_{\zeta_{\perp}})
    -    (\xi_{{(7)}e_i} \eta) (\xi_{{(2)}e_i} Y_{\zeta_{\perp}})
    +    (\xi_{{(10)}e_i} \eta) (\xi_{{(2)}e_i} Y_{\zeta_{\perp}})
\\[2mm]
    &
    +    (\xi_{{(5)}e_i} \eta) (\xi_{{(3)}e_i} Y_{\zeta_{\perp}})
      +    (\xi_{{(8)}e_i} \eta) (\xi_{{(3)}e_i} Y_{\zeta_{\perp}})
    +    (\xi_{{(9)}e_i} \eta) (\xi_{{(3)}e_i} Y_{\zeta_{\perp}})
       \\[1mm]
  &
     + \displaystyle \tfrac{n}{n-1}   (\xi_{{(8)}\xi_{{(4)}e_i} e_i} \eta) (Y_{\zeta_{\perp}})
    -  \displaystyle   (\xi_{{(10)}\xi_{{(4)}e_i} { e_i}} \eta) (Y_{\zeta_{\perp}}).
\end{array} }
\end{equation*}
In particular, if the almost contact metric manifold is of  type
$\mathcal C_5 \oplus \dots   \oplus \mathcal C_{10} $, then
\begin{eqnarray*}
0 & = & (n-1) (\nabla^{\Lie{U}(n)}_{e_i} \xi_{(5)} )_{e_i} +  2
(\nabla^{\Lie{U}(n)}_{e_i} \xi_{(6)} )_{e_i}  + 2
(\nabla^{\Lie{U}(n)}_{e_i} \xi_{(7)} )_{ e_i} \\
&&
 -(\nabla^{\Lie{U}(n)}_{e_i} \xi_{(8)} )_{ e_i}
 - 2(\nabla^{\Lie{U}(n)}_{e_i} \xi_{(9)} )_{ e_i}
 +(\nabla^{\Lie{U}(n)}_{e_i} \xi_{(10)} )_{ e_i}, \\
0 & = &  (n-2) (\nabla^{\Lie{U}(n)}_{e_i} \xi_{(5)} )_{e_i}
   - (\nabla^{\Lie{U}(n)}_{e_i} \xi_{(6)} )_{e_i}
   - (\nabla^{\Lie{U}(n)}_{e_i} \xi_{(7)} )_{ e_i} \\
   &&
   - 2 (\nabla^{\Lie{U}(n)}_{e_i} \xi_{(8)} )_{ e_i}
 +  (\nabla^{\Lie{U}(n)}_{e_i} \xi_{(9)} )_{ e_i}
 + 2(\nabla^{\Lie{U}(n)}_{e_i} \xi_{(10)} )_{ e_i}.
\end{eqnarray*}
\end{lemma}

In the following proofs below we will make extensive use of the
identities contained in Lemmas \ref{xxx} and \ref{ec1bisbisbisxx}.
Likewise, we will need to use the fact ${\xi_{{(4)} \, \xi_{e_i}
e_i}} = 0$ satisfied by almost contact metric structures of type
$\mathcal{C}_1 \oplus \dots \oplus \mathcal C_{10}$.

\begin{proof}[Proof of Theorem \ref{harmclasscontact}]
For (i). The tensor $\xi$ for almost contact metric structures of
type $ \mathcal C_1 \oplus \mathcal C_2 \oplus \mathcal C_5 \oplus
\mathcal C_{6}\oplus  \mathcal C_7 \oplus \mathcal C_{8}$ is such
that   $\xi_{\zeta} =0$, $\xi_{\varphi X} \varphi Y = -
\xi_{X_{\zeta^{\perp}}} Y_{\zeta^{\perp}}$ and $(\xi_{\varphi X}
\eta )(\varphi Y) = (\xi_{X}\eta)(Y)$ \cite{ChineaGonzalezDavila}.
Also in such a case, we have $\xi_{\xi_{e_i} \varphi e_i} \eta =0$.
Therefore, by Lemma \ref{astricciacm1} and  Theorem
\ref{characharmalmcontact} (vii), (i) follows.

 The proof for $\mathcal D=\mathcal
C_1 \oplus \mathcal C_2 \oplus  \mathcal C_9 \oplus  \mathcal
C_{10}$ is similar. In this case
 we have $\xi_{\zeta} =0$,  $\xi_{\varphi X} \varphi Y = -
\xi_{X_{\zeta^{\perp}}} Y_{\zeta^{\perp}}$, $(\xi_{\varphi X} \eta
)(\varphi Y)  = - (\xi_{X}\eta)(Y)$ and $\xi_{\xi_{e_i} \varphi
e_i} \eta=0$.

 For (ii). Because the structure is of type $\mathcal C_1 \oplus  \mathcal C_4 \oplus \mathcal C_5
\oplus \mathcal C_6 \oplus \mathcal C_7 \oplus \mathcal C_8$, the
identities in Lemmas \ref{astricciacm1} and  \ref{xxx} are given
respectively  by
\begin{eqnarray}
  \label{d2omega:part20c1c4}
      0   &  =& \;
  3 \langle (\nabla^{\Lie{U}(n)}_{e_i}\xi_{(1)})_{e_i} X_{\zeta^{\perp}} , Y_{\zeta^{\perp}}\rangle
  +     (n-2)
    \langle (\nabla^{\Lie{U}(n)}_{e_i}\xi_{(4)})_{e_i}X_{\zeta^{\perp}} ,
    Y_{\zeta^{\perp}}\rangle
    \\
    & &
      - \tfrac{n-5}{n-1} \langle \xi_{{(1)} \xi_{{(4)}e_i} e_i} X_{\zeta^{\perp}} ,
      Y_{\zeta^{\perp}} \rangle,
  \nonumber \\
    \label{otraricshc1c4}
      \Ricac_{\mbox{\rm\footnotesize alt}} (X_{\zeta^{\perp}},Y_{\zeta^{\perp}})
       & = &  \langle \xi_{{(1)}\xi_{{(4)}e_i} e_i}X_{\zeta^{\perp}} ,
      Y_{\zeta^{\perp}}\rangle
      - \langle (\nabla^{\Lie{U}(n)}_{e_i}\xi_{(2)})_{e_i}X_{\zeta^{\perp}} , Y_{\zeta^{\perp}} \rangle
      \\
       &&  + \langle (\nabla^{\Lie{U}(n)}_{e_i}\xi_{(4)})_{e_i}X_{\zeta^{\perp}} , Y_{\zeta^{\perp}} \rangle.
\nonumber
\end{eqnarray}

Likewise, the first one of the  characterising  conditions for
harmonic almost contact  structures given in Theorem
\ref{astricciacm1}~(vii) is expressed by
\begin{equation} \label{characharmhermc1c4}
- \langle \xi_{{(1)} \xi_{{(4)}e_i} e_i} X_{\zeta^{\perp}} ,
Y_{\zeta^{\perp}} \rangle = \langle ( \nabla^{\Lie{U}(n)}_{e_i}
 \xi_{(1)})_{e_i} X_{\zeta^{\perp}} , Y_{\zeta^{\perp}} \rangle  + \langle (
\nabla^{\Lie{U}(n)}_{e_i}
 \xi_{(4)})_{e_i} X_{\zeta^{\perp}} , Y_{\zeta^{\perp}} \rangle .
\end{equation}
Now, for $n \geq 3$, it is straightforward to check that equations
\eqref{d2omega:part20c1c4},   \eqref{otraricshc1c4} and
\eqref{characharmhermc1c4} imply the  expression for
$\Ricac_{\mbox{\rm \footnotesize alt}}$ required in (ii).

On the other hand, for almost contact metric  structures of type
$\mathcal C_1 \oplus  \mathcal C_4 \oplus \mathcal C_5 \oplus
\mathcal C_6 \oplus \mathcal C_7 \oplus \mathcal C_8$, we have $(
\xi_{\varphi X} \eta )(\varphi Y)  = (\xi_{X}\eta)(Y)$. Moreover, $
\xi_{\xi_{e_i} \varphi e_i} \eta = \xi_{\xi_{{(4)}e_i} \varphi e_i}
\eta = - \xi_{ \varphi \xi_{e_i}  e_i} \eta$. Therefore,
$$
 \Ricac( \zeta , X ) =
 ( (\nabla^{\Lie{U}(n)}_{e_i} \xi)_{ e_i} \eta )( X)
  -  ( \xi_{\xi_{e_i}  e_i} \eta )( X).
$$
Now, from the second condition in Theorem
\ref{characharmalmcontact} (vii), it follows $\Ricac( \zeta , X )
=  -  2 ( \xi_{\xi_{e_i}  e_i} \eta )( X)$.

 Conversely, it is also direct to see that the expressions
for $\Ricac$ in (ii) and equations \eqref{d2omega:part20c1c4} and
 \eqref{otraricshc1c4}
 imply equation \eqref{characharmhermc1c4} and $(\nabla^{\Lie{U}(n)}_{e_i} \xi)_{ e_i} \eta +
    \xi_{\xi_{e_i}  e_i} \eta=0$. Therefore, the almost
 contact metric structure is harmonic.

  The proof for (iii) is similar to the one for (ii). The only difference is that in this case
 we have  $(\xi_{\varphi X} \eta
)(\varphi Y)  = - (\xi_{X}\eta)(Y)$.

For (iv). If the structure is of type $\mathcal C_2 \oplus
\mathcal C_4 \oplus \mathcal C_5 \oplus \mathcal C_6 \oplus
\mathcal C_7 \oplus \mathcal C_8$, then the identities in Lemmas
\ref{astricciacm1} and  \ref{xxx} are given respectively  by
\begin{equation}
  \label{d2omega:part20c2c4}
  \begin{array}{rl}
    0
    =& \displaystyle
    (n-2)   \langle (\nabla^{\Lie{U}(n)}_{e_i}\xi_{(4)})_{e_i}X_{\zeta^{\perp}},
    Y_{\zeta^{\perp}} \rangle
      - \tfrac{n-2}{n-1}\langle \xi_{{(2)}  \xi_{{(4)}e_i}
      e_i}X_{\zeta^{\perp}}
      ,Y_{\zeta^{\perp}} \rangle
   ,
  \end{array}
\end{equation}
\begin{equation}
    \label{otraricshc2c4}
    \begin{array}{rl}
      \Ricac_{\mbox{\rm\footnotesize alt}} (X_{\zeta^{\perp}},Y_{\zeta^{\perp}})
      = & \langle \xi_{{(2)}\xi_{{(4)}e_i} e_i}X_{\zeta^{\perp}} ,Y_{\zeta^{\perp}}
      - \langle (\nabla^{\Lie{U}(n)}_{e_i} \xi_{(2)})_{e_i}X_{\zeta^{\perp}} ,
      Y_{\zeta^{\perp}}\rangle
      \\
      &
      + \langle (\nabla^{\Lie{U}(n)}_{e_i}\xi_{(4)})_{e_i}X_{\zeta^{\perp}} , Y_{\zeta^{\perp}} \rangle       .
    \end{array}
\end{equation}
For this case, the first  condition  given in Theorem
\ref{characharmalmcontact} (vii) is expressed by
\begin{equation} \label{characharmhermc2c4}
- \langle \xi_{{(2)}\xi_{{(4)}e_i} e_i} X_{\zeta^{\perp}} ,
Y_{\zeta^{\perp}} \rangle = \langle ( \nabla^{\Lie{U}(n)}_{e_i}
 \xi_{(2)})_{e_i} X_{\zeta^{\perp}} , Y_{\zeta^{\perp}} \rangle  + \langle (
\nabla^{\Lie{U}(n)}_{e_i}
 \xi_{(4)})_{e_i} X_{\zeta^{\perp}} , Y_{\zeta^{\perp}} \rangle .
\end{equation}
Now, for  $n \geq 3$, it is straightforward to check that
equations \eqref{d2omega:part20c2c4},   \eqref{otraricshc2c4} and
 \eqref{characharmhermc2c4} imply the first required identity in
(iv). The expression for $\Ricac( \zeta , X )$ follows as in the
proof for  (ii).

 Conversely, it is also direct to see that the expressions
for $\Ricac$ in (iv) and equations \eqref{d2omega:part20c2c4} and
 \eqref{otraricshc2c4}
 imply equation \eqref{characharmhermc2c4}  and $(\nabla^{\Lie{U}(n)}_{e_i} \xi)_{ e_i} \eta +
    \xi_{\xi_{e_i}  e_i} \eta=0$. Therefore, the almost
 contact metric structure is harmonic.

 The proof for (v) is similar to the one for (iv). Note that in this case
 $(\xi_{\varphi X} \eta )(\varphi Y)  = - (\xi_{X}\eta)(Y)$.

For (vi). The intrinsic torsion $\xi$ for normal  structures is
such that  $\xi_{\zeta}=0$, $\xi_{\varphi X} \varphi Y =
\xi_{X_{\zeta^{\perp}}} Y_{\zeta^{\perp}}$ and  $(\xi_{\varphi X}
\eta) (\varphi Y) = (\xi_X \eta)(Y)$ (see
\cite{ChineaGonzalezDavila}). Therefore, the  required identities
in (vi) are  immediate consequences of Lemma \ref{astricciacm1}
and Theorem \ref{characharmalmcontact} (vii).

The proof for (vii) is similar to the one for (vi). In this case
 we have $\xi_{\zeta} =0$,  $\xi_{\varphi X} \varphi Y =
\xi_{X_{\zeta^{\perp}}} Y_{\zeta^{\perp}}$ and $(\xi_{\varphi X}
\eta )(\varphi Y)  =  (\xi_{X}\eta)(Y)$.

Finally, (viii), (ix), (i)$^*$ and (ii)$^*$ are immediate
consequences of Theorem \ref{characharmalmcontact} (vii),  Lemma
\ref{xxx} and Lemma \ref{ec1bisbisbisxx}.
\end{proof}

\begin{proof}[Proof of Corollary \ref{corharmonic}]
The assertions are immediate consequences of   Lemma
\ref{astricciacm1},  Theorem \ref{characharmalmcontact} and
Theorem \ref{harmclasscontact}.
\end{proof}

Now, we focus our attention on studying harmonicity  of almost
contact metric  structures  as a map into $\mathcal{SO}(M) /
(\Lie{U}(n) \times 1)$. Results in that direction were already
obtained by Vergara-Díaz and Wood \cite{VergaraWood}. We will
complete such results by using tools here presented. In next Lemma,
$s^{\rm ac}$ will denote the {\it $ac$-scalar  curvature} defined by
$s^{\rm ac} = \Ricac (e_i, e_i)$. If $\Ricac (X,Y) = \frac{1}{2n}
s^{\rm ac} (\langle X,Y \rangle - \eta (X) \eta (Y))$, the almost
contact metric manifold is said to be {\it weakly-ac-Einstein}. If
$s^{\rm ac}$ is constant, a weakly-ac-Einstein is called {\it
ac-Einstein}.

In Riemannian geometry, it is satisfied $2 d^* \Ric + d s =0$,
 where $s$ is the scalar curvature. The ac-analogue in almost
 contact metric geometry does not hold in general.
\begin{lemma} \label{id:genera}
 For almost contact metric manifolds of type
$\mathcal C_1 \oplus \ldots \oplus \mathcal C_{10}$, we have
\begin{equation*}
\begin{split}
 2d^* (\Ricac)^t(X) +  ds^{\rm ac} (X)  & =
        2 \langle R_{e_i , X} ,    \xi_{\varphi e_i} \varphi\rangle
        -  4 \Ricac ( X,   \xi_{e_i} e_i ) \\
       & \qquad   + 4 \langle \Ricac , \xi^\flat_X \rangle
        - 2 d^* F (\zeta) \Ricac (\zeta , \varphi X),
              \end{split}
\end{equation*}
 where $(\Ricac)^t(X,Y) = \Ricac(Y,X)$ and $\xi_X^{\flat} (Y,Z) =
\langle \xi_X Y, Z \rangle$. In particular, if the manifold is
weakly-ac-Einstein, then
$$
 (n-1) d s^{\rm ac}(X) + (d s^{\rm ac}(\zeta) + s^{\rm ac} d^* \eta) \eta(X) =    2n \langle R_{( e_i, X)} ,    \xi_{\varphi e_i}
 \varphi
 \rangle - 2 s^{\rm ac} \langle \xi_{e_i} e_i , X \rangle .
$$
\end{lemma}

\begin{proof}
Note that $ (\Ricac)^t(X,Y) = \tfrac12  \langle R_{e_i,\varphi
e_i} Y,  \varphi X\rangle$.  Then, we get
\begin{eqnarray*}
d^* (\Ricac)^t(X) & = & - (\nabla_{e_j} (\Ricac)^t ( e_j ,X)
\\
&=& -\tfrac12  e_j \langle R_{e_i, \varphi e_i} X ,\varphi e_j
\rangle
   + \tfrac12 \langle R_{e_i, \varphi e_i}  \nabla_{e_j} X , \varphi e_j \rangle
   + \tfrac12 \langle R_{e_i,\varphi e_i}   X,   \varphi \nabla_{e_j} e_j \rangle
   \\
   &=& - \tfrac12  \langle ( \nabla_{e_j} R)_{e_i,\varphi  e_i} X ,\varphi e_j \rangle
      -  \langle R_{\nabla_{e_j} e_i, \varphi e_i}   X , \varphi e_j \rangle
        - \tfrac12 \langle R_{e_i,\varphi e_i}   X,   (\nabla_{e_j}\varphi) e_j
        \rangle.
\end{eqnarray*}
Now,  since  $(\nabla_{e_j} \varphi)(e_j) = 2 \varphi \xi_{e_j} e_j
+ d^* F (\zeta) \zeta$, using symmetry properties of $R$, it follows
that
\begin{eqnarray*}
d^* (\Ricac)^t(X) & = & - \tfrac12  \langle ( \nabla_{e_j} R)_{X ,
\varphi e_j}  e_i, \varphi e_i\rangle
      +  \langle R_{X , e_j} e_i ,    \nabla_{\varphi e_j} \varphi e_i) \rangle
        -  \langle R_{e_i,\varphi e_i}   X,  \varphi \xi_{e_j} e_j \rangle\\
        && - \tfrac12 d^* F (\zeta) \langle R_{(e_i,\varphi e_i)}   X,  \zeta
        \rangle.
\end{eqnarray*}
Using second Bianchi's identity and taking into account
$$
 \langle R_{X , e_j} e_i ,    \nabla_{\varphi e_j} \varphi e_i) \rangle =
   \langle R_{X , e_j} e_i ,    \nabla^{\Lie{U}(n)}_{\varphi e_j} \varphi e_i) \rangle
 -  \langle R_{X , e_j} e_i ,    \xi_{\varphi e_j} \varphi e_i)
 \rangle,
$$
we get
\begin{eqnarray*}
d^* (\Ricac)^t (X) & = & - \tfrac14  \langle ( \nabla_{X} R)_{e_j,
\varphi e_j }  e_i,\varphi e_i\rangle
      -  \langle R_{X , e_j} ,    \xi_{\varphi e_j} \varphi
      \rangle
      \\
      &&
        -  2 \Ricac ( X,   \xi_{e_j} e_j )- d^* F (\zeta) \Ricac (\zeta , \varphi X) .
\end{eqnarray*}
Note that
$$
 \langle R_{X , e_j} e_i ,    \nabla^{\Lie{U}(n)}_{\varphi e_j} \varphi e_i) \rangle =
  \langle R_{X , e_j} e_i , e_k \rangle
     \langle \nabla^{\Lie{U}(n)}_{\varphi e_j} \varphi e_i , e_k \rangle =0,
$$
because it is a scalar product of a skew-symmetric matrix by a
symmetric matrix.

Finally, it is obtained
\begin{eqnarray} \label{uno}
\qquad 2d^* (\Ricac)^t(X) & = & - \tfrac12  \langle ( \nabla_{X}
R)_{e_j, Je_j }  e_i,\varphi e_i\rangle
      - 2 \langle R_{X , e_j} ,    \xi_{\varphi e_j} \varphi
      \rangle \\
      && \nonumber
        -  4 \Ricac ( X,   \xi_{e_j} e_j ) - 2 d^* F (\zeta) \Ricac (\zeta , \varphi X)  .
\end{eqnarray}

In a  second instance, ones obtains  $ ds^{\rm ac} (X) = \tfrac12
X \langle R_{e_i,Je_i} e_j, Je_j \rangle$. Hence
\begin{eqnarray*}
 ds^{\rm ac} (X) & = & \tfrac12   \langle (\nabla_X R)_{e_i,\varphi e_i} e_j, \varphi e_j \rangle
  + 2 \langle R_{e_i, \varphi e_i} e_j, \nabla_X \varphi e_j \rangle.
\end{eqnarray*}
But we also have  that
\begin{eqnarray*}
 \langle R_{e_i, \varphi e_i} e_j, \nabla_X \varphi e_j \rangle & = &  \langle R_{e_i, \varphi e_i} e_j, e_k \rangle \langle
 \nabla^{\Lie{U}(n)}_X \varphi e_j , e_k \rangle - \langle R_{e_i,\varphi e_i} e_j, e_k \rangle \langle
 \xi_X \varphi e_j , e_k \rangle \\
 & = & \langle R_{e_i,\varphi e_i} e_j, \varphi \xi_X e_j \rangle + \eta ( \xi_X \varphi e_j )
 \langle R_{e_i,\varphi e_i} \zeta , e_j \rangle
\\
&=&  2 \Ricac (e_i , \xi_X e_i) - 2 \eta ( \xi_X \varphi e_j )
\Ricac  (\zeta , \varphi e_j)
  \\
&=&
 2 \Ricac  (e_i , \xi_X e_i) + 2    \Ricac  (\zeta , \xi_X \zeta )
\\
&=& 2 \langle \Ricac   , \xi_X  \rangle.
\end{eqnarray*}
Thus,
 \begin{equation} \label{dos}
  ds^{\rm ac}  (X)  =  \tfrac12   \langle (\nabla_X
R)_{e_i,\varphi e_i} e_j, \varphi e_j \rangle + 4 \langle \Ricac ,
\xi^\flat_X \rangle.
 \end{equation}
 From \eqref{uno} and \eqref{dos}, the required identity is
 obtained.
\end{proof}

\begin{theorem} \label{mapharmcontact}
 For an $2n+1$-dimensional  almost contact metric  manifold $(M,\langle \cdot
,\cdot \rangle, \varphi, \zeta)$,  we have:
\begin{enumerate}
 \item[{\rm (i)}] If $M$ is of type $\mathcal D$, where  $\mathcal D = \mathcal{C}_1 \oplus
\mathcal{C}_2 \oplus \mathcal{C}_5\oplus \mathcal{C}_6\oplus
\mathcal{C}_7\oplus \mathcal{C}_8$ or $\mathcal D = \mathcal{C}_1
\oplus \mathcal{C}_2 \oplus \mathcal{C}_9 \oplus
\mathcal{C}_{10}$, then the almost contact metric structure is a
harmonic map  if and only if it is a harmonic structure and $2d^*
\Ricac + ds^{\rm ac} =0$.
 In particular:
 \begin{enumerate}
  \item[{\rm (a)}]
 If the structure is of type  $\mathcal{C}_1 \oplus
\mathcal{C}_2 \oplus \mathcal{C}_5 \oplus \mathcal{C}_6\oplus
\mathcal{C}_7\oplus \mathcal{C}_8$  and the  manifold is
weakly-$ac$-Einstein, then the almost contact metric structure is a
harmonic map  if and only if $s^{\rm ac}$ satisfies
 $n \, d s^{\rm ac} = s^{\rm ac} \, d^* \eta \, \eta = - s^{\rm
ac} \, \xi_{e_i}^\flat e_i$.
 \item[{\rm (b)}]
  If the structure is of type  $\mathcal D = \mathcal{C}_1 \oplus
\mathcal{C}_2 \oplus \mathcal{C}_9 \oplus \mathcal{C}_{10},
\mathcal{C}_1 \oplus \mathcal{C}_2  \oplus \mathcal{C}_6\oplus
\mathcal{C}_7\oplus \mathcal{C}_8 $ and the manifold is
weakly-$ac$-Einstein, then the almost contact  metric structure is
a harmonic map if and only if $s^{\rm ac}$ is constant.
 \item[{\rm (c)}]
  If the manifold is nearly-K-cosymplectic $(\mathcal{C}_1)$, then almost
 contact metric structure is a harmonic map.
 \end{enumerate}

  \item[{\rm (ii)}] If $M$  is of type $\mathcal D$, where
$\mathcal D = \mathcal{C}_3 \oplus \mathcal{C}_4 \oplus
\mathcal{C}_5\oplus \mathcal{C}_6\oplus \mathcal{C}_7\oplus
\mathcal{C}_8$ or $\mathcal D = \mathcal{C}_3 \oplus \mathcal{C}_4
\oplus \mathcal{C}_9 \oplus \mathcal{C}_{10}$,
  then the almost contact metric  structure  is a harmonic map
 if and only if  it is a harmonic structure and
 \begin{eqnarray*}
2 d^* (\Ricac)^t (X) + ds^{\rm ac} (X) +  4 \Ricac ( X,
\xi_{e_j} e_j ) && \\
 - 4 \langle \Ricac  , \xi^{\flat}_X
\rangle + 2 d^* F (\zeta) \Ricac (\zeta , \varphi X)  &= &0,
\end{eqnarray*}
for all $X \in \mathfrak X (M)$. In particular,
  \begin{enumerate}
  \item[{\rm (a)$^*$}] If $\Ricac$ is symmetric,  then the almost contact
  metric structure is a harmonic map if and only if   $\xi_{\xi_{e_i} e_i}=0$  and
   $2 d^{*} \Ricac + ds^{\rm ac} +  4 \xi_{e_i} e_i  \lrcorner
  \Ricac = 0$. Furthermore, if the manifold is weakly-$ac$-Einstein, then the almost contact metric structure
  is a harmonic map if and only if $\xi_{\xi_{e_i} e_i}=0$ and
   $(n-1) ds^{\rm ac} + (d s^{\rm ac} (\zeta) - s^{\rm ac} d^* \eta) \eta   +  2 s^{\rm ac}  \xi^{\flat}_{e_i} e_i  = 0$.
   \item[{\rm (b)$^*$}] If the almost contact metric structure is of type $\mathcal{C}_3 \oplus \mathcal{C}_i$,
   $i=6,7,10$, then the structure is a
   harmonic map if and only if
   $$
   2 d^* \Ricac + ds^{\rm ac}=0.
   $$
   Furthermore, if the manifold is also weakly-$ac$-Einstein, then
   the almost contact metric structure is a
   harmonic map if and only if $s^{\rm ac}$ is constant.
    \item[{\rm (c)$^*$}] For $n\neq 2$, if the almost contact metric structure is of type
   $\mathcal{C}_4 \oplus \mathcal{C}_i$, $i=5,6,7,9$, then the structure is a
   harmonic map if and only if
   $$
   2 d^* \Ricac + ds^{\rm ac} +  4 \xi_{e_i} e_i  \lrcorner
  \Ricac=0.
  $$
  \end{enumerate}
\end{enumerate}
\end{theorem}
\begin{proof} Most of the  results contained in Theorem are immediate
consequences of Theorem \ref{harmclasscontact} and   Lemma
\ref{id:genera}.
 \end{proof}

\begin{remark} {\rm
If a
 nearly-K-cosymplectic structure   is flat, then it is cosymplectic, i.e.
 $\xi=0$. In fact, in \cite{VergaraWood} it is shown that
$[\lie{u}(n)^{\perp}_{|\zeta^{\perp}} ,
\lie{u}(n)^{\perp}_{|\zeta^{\perp}}] \subseteq \lie{u}(n)$ and, for
nearly-K-cosymplectic structures, $\xi_X \in
\lie{u}(n)^{\perp}_{|\zeta^{\perp}}$, for all $X \in \mathfrak X
(M)$. Therefore, by Proposition \ref{pro:skew} (i), the assertion
follows.}
 \end{remark}


\section{Almost contact metric structures with minimal energy}
{\indent} \setcounter{equation}{0}
 In this section we will apply the simple
and elegant argument used by Bör et al. in \cite{BHLS2}  to almost
contact metric manifolds. For oriented and  compact Riemannian
manifolds $(M,\langle \cdot, \cdot \rangle)$ of dimension $n$
equipped with a $G$-structure, $G \subseteq \SO (n)$ and a
differential $p$-form $\phi$ preserved by the action of $G$,  the
following Bochner type formula was deduced in \cite{BHLS1}
\begin{equation} \label{mainBoch}
\int_M \left( \tfrac{1}{p+1} \| d \phi \|^2 + p \| d^* \phi \|^2 -
\| \nabla \phi \|^2 \right) =  \int_M \langle \widetilde{\mathcal
R} \phi , \phi \rangle,
\end{equation}
where  the volume form $dv$ is omitted  for sake of simplicity
and the operator $\widetilde{\mathcal R}$ is defined as follows.
We firstly remind that, for a covariant $p$-tensor $\alpha$, we
denote
$$
{\rm alt} \,\alpha (x_1,\dots, x_p) = \textstyle \sum_{\tau} {\rm
sign}(\tau) \alpha (x_{\tau(1)}, \dots , x_{\tau(p)}),
$$
where  the sum is extended on the set of permutations $\tau$ of
$\{1, \dots, p\}$ and ${\rm sign}(\tau)$ denotes  the signature of
$\tau$. If $\phi$ is a skew-symmetric  $p$-form, we write
$$
(R^c \phi)(x_1, \dots , x_p) = (R_{x_1 , e_i} \phi)(e_i,x_2,\dots,
x_p).
$$
That is, the operator $R_{x,y}$ acting on the $p$-form $\phi$ in the
usual way and then doing the indicated contraction. Finally, it is
defined
$$
\widetilde{\mathcal R}(\phi) =  {\rm alt} (R^c \phi).
$$

\begin{remark}
{\rm The constant coefficients  in the identity are different in
\cite{BHLS1}. This is due to the different conventions followed by
us for the wedge product, the scalar product of $p$-forms,  etc.}
\end{remark}

In the context of almost contact metric geometry, since $\langle
\widetilde{\mathcal R} \, F , F \rangle = 2 (s- s^{\rm ac} - \Ric
(\zeta , \zeta))$,
 the
Bochner type formula for the fundamental two-form $F$ is given by
$$
\int_M  \left( \tfrac{1}{3} \| d F \|^2 + 2 \| d^* F \|^2 - \|
\nabla F \|^2 \right) = 2 \int_M (s- s^{\rm ac} -  \Ric (\zeta ,
\zeta)).
$$
Moreover, since $\langle \widetilde{\mathcal R} \eta , \eta
\rangle =  \Ric (\zeta , \zeta)$, for the characteristic one-form
$\eta$ we have
$$
\int_M  \left( \tfrac{1}{2} \| d \eta \|^2 +  \| d^* \eta \|^2 -
\| \nabla \eta \|^2 \right) =  \int_M    \Ric (\zeta , \zeta).
$$
\vspace{1mm}

From the expression for the intrinsic torsion given by equation
\eqref{torsion:xialmcont}, it follows that
$$
 4 \| \xi \|^2 = \| \nabla F \|^2 + 6  \| \nabla \eta \|^2.
$$
We recall that $\nabla F \in \mathcal{C}_1 \oplus \dots \oplus
\mathcal{C}_{12}$, $\nabla \eta \in \mathcal{C}_5 \oplus \dots
\oplus \mathcal{C}_{10} \oplus  \mathcal{C}_{12}$. Taking this into
account and the  properties of a tensor in the module
$\mathcal{C}_i$, $i=1,\dots, 12$, (see \cite{ChineaGonzalezDavila}),
it is not hard to deduce
  \begin{gather*}
 4 \| \xi_{(i)} \|^2 = \| (\nabla F)_{(i)} \|^2,  \quad \| (\nabla
 \eta)_{(i)}
 \|^2 =0, \qquad i=1,2,3,4,11; \\
 \| \xi_{(i)} \|^2 = \| (\nabla F)_{(i)} \|^2,  \quad 2 \|
(\nabla \eta)_{(i)}
 \|^2 = \| \xi_{(i)} \|^2, \qquad i=5,6,7,8,9,10,12.
  \end{gather*}

Moreover, if we consider the decomposition
$$
\Lambda^3 {\rm T}^* M = \real{\lambda^{3,0}} \oplus
\real{\lambda^{2,1}_0}
   \oplus \real{\lambda^{1,0}}  \wedge F \oplus   \mathbb{R} \,
    F \wedge \eta \oplus \real{\lambda^{1,1}_0} \wedge \eta \oplus  \real{\lambda^{2,0}}  \wedge
   \eta,
$$
and the $\Lie{U}(n)$-map $\mathrm{alt} \, : \,  {\rm T}^* M
\otimes \lie{u}(n)^\perp  \to \Lambda^3 T^* M$, the exterior
derivative $d F=\mathrm{alt} (\nabla F)$ is in $ \mathcal{C}^a_{1}
\oplus  \mathcal{C}^a_{3} \oplus \mathcal{C}^a_{4} \oplus
\mathcal{C}^a_{5} \oplus \mathcal{C}^a_{8} \oplus
\mathcal{C}^a_{10,11}$, where
\begin{gather*}
\mathcal{C}^a_{1} = \real{\lambda^{3,0}} = \mathrm{alt} (
\mathcal{C}_1), \quad \mathcal{C}^a_{3} = \real{\lambda^{2,1}_0} =
\mathrm{alt} (  \mathcal{C}_{3}), \quad \mathcal{C}^a_{4} =
\real{\lambda^{1,0}}
\wedge F = \mathrm{alt} (\mathcal{C}_{4}), \\
\mathcal{C}^a_{5} = \mathbb{R} \,  F \wedge \eta  = \mathrm{alt} (
\mathcal{C}_5), \quad \mathcal{C}^a_{8} = \real{\lambda^{1,1}_0}
\wedge \eta = \mathrm{alt} (\mathcal{C}_{8}), \quad
\mathcal{C}^a_{10,11} = \real{\lambda^{2,0}}  \wedge
   \eta = \mathrm{alt} (\mathcal{C}_{10} + \mathcal{C}_{11}).
\end{gather*}
To be more precise, the component of $dF$ in
$\mathcal{C}^a_{10,11}$ is given by
$$
(dF)_{(10,11)} = - \eta \wedge \left( 2 (\nabla \eta)_{(10)} \circ
\varphi + \varphi (\nabla_{\zeta} F) \right).
$$
\begin{remark}
{\rm Note that $\varphi (\nabla_{\zeta} F)= \varphi (\nabla_{\zeta}
F)_{(11)} =  - (\nabla_{\zeta} F)_{(11)}$. So that, if an almost
contact metric structure is  of type $\mathcal{C}_{10} \oplus
\mathcal{C}_{11}$ such that $\nabla \eta \circ \varphi =
\nabla_{\zeta} F$, then $dF=0$.}
\end{remark}

Now using the properties satisfied by the different components of
$\nabla F$(see again \cite{ChineaGonzalezDavila}), we have the
following relations
\begin{gather*}
\| (d F)_{(1)} \|^2 = 9 \|  (\nabla F)_{(1)} \|^2; \quad \| (d
F)_{(i)} \|^2 =
 3 \|  (\nabla F)_{(i)} \|^2, \quad  i=3,4;\\
   \| (d F)_{(i)} \|^2 = 6 \|  (\nabla F)_{(i)} \|^2, \quad i=5,8; \quad \| (d F)_{(10,11)}
   \|^2= 3 \|2 (\nabla \eta)_{(10)} \circ \varphi - (\nabla_{\zeta}
   F)_{(11)}\|^2.
\end{gather*}
Let us also  recall that $\|(\nabla_{\zeta}
   F)_{(11)}\|^2= \|(\nabla
   F)_{(11)}\|^2$.

Next we consider the $\Lie{U}(n)$-map $\mathrm{c}_{1\,2} \, : \,
{\rm T}^* M \otimes \lie{u}(n)^\perp : \to T^* M$ defined by the
contraction $\mathrm{c}_{1\,2}(a) (X) = a(e_i,e_i,X)$. The
coderivative $d^* F= - \mathrm{c}_{1\,2} (\nabla F)$ is in
$\mathcal{C}^c_{4,12} \oplus \mathcal{C}^c_{6}$,  where
$\mathcal{C}^c_{4,12} = \eta^{\perp} = \mathrm{c}_{1\,2}
(\mathcal{C}_4 \oplus \mathcal{C}_{12})$ and  $\mathcal{C}^c_{6} =
\mathbb{R}\eta =  \mathrm{c}_{1\,2} (\mathcal{C}_{6})$. In fact,
 $ (d^* F)_{(4,12)} = -  e_i \lrcorner   (\nabla_{e_i} F)_{(4)}   + (\nabla_{\zeta} \eta)_{(12)} \circ \varphi
 $  and $ (d^* F)_{(6)} = d^* F (\zeta) \eta$.

\begin{remark}
{\rm Note that if an almost contact metric structure is  of type
$\mathcal{C}_{4} \oplus \mathcal{C}_{12}$ such that $  e_i \lrcorner
(\nabla_{e_i} F)_{(4)} =   (\nabla_{\zeta} \eta)_{(12)} \circ
\varphi$ or, equivalently, $\varphi e_i \lrcorner (\nabla_{\varphi
e_i} F) = (\nabla_{\zeta} \eta) \circ \varphi$, then $d^*F=0$.}
\end{remark}

Now using the properties of the different components of $\nabla
F$, we have
\begin{gather*}
\| (d^* F)_{(4,12)}  \|^2 =  \|   -  e_i \lrcorner   (\nabla_{e_i}
F)_{(4)} + (\nabla_{\zeta} \eta)_{(12)} \circ \varphi \|^2, \quad
\| (d^* F)_{(6)} \|^2 =
 n \|  (\nabla F)_{(6)} \|^2.
\end{gather*}

We will also recall that  $\|  e_i \lrcorner   (\nabla_{e_i}
F)_{(4)}\|^2 = \tfrac{n-1}{2} \|(\nabla F)_{(4)}\|^2$ and $\|
(\nabla \eta)_{(12)} \|^2 = \| (\nabla_{\zeta} \eta)_{(12)} \|^2$.
\vspace{2mm}

 Now we will analyse the exterior derivative of  the characteristic one-form $\eta$.   The decomposition
 of the space of skew-symmetric two-forms is given by
   $ \Lambda^2 {\rm T}^* M = \mathbb{R} \,  F \oplus
\real{\lambda^{1,1}_0} \oplus \real{\lambda^{2,0}} \oplus \eta
\wedge \eta^{\perp}$, where $\mathbb{R} \,  F \cong  \mathcal{C}_6$,
$\real{\lambda^{1,1}_0} \cong   \mathcal{C}_7$,
$\real{\lambda^{2,0}} \cong \mathcal{C}_{10}$ and $\eta \wedge
\eta^{\perp} \cong \mathcal{C}_{12}$. Therefore, making use of the
alternating map $\mathrm{alt} \, : \,\otimes^2 {\rm T}^* M \to
\Lambda^2 {\rm T}^* M$, for  the exterior derivative $d \eta =
\mathrm{alt} (\nabla \eta)$, we have
\begin{gather*}
 \| (d \eta)_{(i)} \|^2 = 4 \| (\nabla \eta)_{(i)} \|^2, \quad i=6,7,10; \quad   \| (d \eta)_{(12)}
\|^2 = 2 \| (\nabla \eta)_{(12)} \|^2.
\end{gather*}
Finally, it is easy to see that $d^* \eta = - (\nabla_{e_i}
\eta)_{(5)} e_i$. Therefore, $(d^* \eta)^2 = 2n \| (\nabla \eta
)_{(5)} \|^2$.

Taking all of this into account, we obtain the following Bochner
type formulas
\begin{gather}
\int_M  \left( 8 \| \xi_{(1)} \|^2  -4 \| \xi_{(2)} \|^2
 +  \| \xi_{(5)} \|^2 +  (2n-1)\| \xi_{(6)} \|^2
  -   \| \xi_{(7)} \|^2 +  \| \xi_{(8)} \|^2 -   \| \xi_{(9)} \|^2
  -  \| \xi_{(10)} \|^2
   \right.  \nonumber\\
  - 4   \| \xi_{(11)} \|^2  -   \| \xi_{(12)} \|^2
 +  \|2 (\nabla \eta)_{(10)} \circ \varphi - (\nabla_{\zeta}
   F)_{(11)}\|^2
\left. + 2 \|   -  e_i \lrcorner   (\nabla_{e_i} F)_{(4)} +
(\nabla_{\zeta} \eta)_{(12)} \circ \varphi \|^2 \right) \nonumber  \\
\, \qquad \qquad \qquad \qquad \qquad \qquad \qquad \qquad = 2
\int_M \left(s- s^{\rm ac}- \Ric (\zeta , \zeta)\right),
\label{bochnereq1}
\\
\label{bochnereq2} \int_M \left(
   (2n-1) \| \xi_{(5)} \|^2 +   \| \xi_{(6)} \|^2 +  \| \xi_{(7)} \|^2
- \| \xi_{(8)} \|^2 -   \| \xi_{(9)} \|^2 +  \| \xi_{(10)}
\|^2\right)
       =  2 \int_M   \Ric (\zeta , \zeta).
\end{gather}

Now restricting our attention on conformally flat manifolds, we
display the next result where we will denote by $\varphi_{\sigma}$,
$F_{\sigma}$, $\eta_{\sigma}$, $\zeta_{\sigma}$ and
$\xi_{(i)}(\sigma)$  the corresponding tensors associated to an
almost contact metric structure $\sigma$.
\begin{theorem} \label{minimiser}
 Let $(M, \langle \cdot ,\cdot \rangle)$ be a
$2n+1$-dimensional conformally flat  compact Riemannian manifold,
where $n>1$. If  $s$ is the scalar curvature and $C_{\langle \cdot
,\cdot \rangle} = \int_M s$, that is,  a constant depending on the
metric ${\langle \cdot ,\cdot \rangle}$, then every almost contact
metric structure $\sigma$ compatible with $\langle \cdot ,\cdot
\rangle$ satisfies
\begin{eqnarray} \nonumber
\tfrac{2(n-1)}{2n-1}  C_{\langle \cdot ,\cdot \rangle} & =& \int_M
\left( 4 \| \xi_{(1)}(\sigma) \|^2 - 2 \| \xi_{(2)}(\sigma) \|^2
 +  (n-1) \| \xi_{(5)}(\sigma) \|^2 +  \tfrac{(2n+1)(n-1)}{2n-1}\| \xi_{(6)}(\sigma) \|^2 \right. \\
& & - \tfrac{1}{2n-1}   \| \xi_{(7)}(\sigma) \|^2   +
\tfrac{1}{2n-1} \| \xi_{(8)}(\sigma) \|^2 - \tfrac{2(n-1)}{2n-1}\|
\xi_{(9)}(\sigma) \|^2
  - \tfrac{1}{2n-1}    \| \xi_{(10)}(\sigma) \|^2 \nonumber  \\
  && \left.
  -  2 \| \xi_{(11)}(\sigma) \|^2
  -   \| \xi_{(12)}(\sigma) \|^2  +   \tfrac12 \|2 (\nabla \eta_\sigma)_{(10)} \circ \varphi_\sigma -
(\nabla_{\zeta_\sigma}
   F_\sigma)_{(11)}\|^2 \right.  \nonumber \\
   &&
\left. +  \|   -  e_i \lrcorner   (\nabla_{e_i} F_\sigma)_4 +
(\nabla_{\zeta_\sigma} \eta_\sigma)_{(12)} \circ \varphi_\sigma
\|^2\right). \label{BorLam}
\end{eqnarray}
Moreover:
\begin{enumerate}
 \item[${\rm (i)}$] If $\sigma_0$ is an almost contact structure compatible with $\langle
\cdot ,\cdot \rangle$ of type $\mathcal{C}_1 \oplus \mathcal{C}_4$
and $n =3$, then $\sigma_0$ is an energy minimiser such that its
total bending is $B(\sigma_0) = \tfrac{1}{10} C_{\langle \cdot
,\cdot \rangle}$. Furthermore, in this situation any other energy
minimiser is of type $\mathcal{C}_1 \oplus \mathcal{C}_4$.

\item[${\rm (ii)}$] If $n=2$ or $n\geq 4$, and  $\sigma_0$ is an
almost contact structure compatible with $\langle \cdot ,\cdot
\rangle$ of type $\mathcal{C}_4$,  then $\sigma_0$ is an energy
minimiser such that its total bending is $B(\sigma_0) =
\tfrac{1}{2(2n-1)} C_{\langle \cdot ,\cdot \rangle}$. Furthermore,
in this situation any other energy minimiser is of type
$\mathcal{C}_4$.

\item[${\rm (iii)}$] If $\sigma_0$ is an almost contact structure
compatible with $\langle \cdot ,\cdot \rangle$ of type
$\mathcal{C}_2  $, then $\sigma_0$ is an energy minimiser such
that its total bending  is $B(\sigma_0) = - \tfrac{n-1}{2(2n-1)}
C_{\langle \cdot ,\cdot \rangle} $. Furthermore, in this situation
any other energy minimiser is of type $\mathcal{C}_2  $.
\end{enumerate}
 \end{theorem}
\begin{proof}
If $n>1$ and $M$ is  conformally flat, the Ricci curvature tensor
completely determines the Riemannian curvature $R$. Thus we have
$$ R = \tfrac{1}{2n-1} \langle
\cdot ,\cdot \rangle \ovee \Ric  - \tfrac{s}{4n(2n-1)} \langle
\cdot ,\cdot \rangle \ovee \langle \cdot ,\cdot \rangle.
$$
We recall that $\ovee$ is the  Kulkarni-Nomizu product defined by
equation \eqref{KulNom}. Therefore, for the almost contact Ricci
tensor we have
$$
\Ric^\mathrm{ac} =  \tfrac{1}{2n-1}\left( (1 + \varphi)\Ric -
(\zeta \lrcorner \Ric) \otimes \eta - \tfrac{s}{2n} (\langle \cdot
,\cdot \rangle - \eta \otimes \eta)\right).
$$
Hence the almost contact scalar curvature is given by
$s^\mathrm{ac} =  \tfrac{1}{2n-1} ( s - 2\Ric (\zeta,\zeta))$. As
a consequence,
$$
 s- s^{\rm ac} -  \Ric (\zeta , \zeta) = \tfrac{2(n-1)}{2n-1} s -
 \tfrac{2n-3}{2n-1}\Ric (\zeta,\zeta).
$$
Now using the Bochner type formulas \eqref{bochnereq1} and
\eqref{bochnereq2}, we get the equation \eqref{BorLam}.

For ${\rm (i)}$ and  ${\rm (ii)}$. If $n > 3$ and $\sigma_0$ is of
type $\mathcal C_4$,  then for any almost contact metric structure
$\sigma$, we have
\begin{eqnarray*}
2 (n - 1) \int_M \| \xi(\sigma) \|^2 &\geq & \int_M \left( 4 \|
\xi_{(1)}(\sigma)  \|^2 - 2 \| \xi_{(2)}(\sigma)  \|^2
 + 2 (n-1) \| \xi_{(4)}(\sigma)  \|^2  \right. \\
 &&
  + (n-1) \| \xi_{(5)}(\sigma)  \|^2
  + (n - \tfrac{1}{2n-1})\| \xi_{(6)}(\sigma)  \|^2
   -  \tfrac{1}{2n-1}   \| \xi_{(7)}(\sigma) \|^2
 \\
 &&
  \left.  + \tfrac{1}{2n-1}  \| \xi_{(8)} (\sigma) \|^2
  -  (1-  \tfrac{1}{2n-1})  \| \xi_{(9)} (\sigma)\|^2
  + (1 - \tfrac{1}{2n-1})    \| \xi_{(10)} (\sigma) \|^2
   \right)\\
   &  \geq & \tfrac{2(n-1)}{2n-1}  C_{\langle \cdot ,\cdot \rangle} = 2 (n-1) \int_M \| \xi_{(4)}(\sigma_0)
   \|^2 \\
   & = & 2(n-1)  \int_M \| \xi(\sigma_0)\|^2,
\end{eqnarray*}
where we have used the inequalities
\begin{gather*}
  \| \xi_{(10)} (\sigma)  \|^2 +   2 \| \xi_{(11)} (\sigma) \|^2 =
\tfrac12 \| 2(\nabla \eta_\sigma)_{(10)} \circ \varphi_\sigma \|^2
+
\tfrac12 \|(\nabla_{\zeta_\sigma} F_\sigma)_{(11)}\|^2 \\
\qquad \qquad\qquad\qquad\quad \geq \tfrac12 \|2 (\nabla
\eta_\sigma)_{(10)} \circ \varphi_\sigma - (\nabla_{\zeta_\sigma}
   F_\sigma)_{(11)}\|^2, \\
2(n-1)    \| \xi_{(4)} (\sigma)  \|^2 +  \tfrac12    \| \xi_{(12)}
(\sigma) \|^2 =     \| e_i \lrcorner   (\nabla_{e_i}
F_\sigma)_{(4)} \|^2 +     \| (\nabla_{\zeta_\sigma}
\eta_\sigma)_{(12)} \circ \varphi_\sigma \|^2 \\
\qquad \qquad\qquad\qquad\qquad \qquad \quad     \geq  \| - e_i
\lrcorner (\nabla_{e_i} F_\sigma)_{(4)} + (\nabla_{\zeta_\sigma}
\eta_\sigma)_{(12)} \circ \varphi_\sigma \|^2.
\end{gather*}
Also note that $\xi_{(1)} (\sigma)$ only  appears when $n \geq 3$.
Therefore, the case $n=2$ follows by a similar discussion.

Finally, if $\sigma_1$ is a almost contact metric structure which is
an energy minimiser, then we have
\begin{eqnarray*}
 \tfrac{2(n-1)}{2n-1}  C_{\langle \cdot ,\cdot \rangle} & \leq &
\int_M \left( 4 \| \xi_{(1)}(\sigma_1) \|^2 - 2 \|
\xi_{(2)}(\sigma_1) \|^2 + 2(n-1) \| \xi_{(4)}(\sigma_1) \|^2
 + (n-1) \| \xi_{(5)}(\sigma_1) \|^2  \right. \\
&  & + (n - \tfrac{1}{2n-1})\| \xi_{(6)}(\sigma_1) \|^2  -
\tfrac{1}{2n-1}   \| \xi_{(7)}(\sigma_1) \|^2   + \tfrac{1}{2n-1}
\|
\xi_{(8)}(\sigma_1) \|^2 \\
& &  -  (1- \tfrac{1}{2n-1})\| \xi_{(9)}(\sigma_1) \|^2
  + (1 - \tfrac{1}{2n-1} )   \| \xi_{(10)}(\sigma_1) \|^2 \\
  & \leq &   2(n-1)  \int_M  \| \xi(\sigma_1) \|^2 =  \tfrac{2(n-1)}{2n-1}  C_{\langle \cdot ,\cdot
  \rangle}.
\end{eqnarray*}
From this, the inequalities are really equalities. As a
consequence, we get
 \begin{eqnarray*}
0 & = & \int_M \left( 2(n- 3) \| \xi_{(1)}(\sigma_1)\|^2 + 2 n \|
\xi_{(2)}(\sigma_1)\|^2 + 2(n - 1)\| \xi_{(3)}(\sigma_1)\|^2
\right.
 \\
&& +  (n-1) \| \xi_{(5)}(\sigma_1) \|^2  + (n -2 +
\tfrac{1}{2n-1})\| \xi_{(6)}(\sigma_1) \|^2 + (2 (n-1) +
\tfrac{1}{2n-1} )  \| \xi_{(7)}(\sigma_1) \|^2 \\
&&  +  (2(n-1) - \tfrac{1}{2n-1})\| \xi_{(8)}(\sigma_1) \|^2
 +   (2n-1 - \tfrac{1}{2n-1})\| \xi_{(9)}(\sigma_1)\|^2
\\
&&
 \left.  + (2n - 3 + \tfrac{1}{2n-1} )\|
\xi_{(10)}(\sigma_1)\|^2  +  2(n - 1) \| \xi_{(11)}(\sigma_1)\|^2
  +  2(n -1) \| \xi_{(12)}(\sigma_1)\|^2
\right).
 \end{eqnarray*}
Since  for $n>1$ and $n \neq 3$,  all the coefficients  are
positive, it is obtained
\begin{gather*}
\xi_{(1)}(\sigma_1) =\xi_{(2)}(\sigma_1) =\xi_{(3)}(\sigma_1)
=\xi_{(5)}(\sigma_1) =\xi_{(6)}(\sigma_1) = \xi_{(7)}(\sigma_1)
\\=\xi_{(8)}(\sigma_1) =\xi_{(9)}(\sigma_1) = \xi_{(10)}(\sigma_1) =
\xi_{(11)}(\sigma_1) = \xi_{(12)}(\sigma_1) =0.
\end{gather*}
Therefore, the structure $\sigma_1$ is of type $\mathcal{C}_4$.
The proof for $n=3$ of the analogous assertion can be similarly
done. \vspace{2mm}

For ${\rm (iii)}$. If $\sigma_0$ is an almost contact metric
structure of type $\mathcal{C}_2$, then for any almost contact
metric structure $\sigma$ we have {\small
\begin{eqnarray*}
  2  \int_M \| \xi(\sigma) \|^2 &\geq & \int_M \left(  - 4 \|
\xi_{(1)}(\sigma) \|^2 + 2 \| \xi_{(2)}(\sigma) \|^2
 - (n-1) \| \xi_{(5)}(\sigma) \|^2 - (n - \tfrac{1}{2n-1})\| \xi_{(6)}(\sigma) \|^2   \right.
 \\
 &&
 \quad + \tfrac{1}{2n-1}   \| \xi_{(7)} (\sigma) \|^2  - \tfrac{1}{2n-1}  \| \xi_{(8)} (\sigma) \|^2
  +  (1- \tfrac{1}{2n-1})  \| \xi_{(9)} (\sigma) \|^2
  \\
  &&
  \quad + \tfrac{1}{2n-1}    \| \xi_{(10)} (\sigma) \|^2 +  2  \| \xi_{(11)} (\sigma) \|^2
  +  \| \xi_{(12)} (\sigma)  \|^2
\\
&&
  \quad \left. -  \tfrac12 \|2 (\nabla \eta_\sigma)_{(10)} \circ
\varphi_\sigma - (\nabla_{\zeta_\sigma}
   F)_{(11)}\|^2
 -  \|   -  e_i \lrcorner   (\nabla_{e_i} F_\sigma)_{(4)} +
(\nabla_{\zeta_\sigma}
\eta_\sigma)_{(12)} \circ \varphi_\sigma \|^2 \right)  \\
  &=&      -  \tfrac{2(n-1)}{2n-1}  C_{\langle \cdot ,\cdot \rangle} =   \int_M \| \xi_{(2)}(\sigma_0) \|^2 =
    \int_M \| \xi(\sigma_0) \|^2 .
\end{eqnarray*}}

If $\sigma_1$ is a almost contact metric structure which is an
energy minimiser, then we have {\small
\begin{eqnarray*}
    2 \int_M \| \xi(\sigma_1) \|^2  & =&  - \tfrac{2(n-1)}{2n-1}
C_{\langle
\cdot ,\cdot \rangle}  \\
& \geq &   \int_M \left(  - 4 \| \xi_{(1)}(\sigma_1) \|^2 + 2 \|
\xi_{(2)}(\sigma_1) \|^2
 - (n-1) \| \xi_{(5)}(\sigma_1) \|^2 - (n - \tfrac{1}{2n-1})\| \xi_{(6)}(\sigma_1) \|^2 \right.
 \\
 &&
 \quad  + \tfrac{1}{2n-1}   \| \xi_{(7)} (\sigma_1) \|^2  - \tfrac{1}{2n-1}  \| \xi_{(8)} (\sigma_1) \|^2
    +    (1- \tfrac{1}{2n-1})   \| \xi_{(9)} (\sigma_1) \|^2
      \\
   &&
   \quad + \tfrac{1}{2n-1}    \| \xi_{(10)} (\sigma_1) \|^2
   +   2 \| \xi_{(11)} (\sigma_1) \|^2
  +   \| \xi_{(12)} (\sigma_1)  \|^2
  \\
  &&
  \quad \left.
-  \tfrac12 \|2 (\nabla \eta)_{(10)} \circ \varphi - (\nabla_{\zeta}
   F)_{(11)}\|^2
    -   \|   -  e_i \lrcorner   (\nabla_{e_i} F)_{(4)} + (\nabla_{\zeta}
\eta)_{(12)} \circ \varphi \|^2 \right)  \\
  &      = & -  \tfrac{2(n-1)}{2n-1}  C_{\langle \cdot ,\cdot
  \rangle}.
\end{eqnarray*}}
Hence the inequality is really an equality. Therefore,
\begin{eqnarray*}
0 &=& \int_M \left(
 6 \| \xi_{(1)}(\sigma_1)\|^2
 + 2 \|\xi_{(3)}(\sigma_1)\|^2
 +  2 \| \xi_{(4)}(\sigma_1)\|^2
 + ( n+ 1) \| \xi_{(5)}(\sigma_1)\|^2 \right.
\\
&&
 +  (n+  2 -\tfrac{1}{2n-1}) \| \xi_{(6)}(\sigma_1)\|^2
 +   ( 2 - \tfrac{1}{2n-1}) \| \xi_{(7)}(\sigma_1)\|^2
 +  ( 2 + \tfrac{1}{2n-1})  \| \xi_{(8)} (\sigma_1) \|^2
 \\
 &&
 +  ( 1 + \tfrac{1}{2n-1})  \| \xi_{(9)} (\sigma_1) \|^2
 +  (2 - \tfrac{1}{2n-1}) \| \xi_{(10)}(\sigma_1)\|^2
  +     \| \xi_{(12)}(\sigma_1)\|^2
  \\
 &&\left. +  \tfrac12 \|2 (\nabla \eta)_{(10)} \circ \varphi -
(\nabla_{\zeta}
   F)_{(11)}\|^2
 +  \|   -  e_i \lrcorner   (\nabla_{e_i} F)_{(4)} + (\nabla_{\zeta} \eta)_{(12)} \circ \varphi \|^2 \right).
\end{eqnarray*}
Since for $n>1$ all the coefficients  are again positive, it is
obtained
\begin{gather*}
\xi_{(1)}(\sigma_1) =\xi_{(3)}(\sigma_1) =\xi_{(4)}(\sigma_1)
=\xi_{(5)}(\sigma_1) =\xi_{(6)}(\sigma_1) = \xi_{(7)}(\sigma_1)
\\=\xi_{(8)}(\sigma_1) =\xi_{(9)}(\sigma_1) = \xi_{(10)}(\sigma_1) =
\xi_{(11)}(\sigma_1) = \xi_{(12)}(\sigma_1) =0.
\end{gather*}
Thus, we conclude that the structure $\sigma_1$ is of type
$\mathcal{C}_2$.
\end{proof}

Finally, we consider the situation for compact Einstein manifolds of
dimension $2n+1$.
\begin{theorem} \label{minimisereinstein}
 Let $(M, \langle \cdot ,\cdot \rangle)$ be a
$2n+1$-dimensional   compact Einstein manifold. If $s$ is the scalar
curvature,  then every almost contact structure $\sigma$ compatible
with $\langle \cdot ,\cdot \rangle$ satisfies
\begin{eqnarray} \label{bochnereq2borlam}
\int_M
 \left(   (2n-1) \| \xi_{(5)}(\sigma) \|^2 +   \| \xi_{(6)}(\sigma) \|^2 +  \| \xi_{(7)}(\sigma)
 \|^2 \right.
&&   \\
\left.  - \| \xi_{(8)}(\sigma) \|^2 -   \| \xi_{(9)}(\sigma) \|^2
+ \| \xi_{(10)} (\sigma) \|^2 \right)
      & = & \tfrac{2s}{2n+1} {\rm Vol}(M). \nonumber
\end{eqnarray}

Moreover:
\begin{enumerate}
 \item[${\rm (i)}$] For $n=1$, if  $\sigma_0$ is an almost contact structure compatible with $\langle
\cdot ,\cdot \rangle$ of type $\mathcal{C}_5 \oplus \mathcal{C}_6$
$($trans-Sasakian$)$, then $\sigma_0$ is an energy minimiser such
that its total bending is $B(\sigma_0) = \tfrac{s}{3} {\rm
Vol}(M)$. Furthermore, in this situation any other energy
minimiser is trans-Sasakian.

\item[${\rm (ii)}$] For $n>1$, if $\sigma_0$ is an almost contact
structure compatible with $\langle \cdot ,\cdot \rangle$ of type
$\mathcal{C}_5$ $(\alpha$-Kenmotsu, where $2n \alpha = - d^*
\eta_{\sigma_0})$, then $\sigma_0$ is an energy minimiser such
that its total  bending is $B(\sigma_0) = \tfrac{s}{4n^2 -1} {\rm
Vol}(M) $. Furthermore, in this situation any other energy
minimiser is of type $\alpha$-Kenmotsu with $2n \int_M \alpha^2 =
\tfrac{s}{4n^2 -1} {\rm Vol}(M)$.
 \item[${\rm (iii)}$] If $\sigma_0$ is an almost
contact structure compatible with $\langle \cdot ,\cdot \rangle$
of type $\mathcal{C}_8 \oplus \mathcal{C}_9$, such that its total
bending is $B(\sigma_0) = - \tfrac{s}{2n+1} {\rm Vol}(M) $.
Furthermore, in this situation any other energy minimiser is of
type $\mathcal{C}_8 \oplus \mathcal{C}_9$.
\end{enumerate}
 \end{theorem}
 \begin{proof} Equation \eqref{bochnereq2borlam} follows from
the  equation \eqref{bochnereq2}. The remaining assertions can be
deduced by using similar arguments as in the proof of Theorem
 \ref{minimiser}.
 \end{proof}

\section{Examples}\setcounter{equation}{0}

\noindent {\bf Odd dimensional spheres.} Let $\mathbb{C}^{n+1}$ be
the $n+1$-dimensional complex space equipped the natural Kähler
structure denoted by $(\langle \cdot , \cdot \rangle,J)$. If we
consider the sphere $S^{2n+1}(r)$ of radius $r$ and  centered at the
origin and $U$ is a unit normal vector field on $S^{2n+1}(r)$, then
there is an $a$-Sasakian structure on $S^{2n+1}(r)$ fixing the
restrictions of $\langle \cdot , \cdot \rangle$ and the tangent
projection ${\rm tan} \circ J$ to $S^{2n+1}(r)$ as metric and
$(1,1)$-tensor field, respectively, and $\zeta=J U$ as
characteristic vector field. This example for $r=1$ have been shown
by Yano and Kon in \cite{YanoKon}. In this situation
$$
a^2 = \tfrac1{r^2}, \quad \Ric = \tfrac{2n}{r^2} \langle \cdot ,
\cdot \rangle, \quad \Ric^{\rm ac} = \tfrac{1}{r^2} \left( \langle
\cdot , \cdot \rangle - \eta \otimes \eta\right).
$$
Hence $S^{2n+1}(r)$ is ac-Einstein. Now,  using Theorem
\ref{mapharmcontact} (i) (b), the $\pm \tfrac1{r}\,$-Sasakian
structure is
 a harmonic map.

 For $n=1$, the $\pm \tfrac1{r}\,$-Sasakian structure is an energy
 minimiser by Theorem \ref{minimisereinstein}. The energy of this
 structure is
 $\tfrac23 \pi r (3r^2 + 4)$.
 \vspace{1mm}

 \noindent {\bf The product manifold $S^6 \times S^1$.} It is well known
 that $S^6$ can be endowed by an nearly Kähler structure by means
 of the round metric and the standard integrable
 $\Lie{G}_2$-structure on $\mathbb R^7$ \cite{Gray:vector,Gray:six}. Now,  we consider
 the  product manifold $S^6 \times S^1$. Taking the product
 metric, the unit tangent vector $\zeta$ defined by the Maurer-Cartan form
 on $S^1$ and the $(1,1)$ tensor field $\varphi$ defined by the
 almost complex structure on $S^6$ and  $\varphi \zeta
 =0$, we will obtain an almost contact metric structure $\sigma$ on $S^6 \times
 S^1$ of type $\mathcal C_1$. Since $S^6 \times S^1$ is
 conformally flat, the structure $\sigma$ is an energy minimiser
 by Theorem \ref{minimiser}.
\vspace{1mm}

\noindent {\bf The product manifold $S^{2m+1} \times S^1 \times
S^1$.} Hopf manifolds are diffeomorphic to $ S^{2m+1}\times S^1 $
and admit a locally conformal K\"{a}hler structure (type $\mathcal
W_4$ in Gray-Hervella's terms) with parallel Lee form $\theta$
\cite{Va} (see \cite{Gray-H:16} for the classification of almost
Hermitian structures).  In general, if we consider a $1$-Sasakian
manifold $M$ of dimension $2n+1$, then the product manifold $M
\times S^1$ can be equipped with a locally conformal K\"{a}hler
structure. Now, if one proceed as in the previous example, an almost
contact metric structure $\sigma$ of type $\mathcal{C}_4$ can be
defined on $S^{2m+1} \times S^1 \times S^1$. Since $S^{2m+1} \times
S^1 \times S^1$ is conformally flat, the structure $\sigma$ is an
energy minimiser
 by Theorem \ref{minimiser}.
\vspace{1mm}

\noindent {\bf The hyperbolic space.} Let $({\mathbb H}^{2n+1},
\langle \cdot,\cdot \rangle)$ be the $(2n+1)$-dimensional
Poincar\'{e} half-space, i.e.
\[
{\mathbb H}^{2n+1} = \{(x_{1},\dots ,x_{2n+1})\in
\mathbb{R}^{2n+1}\mid x_{1}> 0\}
\]
and $\langle\cdot,\cdot\rangle$ is the Riemannian metric given by
$\langle\cdot,\cdot\rangle =
(cx_{1})^{-2}\sum_{i=1}^{2n+1}(dx_{i})^{2},$ $c>0.$ It has
constant curvature $-c^{2}.$ Moreover, ${\mathbb H}^{2n+1}$ can be
considered as the Lie group with the following product
\[
(x_{1},\dots ,x_{2n+1})\cdot (y_{1},\dots ,y_{2n+1}) =
(x_{1}y_{1}, x_{1}y_{2} + x_{2}, \dots ,x_{1}y_{2n+1} + x_{2n+1}).
\]
So we have a solvable Lie group which is a semi-direct product of
the multiplicative group $\mathbb{R}^{+}_{0}$ and the additive group
$\mathbb{R}^{2n}$ and $\langle\cdot,\cdot\rangle$ becomes into a
left-invariant metric. In fact, the vector fields
\[
X_{i} = cx_{1}\tfrac{\partial}{\partial x_{i}},\;\;\; i = 1,\dots
,2n+1,
\]
form an orthonormal basis of left-invariant vector fields and the
bracket satisfies $[X_{1},X_{j}] = cX_{j},$ $j = 2,\dots ,2n+1,$
the other brackets being zero. Then the Levi Civita connection
with respect to $\langle\cdot,\cdot\rangle$ is determined by
\begin{equation}\label{hypcon}
\nabla_{X_{j}}X_{j} = cX_{1},\;\;\;\; \nabla_{X_{j}}X_{1} =
-cX_{j},\;\;\;\; j = 2,\dots ,2n+1,
\end{equation}
being zero the remaining covariant derivatives with respect to this
basis.
\begin{proposition}
Any left-invariant almost contact structure $(\varphi,\zeta,\eta)$
on ${\mathbb H}^{2n+1}$ compatible with
$\langle\cdot,\cdot\rangle$ with $\zeta = X_{1}$ is a harmonic
structure of type ${\mathcal C}_{5}$ but it does not determine a
harmonic map.
\end{proposition}
\begin{proof} From (\ref{hypcon}) we obtain that this structure is
of type ${\mathcal C}_{5}$ and so, it is harmonic by using Theorem
\ref{harmclasscontact}. Moreover, because ${\mathbb H}^{2n+1}$ has
constant sectional curvature, it follows that it is ac-Einstein and
$s^{\rm ac}=-2nc^2.$ Hence, taking into account that $d^{*}\eta=2nc$
does not vanish, Theorem \ref{mapharmcontact} (i) (a) implies that
it does not determine a harmonic map.
\end{proof}
\vspace{1mm}

\noindent {\bf The generalised Heisenberg groups H(1,r).} Now we
consider some special classes of Lie groups which play an important
r\^{o}le in geometry. Namely,  the generalised Heisenberg groups
$H(1,r)$ introduced in \cite{Haraguchi} and the ones in the sense of
Kaplan $H(p,1)$ \cite{Kaplan}. The Lie group $H(1,r)$ consists of
matrices
\[
a = \left (
\begin{array}{ccc}
I_{r} & A^{t} & B^{t}\\
0 & 1 & c\\
0 & 0 & 1
\end{array}
\right ),
\]
where $I_{r}$ denotes the identity matrix of type $r\times r,$ $A =
(a_{1},\dots ,a_{r})\in {\mathbb{R}}^{r},$ $B = (b_{1},\dots
,b_{r})\in {\mathbb{R}}^{r}$ and $c\in \mathbb{R}.$ Moreover,
$H(1,r)$ is connected, simply connected,  nilpotent, of dimension
$2r+1$ and its center is of dimension $r.$ The following coordinates
$(x^{i},x^{r+i},z),$ $1\leq i\leq r,$ provide a system of global
coordinates on $H(1,r):$
\[
x^{i}(a) = a_{i},\;\;\; x^{r+i}(a) = b_{i},\;\;\; z(a) = c
\]
and a basis of left-invariant one-forms is given by $\alpha_{i} =
dx^{i},$ $\alpha_{r+i} = dx^{r+i} - x^{i}dz,$ $\gamma = dz.$ For
the dual left-invariant vector fields we then have
\[
X_{i} = \tfrac{\partial}{\partial x^{i}},\;\;\; X_{r+i} =
\tfrac{\partial}{\partial x^{r+i}},\;\;\; Z =
\tfrac{\partial}{\partial z} + \textstyle\sum_{j = 1} ^{r}
\displaystyle x^{j}\tfrac{\partial}{\partial x^{r+j}}.
\]
On $H(1,r)$ we consider the Riemannian metric $\langle \cdot,
\cdot \rangle$ for which these vector fields form an orthonormal
basis. Then, the corresponding Levi Civita connection is
determined by
\begin{equation}\label{nablaH}
\begin{array}{lclcl}
\nabla_{X_{i}}X_{r+i} & = & \nabla_{X_{r+i}}X_{i} & = &
-\frac{1}{2}Z,\\[1mm]
\nabla_{X_{i}}Z & = &-\nabla_{Z}X_{i} & = &
\frac{1}{2}X_{r+i},\\[1mm]
\nabla_{X_{r+i}}Z & = & \nabla_{Z}X_{r+i} & = & \frac{1}{2}X_{i},
\end{array}
\end{equation}
where the remaining covariant derivatives with respect to this
basis vanish. For the nonvanishing components of the curvature
tensor $R$ we obtain
\begin{equation}\label{RH}
\begin{array}{lclclcl}
R(X_{i}X_{j},X_{r+i},X_{r+j}) & = & -\frac{1}{4},\;\; i\neq j, & &
R(X_{i},X_{r+j},X_{j},X_{r+i}) & = & \frac{1}{4},\\[0.5pc]
R(X_{i},Z,X_{i},Z) & = & -\frac{3}{4}, & & R(X_{r+i},Z,X_{r+i},Z)
& = & \frac{1}{4}.
\end{array}
\end{equation}

Next, we consider left-invariant almost contact structures
$(\varphi,\zeta,\eta)$ on $H(1,r)$ compatible with $\langle \cdot,
\cdot \rangle$ with $\zeta = Z$ and $\eta = \gamma.$ Denote by
$\varphi^{l}_{k},$ $k,l = 1,\dots ,2r,$ the (constant) components
of $\varphi$ with respect to the basis $\{X_{k},Z\}.$ Then, we
have
\begin{lemma} Any left-invariant almost contact metric structure
$(\varphi,\zeta= Z,\eta= \gamma, \langle \cdot, \cdot \rangle)$ on
$H(1,r)$ is of type $({\mathcal C}_{8}\oplus {\mathcal
C}_{9}\oplus {\mathcal C}_{11})-{\mathcal C}_{11}.$ Moreover, it
is of type
\begin{enumerate}
\item[{\rm (i)}] ${\mathcal C}_{8}\oplus {\mathcal C}_{9}$ if and
only if $\varphi^{r+j}_{i} = \varphi_{j}^{r+i}$ and
$\varphi^{r+j}_{r+i} = \varphi^{j}_{i};$ \item[{\rm (ii)}]
${\mathcal C}_{8}$ if and only if $\varphi^{r+j}_{i} =
\varphi_{j}^{r+i} = 0$ and $\varphi^{r+j}_{r+i} =
\varphi^{j}_{i}.$ Then, $r$ must be even; \item[{\rm (iii)}]
${\mathcal C}_{9}$ if and only if $\varphi^{r+j}_{i} =
\varphi_{j}^{r+i}$ and $\varphi^{r+j}_{r+i} = \varphi^{j}_{i} =
0,$
\end{enumerate}
for all $i,j\in \{1,\dots ,r\}.$
\end{lemma}
\begin{proof} From (\ref{nablaH}), we get
\[
(\nabla_{X_{i}}\varphi)X_{k}  =
-\tfrac{1}{2}\varphi^{r+i}_{k}Z,\;\;\;\;\;\;
(\nabla_{X_{r+i}}\varphi)X_{k} = -\tfrac{1}{2}\varphi^{i}_{k}Z,
\]
where $1\leq i\leq r$ and $1\leq k\leq 2r.$ Also we obtain
$$
\begin{array}{lcl}
(\nabla_{Z}\varphi)X_{i} & = & \frac{1}{2}
\sum_{j=1}^{r}\{(\varphi_{i}^{r+j} + \varphi_{r+i}^{j})X_{j} +
(\varphi_{r+i}^{r+j} - \varphi_{i}^{j})X_{r+j}\},\\[0.5pc]
(\nabla_{Z}\varphi)X_{r+i} & = & \frac{1}{2}
\sum_{j=1}^{r}\{(\varphi_{r+i}^{r+j} - \varphi_{i}^{j})X_{j} -
(\varphi_{r+i}^{j} + \varphi_{i}^{r+j})X_{r+j}\}.
\end{array}
$$
Then $d^{*}F(Z) = d^{*}\eta =0$ and $(\nabla_{Z}\varphi) = 0$ if and
only if $\varphi^{r+j}_{i} = \varphi_{j}^{r+i}$ and
$\varphi_{r+i}^{r+j} = \varphi_{i}^{j}.$ Moreover, these last
conditions imply
$$
\begin{array}{l}
 (\nabla_{X_{i}}\varphi)X_{j} =
(\nabla_{X_{j}}\varphi)X_{i} = -(\nabla_{X_{r+i}}\varphi)X_{r+j} =
-(\nabla_{X_{r+j}}\varphi)X_{r+i}\\
\hspace{1.85cm}  = -(\nabla_{\varphi X_{i}}\varphi)\varphi X_{j} =
(\nabla_{\varphi X_{r+i}}\varphi)\varphi X_{r+j} =
-\frac{1}{2}\varphi^{r+j}_{i}Z;\\[1mm]
(\nabla_{X_{i}}\varphi)X_{r+j} = - (\nabla_{X_{r+j}}\varphi)X_{i}
=(\nabla_{\varphi X_{i}}\varphi) \varphi X_{r+j} = -
(\nabla_{\varphi X_{r+j}}\varphi)\varphi X_{i} =
\frac{1}{2}\varphi_{i}^{j}Z.
\end{array}
$$

 This proves the Lemma.
\end{proof}

\begin{remark}{\rm The almost contact metric structure
$(\varphi,\zeta = Z,\eta= \gamma,\langle\cdot,\cdot \rangle)$ with
$\varphi X_{i} = X_{r+i}$ and $\varphi X_{r+i} = -X_{i},$
$i=1,\dots ,r,$ is of type ${\mathcal C}_{9}$ and so, it is almost
cosymplectic but non-cosymplectic. Moreover, for $r$ even and
taking $\varphi X_{2i -1} = X_{2i}$ and $\varphi X_{2i} = -
X_{2i-1},$ it follows that $(\varphi,\zeta = Z,\eta=
\gamma,\langle\cdot,\cdot \rangle)$ is of type} ${\mathcal
C}_{8}.$
\end{remark}

\begin{proposition}\label{harmonicmap} Any harmonic left-invariant almost contact metric structure
$(\varphi,\zeta =Z,\eta= \gamma,\langle\cdot ,\cdot \rangle)$ on
$H(1,r)$ determines a harmonic map.
\end{proposition}
\begin{proof} From (\ref{nablaH}), the intrinsic torsion $\xi$ of $(\varphi,\zeta
= Z,\eta= \gamma,\langle\cdot,\cdot \rangle)$ satisfies
\begin{equation}\label{intrinsicH}
\begin{array}{lclclcl}
\xi_{X_{i}}Z & = & -\frac{1}{2}X_{r+i}, & & \xi_{X_{r+i}}Z & = &
-\frac{1}{2}X_{i},\\
\xi_{X_{i}}X_{j} & = & \xi_{X_{r+i}}X_{r+j} =0, & &
\xi_{X_{i}}X_{r+j} & = & \xi_{X_{r+i}}X_{j} =0,\;i\neq j.
\end{array}
\end{equation}
Then, for the one-form $\nu$ defined in \eqref{harmmap1},  one
obtains by using (\ref{RH})
$$
\begin{array}{lcl}
\nu(X_{i})&=&\langle\xi_{X_{j}}X_{r+j},R(X_{j},X_{i})X_{r+j}\rangle
+
\langle\xi_{X_{r+j}}X_{j},R(X_{r+j},X_{i})X_{j}\rangle\\
 & & + \langle\xi_{Z}X_{i},R(Z,X_{i})X_{i}\rangle =
 \frac{3}{4}\langle\xi_{Z}X_{i},Z\rangle = 0,\\[0.5pc]
\nu(X_{r+i})&=&\langle\xi_{X_{j}}X_{r+j},R(X_{j},X_{r+i})X_{r+j}\rangle
+
\langle\xi_{X_{r+j}}X_{j},R(X_{r+j},X_{i})X_{j}\rangle\\
 & & + \langle\xi_{Z}X_{r+i},R(Z,X_{r+i})X_{r+i}\rangle =
- \frac{1}{4}\langle\xi_{Z}X_{r+i},Z\rangle = 0,\\[0.5pc]
\nu(Z)&=&\langle\xi_{X_{j}}X_{r+j},R(X_{j},Z)X_{r+j}\rangle +
\langle\xi_{X_{r+j}}X_{j},R(X_{r+j},Z)X_{j}\rangle = 0.
\end{array}
$$
Hence, $\nu$ is vanished and we have the result.
\end{proof}

Now, from Theorem \ref{harmclasscontact}, we can conclude
\begin{corollary} Any left-invariant almost contact metric structure
$(\varphi,\zeta =Z,\eta= \gamma,\langle\cdot ,\cdot \rangle)$ on
$H(1,r)$ of type ${\mathcal C}_{8}\oplus {\mathcal C}_{9}$
determines a harmonic map.
\end{corollary}
Using (\ref{intrinsicH}), one directly obtains $d^{*}\xi(Z) = 0$
and
$$
d^{*}\xi(X_{i}) =  -(\nabla_{Z}\xi_{Z}X_{i} +
\tfrac{1}{2}\xi_{Z}X_{r+i}), \;\;\;\; d^{*}\xi(X_{r+i}) =
\tfrac{1}{2}\xi_{Z}X_{i} - \nabla_{Z}\xi_{Z}X_{r+i}.
$$
Hence, we have
\begin{proposition}\label{pintrinsicH} On $H(1,r)$ the following conditions are equivalent:
\begin{enumerate}
\item[{\rm (i)}] The left-invariant almost contact metric structure
$(\varphi, \zeta = Z,\eta = \gamma,\langle \cdot ,\cdot\rangle)$
is harmonic; \item[{\rm (ii)}] its intrinsic torsion $\xi$
satisfies
\[
\langle\xi_{Z}X_{i},X_{j}\rangle =
\langle\xi_{Z}X_{r+i},X_{r+j}\rangle,\;\;\;
\langle\xi_{Z}X_{i},X_{r+j}\rangle = \langle\xi_{Z}X_{j},
X_{r+i}\rangle;
\]
\item[{\rm (iii)}]$\sum_{k=1}^{r}\{(\varphi_{j}^{r+k} +
\varphi_{r+j}^{k})(\varphi_{i}^{k} + \varphi_{r+i}^{r+k}) +
(\varphi_{r+j}^{r+k} -\varphi_{j}^{k})(\varphi_{i}^{r+k} -
\varphi_{r+i}^{k})\}  =  0,\\[0.5pc]
\sum_{k = 1}^{r}\{\varphi_{i}^{k}\varphi_{r+j}^{r+k} -
\varphi_{i}^{r+k}\varphi_{r+j}^{k} +
\varphi_{j}^{r+k}\varphi_{r+i}^{k} -
\varphi_{j}^{k}\varphi_{r+i}^{r+k}\}  =  0;$
\end{enumerate}
for all $i,j\in \{1,\dots ,r\}.$
\end{proposition}
\noindent Here $ds^{\rm ac} = 0$ and, for $X$ in the Lie algebra
of $H(1,r),$ we get
$$
\begin{array}{lcl}
(d^{*}\Ricac)(X) & = &
\sum_{k=1}^{2r}\Ricac(X_{k},\nabla_{X_{k}}X) =
\frac{\eta(X)}{2}\sum_{i=1}^{r}(\Ricac(X_{i},X_{r+i}) +
\Ricac(X_{r+i},X_{i}))\\[0.5pc]
 & = & \frac{\eta(X)}{8}
 \sum_{i,j = 1}^{r}(\varphi^{r+j}_{i} + \varphi^{j}_{r+i})(\varphi^{i}_{j} +
 \varphi^{r+i}_{r+j}).
\end{array}
$$
Then, from Propositions \ref{harmonicmap} and \ref{pintrinsicH},
we have
\begin{corollary} Left-invariant almost contact metric structures
$(\varphi, \zeta = Z,\eta = \gamma,\langle \cdot ,\cdot\rangle)$ on
$H(1,r)$ such that
$ \varphi_{i}^{r+j} = \pm \varphi^{r+i}_{j}$, $\; \varphi_{i}^{j} =
\pm \varphi_{r+i}^{r+j},
$ for all $i,j\in \{1,\dots ,r\},$ determine harmonic maps and
$d^{*}\Ricac = 0.$
\end{corollary}

\begin{remark} {\rm On $H(1,2)$ we consider the following examples of left-invariant almost contact metric
structures $(\varphi, \zeta = Z,\eta = \gamma,\langle \cdot
,\cdot\rangle)$ where
\begin{enumerate}
\item[{\rm (A)}] $\varphi X_{1} = \frac{\sqrt{2}}{2}(X_{2} +
X_{4}),\;\;\;\; \varphi X_{2} = \frac{\sqrt{2}}{2}(- X_{1} +
X_{3}),\;\;\;\; \varphi X_{3} = \frac{\sqrt{2}}{2}(-X_{2} + X_{4}),$\\
$\varphi X_{4} = -\frac{\sqrt{2}}{2}(X_{1} + X_{3}),$ \item[{\rm
(B)}] $\varphi X_{1} = \frac{\sqrt{2}}{2}(X_{2} + X_{4}), \;\;\;\;
\varphi X_{2} = -\frac{\sqrt{2}}{2}(X_{1} + X_{3}), \;\;\;\;
\varphi X_{3} = \frac{\sqrt{2}}{2}(X_{2} -
X_{4}),$\\
$\varphi X_{4} = \frac{\sqrt{2}}{2}(-X_{1} + X_{3}),$ \item[{\rm
(C)}] $\varphi X_{1} = \frac{\sqrt{2}}{2}(X_{2} + X_{3}),\;\;\;\;
\varphi X_{2} = \frac{\sqrt{2}}{2}(-X_{1} + X_{4}), \;\;\;\;
\varphi X_{3} = -\frac{\sqrt{2}}{2}(X_{1} +
X_{4}),$\\
$\varphi X_{4} = \frac{\sqrt{2}}{2}(-X_{2} + X_{3}).$
\end{enumerate}
Then, the example (A) is of type ${\mathcal C}_{8}\oplus {\mathcal
C}_{9} -({\mathcal C}_{8}\cup {\mathcal C}_{9})$ and (B) and (C)
are of type ${\mathcal C}_{8}\oplus {\mathcal C}_{9}\oplus
{\mathcal C}_{11} -({\mathcal C}_{8}\cup {\mathcal C}_{9}\cup
{\mathcal C}_{11}).$ Moreover, using Propositions
\ref{harmonicmap} and \ref{pintrinsicH}, these structures
determine harmonic maps if and only if
$\varphi_{1}^{3}(\varphi_{3}^{4} -\varphi_{1}^{2}) =
(\varphi_{1}^{3} +
\varphi_{2}^{4})(\varphi_{1}^{4}-\varphi_{2}^{3}) = 0.$ Hence, it
follows that the examples (A) and (B) determine harmonic maps and
(C) is not a harmonic almost contact structure.}
\end{remark}

\noindent {\bf The generalised Heisenberg groups H(p,1).} The Lie
group $H(p,1)$ consists of matrices
\[
a = \left (
\begin{array}{ccc}
1 & A & c\\
0 & I_{p} & B^{t}\\
 0 & 0 & 1
\end{array}
\right ),
\]
where $I_{p}$ denotes the identity matrix of type $p\times p,$ $A =
(a_{1},\dots ,a_{p})\in {\mathbb{R}}^{p},$ $B = (b_{1},\dots
,b_{p})\in {\mathbb{R}}^{p}$ and $c\in \mathbb{R}.$ Moreover,
$H(p,1)$ is a connected, simply connected,  nilpotent,  of dimension
$2p+1$ and its center is one-dimensional.  A global system of
coordinates $(x^{i},x^{p+i},z),$ $1\leq i\leq p,$ on $H(p,1)$ is
defined by
\[
x^{i}(a) = a_{i},\;\;\; x^{p+i}(a) = b_{i},\;\;\; z(a) = c
\]
and a basis of left-invariant one-forms is given by $\alpha_{i} =
dx^{i},$ $\alpha_{p+i} = dx^{p+i},$ $\gamma = dz -
\sum_{j=1}^{p}x^{j}dx^{p+j}.$  Therefore, for the dual
left-invariant vector fields, we  have
\[
X_{i} = \tfrac{\partial}{\partial x^{i}},\;\;\; X_{p+i} =
\tfrac{\partial}{\partial x^{p+i}} + x^{i}\tfrac{\partial}{\partial
z} ,\;\;\; Z = \tfrac{\partial}{\partial z}.
\]
On $H(p,1)$ we consider the Riemannian metric $\langle \cdot,
\cdot \rangle$ for which these vector fields form an orthonormal
basis. Then, the corresponding Levi Civita connection is
determined by
\begin{equation}\label{nablaHp}
\begin{array}{lclcl}
\nabla_{X_{i}}X_{p+i} & = & - \nabla_{X_{p+i}}X_{i} & = &
\frac{1}{2}Z,\\[1mm]
 \nabla_{X_{i}}Z & = &\nabla_{Z}X_{i} & = &
-\frac{1}{2}X_{p+i},\\[1mm]
\nabla_{X_{p+i}}Z & = & \nabla_{Z}X_{p+i} & = & \frac{1}{2}X_{i},
\end{array}
\end{equation}
where the remaining covariant derivatives with respect to this basis
vanish. For the nonvanishing components of the curvature tensor $R$,
we obtain
\begin{equation}\label{RHp}
\begin{array}{lclclcl}
R(X_{i},X_{j},X_{p+i},X_{p+j}) & = & \frac{1}{4},\;\; i\neq j, & &
R(X_{i},X_{p+i},X_{i},X_{p+i}) & = & \frac{3}{4},\\[0.5pc]
R(X_{i},X_{p+i},X_{j},X_{p+j}) & = & \frac{1}{2},\;\; i\neq j, & &
R(X_{i},X_{p+j},X_{j},X_{p+i}) & = & \frac{1}{4},\;\; i\neq
j,\\[0.5pc]
R(X_{i},Z,X_{i},Z) & = & -\frac{1}{4}, & & R(X_{p+i},Z,X_{p+i},Z)
& = & - \frac{1}{4}.
\end{array}
\end{equation}

As on $H(1,r),$ we consider left-invariant almost contact
structures $(\varphi,\zeta,\eta)$ on $H(p,1)$ compatible with
$\langle \cdot, \cdot \rangle$ with $\zeta = Z$ and $\eta =
\gamma.$ Denote by $\varphi^{l}_{k},$ $k,l = 1,\dots ,2p,$ the
(constant) components of $\varphi$ with respect to the basis
$\{X_{k},Z\}.$ Then, we have
\begin{lemma} Any left-invariant almost contact metric structure
$(\varphi,\zeta= Z,\eta= \gamma, \langle \cdot, \cdot \rangle)$ on
$H(p,1)$ is of type $({\mathcal C}_{6}\oplus {\mathcal
C}_{7}\oplus {\mathcal C}_{11})-{\mathcal C}_{11}.$ Moreover, it
is of type
\begin{enumerate}
\item[{\rm (i)}] ${\mathcal C}_{6}\oplus {\mathcal C}_{7}$ if and
only if $\varphi^{p+j}_{i} = \varphi_{j}^{p+i}$ and
$\varphi^{p+j}_{p+i} = \varphi^{j}_{i},$ for all $i,j\in \{1,\dots
,p\};$ \item[{\rm (ii)}] ${\mathcal C}_{6}$ if and only if
$\varphi^{p+i}_{i} = \lambda,$ with $\lambda= \pm 1,$ for all
$i\in \{1,\dots ,p\},$ being zero the remainder components of
$\varphi.$ Then, it is $\pm\frac{1}{2}$-Sasakian; \item[{\rm
(iii)}] ${\mathcal C}_{7}$ if and only if
$\sum_{i=1}^{p}\varphi^{p+i}_{i}= 0,$ $(p\geq 2).$
\end{enumerate}
\end{lemma}
\begin{proof} From (\ref{nablaHp}), we obtain
\[
(\nabla_{X_{i}}\varphi)X_{k}  =
\tfrac{1}{2}\varphi^{p+i}_{k}Z,\;\;\;\;\;\;
(\nabla_{X_{p+i}}\varphi)X_{k} = -\tfrac{1}{2}\varphi^{i}_{k}Z,
\]
where $1\leq i\leq p$ and $1\leq k\leq 2p,$ and
$$
\begin{array}{lcl}
(\nabla_{Z}\varphi)X_{i} & = & \frac{1}{2}
\sum_{j=1}^{p}\{(\varphi_{i}^{p+j} + \varphi_{p+i}^{j})X_{j} +
(\varphi_{p+i}^{p+j} - \varphi_{i}^{j})X_{p+j}\},\\[0.5pc]
(\nabla_{Z}\varphi)X_{p+i} & = & \frac{1}{2}
\sum_{j=1}^{p}\{(\varphi_{p+i}^{p+j} - \varphi_{i}^{j})X_{j} -
(\varphi_{p+i}^{j} + \varphi_{i}^{p+j})X_{p+j}\}.
\end{array}
$$
Then, $d^{*}F(Z) = \sum_{i=1}^{p}\varphi_{i}^{p+i},$ $d^{*}\eta
=0$ and $(\nabla_{Z}\varphi) = 0$ if and only if
$\varphi^{p+j}_{i} = \varphi_{j}^{p+i}$ and $\varphi_{p+i}^{p+j} =
\varphi_{i}^{j}.$ Moreover, these last conditions imply
$$
\begin{array}{l}
(\nabla_{X_{i}}\varphi)X_{j} = (\nabla_{X_{j}}\varphi)X_{i} =
(\nabla_{X_{p+i}}\varphi)X_{p+j} =
(\nabla_{X_{p+j}}\varphi)X_{p+i}\\
 \qquad \quad  \quad \; \,\,\, = (\nabla_{\varphi
X_{i}}\varphi)\varphi X_{j} = (\nabla_{\varphi
X_{p+i}}\varphi)\varphi X_{p+j} =
\frac{1}{2}\varphi^{p+j}_{i}Z;\\[1mm]
(\nabla_{X_{i}}\varphi)X_{p+j} = (\nabla_{X_{p+j}}\varphi)X_{i}
=(\nabla_{\varphi X_{i}}\varphi) \varphi X_{p+j} =
(\nabla_{\varphi X_{p+j}}\varphi)\varphi X_{i} = -
\frac{1}{2}\varphi_{i}^{j}Z.
\end{array}
$$
This proves the Lemma.
\end{proof}
From (\ref{nablaHp}), the intrinsic torsion $\xi$ of
$(\varphi,\zeta = Z,\eta= \gamma,\langle\cdot,\cdot \rangle)$
satisfies
\begin{equation}\label{intrinsicHp}
\begin{array}{lclclcl}
\xi_{X_{i}}Z & = & \frac{1}{2}X_{p+i}, & & \xi_{X_{p+i}}Z & = &
-\frac{1}{2}X_{i},\\
\xi_{X_{i}}X_{j} & = & \xi_{X_{p+i}}X_{p+j} =0,
& & \xi_{X_{i}}X_{p+j} & = & \xi_{X_{p+i}}X_{j} =0,\;i\neq j.
\end{array}
\end{equation}
Then, as in the case $H(1,r),$ we obtain that the one-form $\nu$
defined in \eqref{harmmap1} vanishes and we have
\begin{proposition} Any harmonic left-invariant almost contact metric structure
$(\varphi,\zeta =Z,\eta= \gamma,\langle\cdot ,\cdot \rangle)$ on
$H(p,1)$ determines a harmonic map.
\end{proposition}
Moreover, from Theorem \ref{harmclasscontact}, we can conclude
\begin{corollary} Any left-invariant almost contact metric structure
$(\varphi,\zeta =Z,\eta= \gamma,\langle\cdot ,\cdot \rangle)$ on
$H(p,1)$ of type ${\mathcal C}_{6}\oplus {\mathcal C}_{7}$
determines a harmonic map.
\end{corollary}

\end{document}